\documentclass[11 pt]{amsart}
\usepackage{amsmath,amsthm,amssymb,xypic}

\newtheorem{Theorem}{Theorem}[section]
\newtheorem{Lemma}[Theorem]{Lemma}
\newtheorem{Proposition}[Theorem]{Proposition}
\newtheorem{Corollary}[Theorem]{Corollary}
\newtheorem{PD}[Theorem]{Proposition--Definition}
\newtheorem{exercise}[Theorem]{Exercise}
\theoremstyle{definition}
\newtheorem{Definition}[Theorem]{Definition}
\newtheorem{CO}[Theorem]{Crucial observation}
\theoremstyle{remark}
\newtheorem{remark}[Theorem]{Remark}
\newtheorem{claim}{Claim}
\newtheorem{conjecture}[Theorem]{Conjecture}
\newtheorem{example}[Theorem]{Example}

\numberwithin{equation}{section}

\DeclareMathOperator{\mult}{mult}
\newcommand{\hint}{ {\sf hint}\rm:  }
\newcommand{\NE}{\overline{NE}}
\newcommand{\Ol}{{\mathcal O}}
\renewcommand{\O}{{\mathcal O}}
\newcommand{\Os}{{\mathcal O}}
\newcommand{\D}{{\mathcal D}}
\newcommand{\F}{{\mathbb F}}
\newcommand{\E}{{\mathcal E}}
\renewcommand{\H}{{\mathcal H}}
\newcommand{\I}{{\mathcal I}}

\newcommand{\Proj}{{\mathbb P}}
\renewcommand{\P}{{\mathbb P}}
\newcommand{\J}{{\mathcal J}}
\newcommand{\Q}{{\mathbb  Q}}
\newcommand{\N}{{\mathbb  N}}
\newcommand{\Z}{{\mathbb  Z}}
\newcommand{\C}{{\mathbb  C}}
\newcommand{\G}{{\mathbb G}}
\newcommand{\R}{{\mathbb  R}}
\newcommand{\Sn}{{\bf S}}
\newcommand{\Sc}{{\bf S}}
\newcommand{\T}{T^{\#}}

\newcommand{\Hs}{H^{\#}}

\newcommand{\X}{X^{\#}}
\newcommand{\h}{{\mathcal H}}
\newcommand{\hh}{{\mathcal H}^{\#}}

\newcommand{\f}{\varphi}
\newcommand{\ra}{\rightarrow}
\newcommand{\raa}{\longrightarrow}
\newcommand{\iso}{\simeq}
\newcommand{\flip}{\dasharrow}

\newcommand{\lel}[1]{\sim_{#1}}
\newcommand{\nel}[1]{\equiv_{#1}}

\begin{document}

\title[Morphisms of Projective Varieties]{Morphisms of Projective Varieties
from the viewpoint of Minimal Model Theory}

\author{Marco Andreatta}
\address{Dipartimento di
Matematica,Universit\`a di Trento, 38050 Povo (TN), Italia}
\email{andreatt@science.unitn.it}

\author{Massimiliano Mella}
\address{Dipartimento di
Matematica, Universit\`a di Ferrara, Via Machiavelli 35,
  44100 Ferrara, Italia}
\email{mll@unife.it}

\maketitle
\tableofcontents
\part*{Introduction}
One of the main results of last decades algebraic geometry was the foundation
of Minimal Model Program, or MMP, and its proof in dimension three.
Minimal Model Theory shed a new light on what is
  nowadays called higher dimensional geometry. In mathematics high numbers
are really a matter of circumstances and here we mean greater than or equal to 3.
The impact of MMP has been felt in almost all areas of algebraic geometry.
In particular the philosophy and some of the main new objects like
extremal rays,  Fano-Mori contractions or
spaces and log varieties started to play around and give fruitful answer to
different problems.

The aim of Minimal Model Program is to choose, inside of a birational class
of varieties, ``simple'' objects. The first main breakthrough of the theory
is the definition of these objects: minimal models and Mori spaces.
This is related to numerical properties of the intersection
of the canonical class of a variety with effective cycles. After this,
old objects, like the Kleiman cone of effective curves and rational curves
on varieties,
acquire a new significance. New ones, like Fano-Mori contractions,
start to play an important role.
And the tools developed to tackle these problems allow the study of
formerly untouchable varieties.

Riemann surfaces were classified, in the $\mbox{{\sc XIX}}^{\rm th}$ 
century, according to the curvature
of an holomorphic metric. Or, in other words, according to the 
Kodaira dimension.
Surfaces needed a harder amount of work. For the first time 
birational modifications
played an important role. The theory of (-1)-curves studied by
the Italian school of  Castelnuovo, Enriques and Severi, at the beginning of
$\mbox{{\sc XX}}^{\rm th}$ century, allowed to define minimal surfaces.
Then the first rough classification of the latter, again by Kodaira 
dimension, was fulfilled.
Minimal Model Program is now a tool to start investigate this 
question in dimension
3 or higher.

In writing these notes we want to give our point of view on this area
of research. We are not trying to give a treatment of the whole 
subject.  Very nice books appeared
recently for this purpose, and we often refer to them in the paper.
We would like to present, in a sufficiently self contained
way, our contributions and interests in this field of mathematics.
We will study Fano-Mori spaces both from the biregular and the 
birational point of view.
For the former we will recall and develop Kawamata's Base Point Free
technique and some of Mori's deformation arguments.
For the latter we lean on Sarkisov and \#-Minimal Model Programs.

The content of the parts is the following.
In Part 1 we collect the main definitions and Theorems we are going 
to use afterwords.

Part 2 is studied to get the reader acquainted with the Base Point 
Free technique,
BPF. For this purpose we give, or sketch, proofs of
Kawamata's Base Point Free Theorem, trying to hide the technicalities.

In Part 3 we introduce the main actor of the book, Fano-Mori 
contractions and more
generally Fano-Mori spaces. Using BPF we describe then
the main properties that will allow us to study them.

Part 4 is the applications of all the above to smooth varieties.
  Namely we give a biregular classification of Fano-Mori spaces of 
dimension less than or equal to four
and Mukai manifolds.

In Part 5 we present an other side of the moon, the birational world. 
Here a beautiful old Theorem
of Castelnuovo and Noether is proved in a modern language. Philosophy 
and applications of
Minimal Model Program for threefolds are outlined.

These notes collect some topics we presented in three mini-courses
  which were held in Wykno (Pl) (1999), Recife (Br) (2000)
and Ferrara (It) (2000), respectively.
We discussed with many people about this subject and we are very thankful
to them all. But we would like to distinguish Jaros\l aw Wi\'sniewski 
and thank him deeply.

\part{Preliminaries}

In this part we collect all definitions which are more or less
standard in the algebraic
geometry realm in which we live.

\section{The Kleiman--Mori cone of a projective variety.}

First we fix a good category of objects (real differentiable
varieties are not the good ones to extend the Riemann and Poincar\'e approach).
Let $X$ be a normal variety over an algebraically closed field $k$
of dimension $n$,
that is an integral separated scheme which is of finite type over $k$.
We actually assume also that $char(k) =0$,
nevertheless many results at the beginning of the theory hold also
in the case of positive characteristic.

We have to introduce some basic objects on $X$.

Let $Div(X)$ be the group of Cartier divisors on $X$
and $Pic(X)$ be the group of line bundles on $X$.
Let also $Z^1(X)$ be the group of Weil divisors
and $Z_1(X)$ be the group of 1-cycles on $X$
i.e. the free abelian group generated, respectively, by prime divisors, and
  reduced irreducible curves.

We will often use $\Q$-Cartier divisors, that is linear
combinations with rational coefficients of Cartier divisors. For these
objects it
is useful to introduce the following notations.
Let $D=\sum d_i D_i \in Div(X)\otimes \Q$ be a $\Q$-Cartier divisor, 
then\label{pg:nota}
$\lfloor D\rfloor:=\sum \lfloor d_i\rfloor D_i$, $\lceil D\rceil:=-\lfloor
-D\rfloor$ and
$\langle D\rangle:=D-\lfloor D\rfloor$, where $\lfloor d_i\rfloor$ is the
integral part of $d_i$.

Then there is a pairing
$$Pic(X) \times Z_1(X) \ra \Z$$
defined, for an irreducible reduced curve $C\subset X$, by $(L, C) 
\ra L\cdot C := deg_C(L_{|C})$,
and extended by linearity.

Two invertible sheaves $L_1, L_2 \in Pic(X)$ are {\sl numerically equivalent},
denoted by $L_1 \equiv L_2$, if $L_1\cdot C = L_2\cdot C$ for every 
curve $C \subset X$.
Similarly, two $1$-cycles $C_1, C_2$ are {\sl numerically equivalent}, $C_1
\equiv
C_2$ if $L\cdot C_1 = L\cdot C_2$ for every $L\in Pic(X)$.

Define
$$N^1 X = (Pic(X)/\equiv)\otimes \R   {\hbox { and }}
N_1 X = (Z_1(X)/\equiv)\otimes \R ;$$

obviously, by definition, $N^1(X)$ and $N_1(X) $ are dual $\R$-vector spaces
and $\equiv$ is the smallest equivalence relation for which this holds.

In particular for any divisor $H \in Pic(X)$ we can view the class of $H$ in
$N^1(X)$ as a linear form on $N_1(X)$. We will use the following notation:
$$H _{\geq 0}: = \{ x \in N_1(X) : H\cdot x \geq 0\} {\hbox{   and 
similarly for }}
>0, \leq 0, <0$$
and
$$H^{\perp}: = \{ x \in N_1(X) : H\cdot x =0\}.$$

The fact that $\rho := dim_{\R} N^1(X)$ is finite is the Neron-Severi 
theorem, \cite[pg 461]{GH}.
The natural number $\rho$ is called the Picard number of the variety $X$.
(Note that for a variety defined over $\C$ the finite dimensionality of
$N_1(X)$ can be read from the fact that $N_1(X)$ is a subspace of $H_2(X,
\R)$).

More generally, given a projective morphism $f:X\ra Y$ and $A,B\in
Div(X)\otimes \Q$,
then $A$ is  {\sl $f$-numerically equivalent} to $B$
( {\sf $A\nel{f} B$} )
if $A\cdot C=B\cdot C$ for any curve $C$ contracted by $f$.
$A$ is {\sl $f$-linearly} equivalent to $B$
( {\sl $A\lel{f} B$} )
if $A-B\sim f^*M$, for some line bundle $M\in Pic(Y)$, we will
suppress the subscript when no confusion is likely to arise.

\smallskip
Note that if $X$ is a surface then $N^1(X) = N_1(X)$ ;
using M. Reid words (see \cite{Reid}), ''Although very simple,
this is one of the key ideas of Mori theory, and came as a surprise to anyone
who knew the theory of surfaces before 1980: the
quadratic intersection form of the curves on a nonsingular surface can
for most purpose be replaced by the bilinear pairing between $N^1$ and
$N_1$, and in this form
generalizes to singular varieties and to higher dimension.''

\smallskip
We notice also that {\sl algebraic equivalence},see \cite[pg 
461]{GH}, of $1$-cycles
implies numerical equivalence. Moreover, if $X$ is a variety over $\C$
then, in terms of Hodge Theory, $N^1(X) = (H^2(X,\Z)/(Tors) \cap 
H^{1,1}(X))\otimes \R$.

\smallskip
We denote by $NE(X) \subset N_1(X)$
  the cone of effective $1$-cycles, that is
$$NE(X) = \{ C \in N_1(X) : C = \sum r_iC_i \mbox{\rm   where  } r_i \in \R,
r_i \geq 0\},$$
where $C_i$ are irreducible curves.
Let $\overline{NE(X)} $ be the closure of $NE(X)$ in the
real topology of $N_1(X)$. This is called the {\sl Kleiman--Mori cone}.

We also use the following notation:
$$\overline{NE(X)}_{H \geq 0} : = \overline{NE(X)} \cap H_ { \geq 0}
{\hbox{   and similarly for }} >0, \leq 0, <0.$$

One effect of taking the closure is the following trivial observation,
which has many important use in applications: if $H \in N^1(X) $
is positive on $\overline{NE(X)} \setminus {0}$ then the section
$(H\cdot z = 1) \cap \overline{NE(X)}$ is compact. Indeed, the projectivised
of the closed cone $\overline{NE(X)}$ is a closed subset of $\mathbb
P^{\rho -1} = P(N_1(X))$,
and therefore compact, and the section $ (H\cdot z = 1)$ projects
homeomorphically to it.
The same holds for any face or closed sub-cone of $\overline{NE(X)}$.

\medskip
An element $H \in N^1(X)$ is called {\sl numerically
eventually free} or {\sl numerically effective}, for short {\sl nef}, if
$H\cdot C \geq 0$ for every curve $C\subset X$
(in other words if
$H \geq 0$  on $\overline {NE(X)}$).

The relation between nef and ample divisors
is the content of the following Kleiman criterion that
is a corner stone of Mori theory.

\begin{Theorem} [\cite {Kle}]
For $H \in Pic(X)$, view the class of $H$ in $N^1(X)$ as a linear form
on $N_1(X)$. Then
$$ H {\hbox { is ample} } \Longleftrightarrow
HC > 0 {\hbox { for all } } C \in \overline{NE(X)}\setminus \{0\}.$$
\end{Theorem}

In other words the theorem says that the cone of ample divisors is the interior
of the nef cone in $N^1(X)$, that is the cone spanned by all nef divisors.

Note that it is not true that $HC > 0$  for every curve $C \subset X$
implies that $H$
is ample, see for instance \cite[Example 4.6.1]{CKM}. The condition 
in the theorem is stronger.

This is only a weak form of Kleiman's criterion, since $X$ is a priori
assumed to be projective. The full strength of Kleiman's criterion
gives a necessary and sufficient condition for ampleness in terms of the
geometry
of $\overline{NE(X)}$.

\bigskip
Assume that $X$ is smooth and  denote by $K_X$ the {\sl canonical divisor}
of $X$, that is
an element of $Div(X)$ such that $\Ol_{X}(K_X) = \Omega^n_{X}$, where 
$\Omega_X$ is the sheaf of one forms
on $X$.

The first main theorem of Mori theory is the following description of the
negative part, with respect to $K_X$, of the Kleiman-Mori cone.

We recall that, by definition, a {\sl rational curve} is an
irreducible, reduced curve defined over $k$ whose normalization is $\mathbb
P^1$.

\begin{Theorem} [\cite{Mo}, cone theorem] Let $X$ be a non singular 
projective variety.

1) There are countably many rational curves $C_i \subset X$ such
that $0 < -C_i K_X \leq dimX +1$ and
$$\overline{NE(X)} = \overline{NE(X)}_{{K_X} \geq 0} + \sum \R_{\geq 0}[C_i].$$

2) For any $\epsilon > 0$ and ample divisor $H$,
$$\overline{NE(X)} = \overline{NE(X)}_{{K_X+\epsilon H} \geq 0} + \sum
_{\hbox{finite}}
\R_{\geq 0}[C_i].$$
\label{cone}
\end{Theorem}

In simple words the theorem says the following.
Consider the
linear form on $N_1(X)$ defined by $K_X$; the part of the Kleiman-Mori
cone $\overline{NE(X)}$
which sits in the negative semi-space defined by $K_X$ (if not
empty) is locally
polyhedral and it is spanned by a countable number of {\sl extremal rays},
$\R_{\geq 0}[C_i]$.
Moreover each extremal ray is spanned in $N_1(X)$ by a rational
curve with bounded intersection with the linear form $ - K_X$ and
moving an $\epsilon$ away from the hyperplane $K_X= 0 $ (in
the negative direction)
the number of extremal rays becomes finite.

There are essentially two ways of proving this theorem; the original one,
which is due to Mori, is very geometric and valid in any characteristic.
It is presented in the paper \cite{Mo} and in many other places, for example
in \cite{KoMo} and \cite{De}. It is based on the study of 
deformations of rational curve on an
algebraic variety, it makes use of the theory of Hilbert schemes
and of theorems like  for instance \ref{hilb}.

Another proof was provided by Y. Kawamata  (\cite{Ka}); it gives the cone
theorem as a consequence of the following rationality theorem.

\begin{Theorem} [\cite{KMM} 4.1.1, rationality theorem]
Let $X$ be a $n$-dimensional variety
defined over $\C$ which is smooth or, more generally
with LT singularities (see Definition \ref{lc}),
for which $K_X$ is not nef.

Let $L$ be an
ample line bundle on $X$ and define the
nef value (or nef threshold) of the pair $(X,L)$ by
$$r=\mbox{\rm inf}\{t\in \R: K_X+tL \textrm{ is nef} \}.$$

Then the nef value is a
rational number.

Moreover if
$a: = min\{ e \in \N: eK_X \textrm{ is Cartier}\}$,
and $ar :=v/u$, with $(v,u)= 1$, then $v \leq a(n+1)$.

\label{Ka-rat}
\end{Theorem}

The proof of this proposition uses the technique of the base point free theorem
which we will introduce in the next section. In particular
it makes use of vanishing theorems and it is therefore valid
only in characteristic zero.

It was noticed by M. Reid and Y. Kawamata that the rationality theorem
and the Base point free theorem implies immediately the Mori's cone theorem,
in the more general case of varieties with LT singularities.

\smallskip
A very nice presentation of the above theorems (Kleiman-Mori-Kawamata),
together with complete proofs, in the {\bf case of surfaces}
is in \cite{Reid}, Chapter D.
This material can be presented in few hours (3-4) to an audience
with a limited knowledge of basic algebraic geometry and it can provide
a good insight in the field; this is our experience at the Ferrara's course.

The surface case is a perfect tutorial case
in order to understand the Minimal Model Program.
This was first pointed out by S. Mori who
worked out a complete description of extremal rays
in the case of a smooth surface (see \cite[Chapter 2]{Mo},
see also \cite[pg 21-23, \S 1.4]{KoMo}).
Moreover he also showed how it is possible to associate to
each extremal ray a morphism from the surface. When the ray is 
spanned by a rational curve
with self intersection $-1$, this is a celebrated theorem
of Castelnuovo. Castelnuovo's proof is also very enlightening and
it can be found in \cite[th\'eor\'em II.17]{Be},  or in \cite[theorem 
V. 5.7]{Ha}.

\section{Fujita $\Delta-genus$}
\label{sec:fuj}
A classical approach to the classification of projective varieties,
which dates back to the Italian school, consists of the following:
a) take an hyperplane section, b) characterise it by induction,
c) describe the original variety by ascending the properties of the hyperplane
section. To stress its
classical flavor T. Fujita called it {\sl Apollonius method}; we will 
now introduce some
definitions and techniques as presented in the work of T. Fujita (\cite{Fu});
see also section \ref{chap:ladder}.

\begin{Definition} Let $F$ be a variety of dimension $d$ and
let $L$ be an ample line bundle on $F$. The pair $(F,L)$ is called a
polarized variety.
We will denote by
$$\chi(F,tL) = \Sigma \chi_j\frac{(t(t+1)...(t+j-1)}{j!}$$
the Hilbert polynomial
of $(F,L)$.

Then the $\chi_j$'s are integers and  $\delta (F,L) := \chi_n = 
L^d>0$ is called
the degree of $(F,L)$, while $g(F,L) := 1-\chi_{(n-1)}$ is called the
sectional genus.

The $\Delta$-genus of $(F,L)$ is defined by the formula
$$\Delta (F,L) = d+\delta -h^0(F,L).$$
\end{Definition}

\begin{Definition} Let $(F,L)$ be a polarized variety.
Let $D$ be a member of $|L|$ and suppose that
$D$, as a subscheme of $F$, is irreducible and reduced.
In such a case $D$ is called a rung of $(F,L)$.
Let $r:H^0(F,L) \ra H^0(D,L_D)$ be the restriction map. If
$r$ is surjective the rung is said to be a regular rung.

A sequence $F = F_d \supset F_{d-1} \supset....\supset F_1$
of subvarieties of $F$ such that $F_i$ is a rung (a regular rung)
of $(F_{i+1}, L_{i+1})$ is called a ladder (a regular ladder).
\end{Definition}

\begin{remark} If $D$ is a rung then the pair $(D,L_D)$ is a 
polarized variety of
dimension $d-1$. The structure of $(F,L)$ is reflected in that
of $(D,L_D)$. One can study $(F,L)$ via $(D,L_D)$ using induction on $d$.
This is the main idea of the Apollonius method.
In particular we have that $\chi(D,tL_D) = \chi(F,tL)-\chi(F,(t-1)L)$,
$g(D,L_D) = g(F,L)$, $\delta(D,L_D) = \delta(F,L)$
and $\Delta(F,L) -\Delta(D,L_D) = dim Coker(r)$.
In particular if the rung is regular the two $\Delta$-genera are the same.

In classical geometry the number $dimCoker(r)$ was called deficiency.

\medskip
Assume that $L$ is very ample and let $\f_L$ be the map associated to the
elements of the complete linear system $|L|$.
Then it is a classical result that
$\Delta(F,L) \geq 0$ and equality holds for the so called ''Varieties
of Minimal Degree'', \cite[pg 173]{GH}.
  This varieties are classified as projective spaces, hyperquadrics,
scroll over rational normal curves or generalised cones over them.

In the case of surfaces a precise statement is the following:
\begin{Proposition}
Let $(S,L)$ be a pair with $S$ a surface and $L$ an
ample line bundle on $S$. If $\Delta(S,L) = 0$
then the pair is  among the following:
\item{(1)} $(\Proj^2,\O(e))$, with $e = 1,2$,
\item{(2)} $(\F_r, C_0 + kf)$ with $k \geq r+1$, $r\geq 0$,
\item{(3)} $(\Sc_r,\O_{\Sc_r}(1))$ with $r\geq 2$.

\noindent
(Where
$\F_r$ is a Hirzebruch surface, i.e. a $\Proj^1$-bundle $\Proj(\O(r)\oplus\O)$
over projective line $\Proj^1$ with a unique section $C_0\subset \F_r$
(isomorphic to $\Proj^1$) such that $C_0^2=-r\leq 0$ and a fiber of 
the projection
$\F_r\ra\Proj^1$ which we will denoted by $f$.
While $\Sc_r$ is a (normal) cone defined by contracting
$C_0\subset \F_r$ to a normal point; in terms of projective geometry 
$\Sc_r$ is a cone over
$\Proj^1\hookrightarrow\Proj^r$ embedded via Veronese map ($r$-uple embedding).
The restriction of the hyperplane section line bundle
from $\Proj^{r+1}$ to $\Sc_r$ will be denoted by $\O_{\Sc_r}(1)$)
\label{surface}
\end{Proposition}

\begin{exercise}
Prove the classification of surfaces 
of Minimal degree .
The first step consists in showing that if a
line meets a surface of minimal degree in three
or more points then it lies on the surface (see it for instance at \cite[pg 525]{GH}).
\end{exercise}

\bigskip
A modern approach to this classification which extends to the case when $L$
is merely ample is due to
T. Fujita (\cite{Fu}).
\end{remark}

\begin{Proposition} Let $(F,L)$ be a polarized variety and assume that there
exists a ladder for this pair. Then $\Delta(F,L) \geq 0$
(this is actually always true, without the assumption of the existence of a
ladder).
If moreover the ladder is regular and for a divisor
$D_1 \in |L_{|F_1}|$ the map $H^0(F_1,L_{|F_1}) \ra 
H^0(D_1,L_{|D_1})= \C^{\delta}$
is surjective (we will call this a complete regular ladder)
then $\Delta(F,L) = 0$ and the pair $(F,L)$ is a
variety of minimal degree; in particular $F$ is normal,
$g(F,L) = 0$ and $L$ is very ample.
\end{Proposition}

\begin{proof} The proof follows immediately from the above observations plus
the fact that the surjectivity of  $H^0(F_1,L_1) \ra
H^0(D_1,L_{D_1})= \C^{\delta}$
implies that $D_1$ is a rational normal curve.
\end{proof}

\begin{exercise}Let $F_1$ be a curve and $L$ a line bundle on
$F_1$. Assume that  $D\in |L|$ is  an effective divisor such that
$H^0(F_1,L) \ra H^0(D,L_D)= \C^{\delta}$ is surjective.
Prove that $F_1$ is a rational normal curve, when embedded by $|L|$.
\end{exercise}

\part{Base point free technique}
\setcounter{section}{0}
\label{ch:BPF}
In this part we introduce the Base point free technique (for short
BPF). This theory has
been mainly developed by
Kawamata, Reid, Shokurov in a series of papers, see \cite{KMM},
\cite{Ko0}
and \cite{Ka1}.
The aim of BPF is to show that an adjoint linear system, under some conditions,
is free from fixed points.
In the first section we will try to prevent the reader from a too technical
approach,
giving the main ideas and results, without too many definitions and details.
The latter are left for the interested reader, together with examples
and exercises.

\section{Base point freeness}
We start with the easy case of the curve:
let $C$ be a compact Riemann surface of genus $g$ and let $K_C$ be the
canonical
bundle of $C$. To give a morphism $C \ra \Proj ^N$ is equivalent to give
a line bundle $H$ without base points.
For this we have the very well known

\begin{Theorem} If $degH \geq 2g$ then $H$ has no base point.
\end{Theorem}

\begin{proof}
Let $L:= H-K_C$ and let $x\in C$ be a point on $C$.
Note that by assumption $deg L \geq 2$ and thus
\begin{equation}
\label{curvevan}
H^1(C,K_C+L-x)=H^0(C,x-L)=0,
\end{equation}
the first equality coming from Serre duality.

Then we consider the exact sequence
$$0\ra\O_C(K_C+L-x)\ra \O_C(K_C+L)\ra \O_x(K_C+L)\ra 0,$$
which comes by tensoring the structure sequence of $x$ on $C$,
$$0\ra \I_x\ra \O_C \ra \O_x\ra 0,$$
by the line bundle $K_C+L$.

The sequence gives rise to a long exact sequence in cohomology whose
first terms are, keep in mind equation (\ref{curvevan}),
$$0\ra H^0(C,K_C+L-x)\ra H^0(C,K_C+L)\ra H^0(x,K_C+L)\ra 0.$$

In particular we have the surjective map
$$H^0(C,H)\stackrel{\alpha}{\ra} H^0(x,H)\ra 0.$$

Furthermore
$x$ is a closed point and therefore
\begin{equation}
\label{curvenonvan}
H^0(x,H)=\C\neq 0.
\end{equation}

The surjectivity of $\alpha$ translates into the existence of
a section of $\O(H)$ which is not vanishing at $x$. That is the pull back
via $\alpha$ of $1$.
\end{proof}

What we have done can be summarized in the following slogan, which is
somehow the manifesto of the base point free technique.

{\it Construct
section of an adjoint line bundle proving a vanishing statement,
(\ref{curvevan}), and
a non vanishing on a smaller dimensional variety, (\ref{curvenonvan}).}

What happens if we try to generalise this to higher dimensional varieties
and which problems shall we encounter ?

\noindent{\it Simple observation}:
  the point $x\in C$ is a smooth Cartier divisor, that is
why with an abuse of language we wrote $H^i(C,K_C+L-x)$. This is no more
true for a point on a variety $X$ of higher dimension.

Let $x \in X$ be  a point of a smooth projective variety $X$
of dimension $n$. If we just want to mimic the above arguments, then
in (\ref{curvevan})
we are concerned with cohomology groups of non locally free sheaves which
are of difficult interpretation.
Note that there is a way to
make a divisor out of a point: blow it up! Do it and get a
morphism $\pi:Y\ra X$ with exceptional divisor $E$ and
$$\begin{array}{l}
\pi^*K_X= K_Y -(n-1)E\\
\end{array}
$$
Assume now that we want to prove that $x$ is not a base
point of a divisor of the type $H := K_X+ L$; we pull back
the divisors on $Y$ and we have an exact sequence, coming from
the structure sequence associated to $E$, of the type
$$ H^0(Y,\pi^*(K_X+L)) \ra H^0(E,\pi^*(K_X+L))\ra H^1(Y,\pi^*(K_X+L)-E).$$

Since $H^0(Y,\pi^*(K_X+L)) = H^0(X,K_X+L)$ (Hartogs theorem)
and $H^0(E,\pi^*(K_X+L)) = \C$, we have to prove ''only'' the vanishing
of
$$H^1(Y,\pi^*(K_X+L) -E) = H^1(Y,K_Y+\pi^*L-nE).$$

This is of course in general not true and one has to choose
carefully good assumptions on $L$ to have a
vanishing theorem of this type; let us state the best
available version of it (apart generalising it to a relative
or to a singular situation) which is due to Kawamata-Viehweg
(they worked on previous versions of Enriques, Kodaira, Ramanujan, ...)

\begin{Theorem}[Vanishing Theorem,
see \cite{KMM} or \cite{EV}]\

Let $X$ be a smooth variety and let $D = \sum a_iD_i$ be
a $\Q$-Cartier divisor satisfying the following conditions:
\item{i)} $D$ is nef and big, that is $D$ is nef and $D^n>0$, where $n:=dim X$.
\item{ii)} $\langle D\rangle$, keep in mind notations at page 
\pageref{pg:nota},
  has support with only normal crossing
(that is each $D_i$ is smooth and they intersect
everywhere transversally).

Then
$$H^j(X,K_X + \lceil D\rceil) = 0\textrm{   for   } j>0.$$
\label{van}
\end{Theorem}

Let us show how to use this vanishing theorem under {\bf very
special hypothesis}.

\noindent
{\bf Assume} {\sl that $L$ is ample (or nef and big) and that we
can find a divisor $D_1 \in |L|$ such that
$$
\pi^*(D_1)= \tilde{D_1} + aE
$$
with $c:= \frac{n}{a} <1$ and $\tilde{D_1}$ smooth. For instance 
assume that the only singularity of $D_1$
is an ordinary $(n+1)$-uple point at $x$.
Then $x$ is not a base point of $K_X + L$}

\smallskip \noindent
{\bf Proof} Note that for every $0 < \delta\ll 1$ we have that
$\pi^* (L)  - \delta E := A$ is ample.
Then we can write
$$ \pi^*(K_X+L) -E \equiv
K_Y + c\tilde{D_1} + caE -nE + (1-c)A + (1-c)\delta E
$$
equivalently
$$\pi^*(K_X+L) -E - (1-c)\delta E -c\tilde{D_1} -K_Y$$
is an ample $\Q$-divisor on $Y$.

We can apply the vanishing theorem on $Y$ and conclude that
$$H^1(Y,\pi^*(K_X+L) -E) = 0$$ since  $\lceil - (1-c)\delta E 
-c\tilde{(D_1)}\rceil=0 $.
Thus
$$H^0(Y,\pi^*(K_X+L))\to H^0(E,\pi^*(K_X+L))\iso\C,$$
is surjective and $x\not\in Bsl(K_X+L)$.
\bigskip

Unfortunately it is very unlikely that our special hypothesis are 
satisfied. Now it comes the moment to
give a precise general statement and to outline its proof.

\begin{Theorem}[Base point freeness, \cite{Sh1}, \cite{Ka} or \cite{KMM}]
\label{TeoA}
Let $X$ be a variety of dimension $n$, with ''good singularities'' 
(i.e. smooth or LT singularities, see the
definition
\ref{lc})
and $H$ a Cartier divisor.
Assume that $H$ is nef and $aH-K_X=:L$ is ample for some $a\in \N$. Then
for $m\gg 0$ the line bundle $mH$ is generated by global
sections,
i.e. there exists an integer $m_0$ and a  regular map $\f :X \ra W$
given by elements in $H^0(X,mH)$ for any $m\geq m_0$.

\end{Theorem}

\begin{remark} The above theorem was proved by Y. Kawamata and VV. Shokurov
(see \cite{Ka} and \cite{Sh1}) with a method which builds up from the 
classical methods of the
Italians and which was developed in the case of surfaces by Kodaira-Ramanujan-Bombieri.

A very significant step in the understanding and in the spreading out
of the technique was given
in a beautiful paper of M. Reid (see \cite{Re}) which we strongly suggest to
the reader.

This type of results are fundamental in algebraic geometry and
they are constantly under improvement, recently important steps were
achieved. among others by Kawamata,Shokurov, Koll\'ar and Ein-Lazarsfeld.

A big drawback is that the method, as it stands, is not
effective,
i.e. it does not give a good bound for $m$ (contrary to the case of curves
and surfaces).
Some bound can however be achieved, namely
one can show that $m_0 \leq 2(n+2)!(a+n)$
(Effective Base point freeness: see \cite{Ko}).
\end{remark}

We will only outline the proof and we refer to \cite{KMM} for many
technical, and often very relevant, parts which we now state and briefly
comment.

\medskip
First we observe that the ''perfect''
assumptions we have given above, are difficult to achieve
in general. So more than one blow up is
required and for this we need the following.

\begin{Definition}
For a pair $(X,H)$
of a variety
$X$ and a $\Q$-divisor $H$, a {\sf log resolution} is a proper
birational
morphism $f :Y\ra X$ from a smooth variety $Y$ such that the union of the
support of $f_*^{-1}H$ and of the exceptional locus is a normal
crossing divisor.
\end{Definition}

\begin{Theorem}Let $X$ be a variety with LT singularities,
$B$ an effective and nef $\Q$ divisor and $L$ an ample divisor on $X$.
Then there exists a log resolution $f:Y\ra X$
such that
$$\begin{array}{l}

K_Y=f^*K_X+\sum e_iE_i\\
f^* (B)=B'+\sum b_iE_i\\
f^*(L)=A+\sum p_iE_i\\

\end{array} $$

where all relevant divisors in $Y$ are smooth and normal crossing, 
all $E_i$ are
exceptional,
$A$ is an $f$-ample $\Q$-divisor, $0\leq p_i \ll 1$ and $e_i > -1$.
\label{hiro}
\end{Theorem}

The theorem follows essentially from the work of Hironaka
on resolution of singularities. The statement on the $e_i$
is the definition of LT singularities (see Definition \ref{lc})
while the ampleness of
$A$ is usually called Kodaira's lemma;
for a proof see \cite[corollary 0.3.6]{KMM}.

Using a log-resolution instead of the blow-up
we will achieve our assumption but we will very likely loose the non
vanishing part
(namely $H^0(E,\pi^*(K_X+L)) = \C$). For this we need the next very
important result, due to V.V. Shokurov.

\begin{Theorem} {\bf Non vanishing theorem}
Let $X$ be a non singular projective variety; let
$N$ be a Cartier divisor and $A$ a $\Q$-divisor on $X$ such that
\item{i)} $N$ is nef
\item{ii)}$\lceil A\rceil \geq 0$  and $\langle A \rangle$ has support with
only normal
crossing
\item{iii)} $dN+ A- K_X = M $ where $M$ is nef and big, for some positive
$d \in N$.

Then $H^0 (X, mN+\lceil A\rceil) \not= \emptyset $ for all $m\gg0$.
\label{snv}
\end{Theorem}

A proof of this theorem can be found in \cite[2.1.1]{KMM}. It is
a combination of the Riemann-Roch formula and the vanishing theorem \ref{van}

\bigskip \noindent
{\bf Sketch of the proof of \ref{TeoA}.}

\noindent
By the non vanishing theorem there exists an effective divisor
$B \in |m H|$ for all $m \geq m_0\gg0$.

\noindent
{\sl Noetherian argument:}
Let $B(\gamma)$ denote the reduced base locus of $|\gamma H|$.
Clearly $B(\gamma^s) \subseteq B(\gamma^t)$ for any positive integers $s>t$.
Noetherian induction implies that the sequence $B(\gamma^i)$ stabilises and we
call the limit
$B_{\gamma}$. So either $B_{\gamma}$ is non empty for some $\gamma$ 
or $B_{\gamma}$ and
$B_{\gamma^{\prime}}$ are
empty for
two relative prime integers $\gamma$ and $\gamma^{\prime}$. In  the 
latter case, take
$s,t$
such that $B(\gamma^s)$
and $B((\gamma^{\prime})^t)$ are empty and use the fact that every sufficiently 
large
integer is a linear
combination of $\gamma^s$ and $(\gamma^{\prime})^t$ with non negative 
coefficients to
conclude that $|mH|$ is base
point free for all $m\gg 0$.

So we must show that the assumption that some $B_{\gamma}$ is non 
empty leads to a
contradiction.
Let $m = \gamma^s$ such that $B_{\gamma}= B(m)$ and assume that this 
is not empty.

With $L$ as in the statement of the theorem and $B$ as at the beginning of
the proof,
let $e_i, b_i,p_i$ as in the theorem
\ref{hiro} and define
$$c:=\hbox{ min}\{\frac{e_i+1-p_i}{b_i}\}$$

By taking $m$ big enough we can assume that there exists a  divisor 
$B\in |mH|$ with
arbitrarily high multiplicity
along $B_{\gamma}$, in other words $0< c<1$.  By changing the
coefficients $p_i$ a little
we can assume that the minimum is achieved for
exactly one index. Let us denote the corresponding divisor by
$E_0$ and let $Z = f(E_0)$.

By Bertini theorem we can assume that
$Z$ is contained in the base locus of $mH$, i.e. in $B_{\gamma}$.

Then
$$\begin{array}{l}
K_Y + A + cB'+ \sum (cb_i - e_i + p_i)E_i + f^*(m-cm)H\\
\equiv f^*(K_X + L +mH)= f^*(m+a)H
\end{array}$$

and

$$\sum (cb_i - e_i + p_i)E_i = E_0 - D + Fr$$
where $E_0$, $D$ are effective divisor
without common irreducible components and $Fr$ is the fractional divisor
with rational coefficients between 0 and 1, defined by
$Fr = \sum \{cb_i - e_i + p_i \} E_i$, where $\{r\}$ is the fractional
part of the rational number $r$.

Thus
$$f^*((m+a)H) + D -E_0 -Fr - cB'-K_Y\equiv A + f^*(m-cm)H$$ is ample.
Denote $N(m) := f^*((m+a)H) +D$ for brevity,
by the vanishing theorem we have then
$$H^i(Y, N(m)-E_0) = 0 \hbox{  for } i > 0$$
(incidentally observe that also the following vanishing is true
$$H^i(E_0, N(m)) = 0 \hbox{  for } i > 0).$$

\par
By the first vanishing the restriction map
$$H^0(Y,N(m))\ra H^0(E_0,N(m)_{\vert E_0})$$
is surjective.
\par
By the non vanishing theorem for any $m_1 >> m_0$ there exists a non-zero
section $s$ of $N(m_1)_{\vert E_0}$.
By surjectivity this extends to a non-zero section of $N(m_1)$ on $Y$;
which is not identically zero along $E_0$.
Moreover  $H^0(Y,N(m_1)) = H^0(X,(m_1+a) H)$ since $D$ is $f$-exceptional.
The section  $s$ descends to a section of $(m_1+a)H$ which does not 
vanish along
$f(E_0) =Z \subset B_c$, which is a contradiction.

\section{Singularities and log singularities.}

In the previous sections we did not introduce any technical
definitions of singularities or of singular pairs.
Let us do it now for the interested reader.

\begin{Definition} Let $X$ be a normal variety and $D=\sum_id_iD_i$
be an effective
$\Q$-divisor such that $K_X+D$ is $\Q$-Cartier. If $\mu:Y\ra X$ is a
log resolution of the pair $(X,D)$, then we can write
$$K_Y+\mu_*^{-1}D=\mu^*(K_X+D)+F$$
with $F=\sum_jdisc(X,E_j,D)E_j$ for the exceptional divisors $E_j$. We call
$e_j:=disc(X,E_j,D)\in \Q$ the {\sf discrepancy coefficient} for $E_j$, and
regard
  $-d_i$
as the discrepancy coefficient for $D_i$.

The variety $X$ is said to have {\sf terminal} (respectively
{\sf canonical}, {\sf log terminal} (LT))
singularities if $e_j>0$ (resp. $e_j\geq 0$, $e_j>-1$),
for any $j$.

The pair $(X,D)$ is said to have {\sf log canonical} (LC)
(respectively {\sf Kawamata log terminal} (KLT)) singularities if
$d_i\leq 1$ (resp. $d_i< 1$) and $e_j\geq -1$ (resp.
$e_j>-1$) for any $i,j$ of a log resolution $\mu:Y\ra
X$.

The {\sf
log canonical threshold} of a pair $(X,D)$ is $lct(X,D):=
sup\{t\in\Q$: $(X,tD)$ is LC$\}$.
\label{lc}
\end{Definition}

\begin{Definition}[\cite{Ka1}] Let $X$ be a normal variety and $D=\sum
d_iD_i$ an effective $\Q$-divisor such that $K_X+D$ is $\Q$-Cartier.
A subvariety $W$ of $X$ is said to be a {\sf center of log canonical
singularities} for the pair $(X,D)$, if there is a birational morphism
from a normal variety $\mu:Y\ra X$ and a prime divisor $E$ on $Y$, 
not necessarily
$\mu$-exceptional, with
the discrepancy coefficient $e\leq -1$ and such that $\mu(E)=W$.
For another such $\mu^{\prime}:Y^{\prime}\ra X$, if the strict
transform $E^{\prime}$ of $E$ exists on $Y^{\prime}$, then we have the
same discrepancy coefficient for $E^{\prime}$. The divisor
$E^{\prime}$ is considered to be equivalent to $E$, and the
equivalence class of these prime divisors is called a {\sf place of
log canonical singularities} for $(X,D)$.
The set of all centers (respectively places) of LC singularities is
denoted
by $CLC(X,D)$ (resp. $PLC(X,D)$), the locus of all centers of LC
singularities is denoted by $LLC(X,D)$.
\end{Definition}

The study of these objects has been developed by Kawamata and we can summarise
the main results in the following Theorem.

\begin{Theorem}[\cite{Ka1},\cite{Ka2}] Let $X$ be a normal variety and $D$ an
effective $\Q$-Cartier divisor such that $K_X+D$ is $\Q$-Cartier. Assume
that $X$ is LT and $(X,D)$ is LC.
\begin{itemize}
\item[i)]
If $W_1,W_2\in CLC(X,D)$ and $W$ is
an irreducible component of $W_1\cap W_2$, then $W\in CLC(X,D)$. In
particular, there
exist minimal elements in $CLC(X,D)$ with respect to inclusion.
\item[ii)]
If $W\in CLC(X,D)$ is a minimal center then $W$ is normal
\item[iii)] (subadjunction formula)Let $H$ be an ample Cartier divisor
and $\epsilon$ a positive rational number. If $W$ is a minimal center
for $CLC(X,D)$ then there exists an effective
$\Q$-divisor $D_W$ on $W$ such that $(K_X+D+\epsilon H)_{|W}\equiv
K_W+D_W$ and
$(W,D_W)$ is KLT.
\end{itemize}
\label{clc}
\end{Theorem}
\begin{remark} The first two statements are, essentially, a consequence
of Shokurov Connectedness Lemma, which is itself a direct consequence of the
Vanishing theorem \ref{van}.
The subadjunction formula is quite of a different flavor and is
related to semipositivity results for the relative dualising sheaf of a
morphism.
\end{remark}

In particular \ref{clc} tell us that the minimal center $W$ is not 
too bad and there
is some hope to be able to work on it.

\begin{exercise} It is in fact not so difficult to work out all possible
minimal centers $W\in CLC(X,D)$, where
$X$ is a smooth surface and $D$ any divisor (i.e a curve).
The same, a little harder, if $X$ is a smooth threefold; one should 
keep in mind
that KLT singularities are rational singularities.
\end{exercise}

Let $(X,D)$ be a log variety and assume that $(X,D)$ is LC and $W\in CLC(X,D)$
is a minimal center. The Weil divisor $D$ is usually called boundary 
of the log pair.
Then we have a log resolution $\mu:Y\ra X$ with
$$K_Y=\mu^*(K_X+D)+\sum e_i E_i,$$
this time we put also the strict transform of the boundary on the right
hand side. Since $(X,D)$ is LC and $W\in CLC(X,D)$ then $e_i\geq -1$ and
there is at least one $e_j=-1$ such that $\mu(E_j)=W$.

A first problem
is that to apply Kawamata's BPF method we need to have one and only one
exceptional divisor with discrepancy $-1$ and  $W$ as
a center. To fulfill this requirement we need a

\noindent {\sl Perturbation argument:}
Choose a generic
very ample $M$ such that $W\subset Supp(M)$ and no other $Z\in CLC(X,D)
\setminus\{W\}$ is contained in $Supp(M)$, this is always possible
since $W$ is minimal in a dimensional sense. We then perturb $D$ to a divisor
$D_1:=(1-\epsilon_1)D+\epsilon_2 M$, with $0<\epsilon_i\ll1$ in such a way
that
\begin{itemize}
\item[-] $(X,D_1)$ is LC
\item[-] $CLC(X,D_1)=W$
\item[-] $\mu^*\epsilon_2 M=\sum m_i E_i +P$, with $P$ ample; this
is possible by Kodaira Lemma.
\end{itemize}
After this perturbation the log resolution looks like the following
$$K_Y+\sum_{j=0} E_j+\Delta-A=\mu^*(K_X+D_1)-P,$$
where the $E_j$'s are integral irreducible divisors and $\mu(E_j)=W$,
$A$ is a $\mu$-exceptional integral divisor and
$\lfloor \Delta\rfloor=0$.
It is now enough to use the ampleness of $P$ to choose just one of the $E_j$.
Indeed for small enough $\delta_j>0$ $P^{\prime}:=
P-\sum_{j=1} \delta_j E_j$ is still ample
therefore we produce the desired resolution
\begin{equation}
\label{pert}
K_Y+E_0+\Delta^{\prime}-A=\mu^*(K_X+D^{\prime})-P^{\prime};
\end{equation}
here and all through the paper after a perturbation
we will always gather together all the fractional
part with negative log discrepancy in $P$ and
$\Delta$, respectively the ample part of it and the remaining.

\noindent
If instead of an ample $M$ we choose
a nef and big divisor, we can repeat the above argument with Kodaira Lemma,
but this time we cannot choose
the center $\mu(E_0)$ like before,
and in particular we cannot assume that at the end
we are on a minimal center for $(X,D)$.

\subsection{How to use singularities and the CLC locus
to prove base point free-type theorems}
\label{bpf-method}

Assume now that $X$ is a variety with log terminal and Gorenstein singularities
and let $L$ be an
ample line bundle on $X$.

Let $D$ be an effective $\Q$-Cartier divisor
such that $D\equiv tL$ for a rational number $t<1$.
Let $W\in CLC(X,D)$ be a minimal center. Perturb $D$  using the very 
ample line bundle
$M:= mL$ for $m\gg 0$. So that we can assume that there exists only 
one exceptional divisor in any
log resolution of $(X,D)$ with
discrepancy $-1$ and  $W$ as a center.
Thus taking an embedded log resolution of the pair $(X,D)$,
$\mu:Y\ra X$ we have
$$K_Y+E+F=\mu^*(K_X+D)$$
where $E$ is a reduced divisor such that $\mu (E) =W$
and $F= \sum f_iF_i$ with $f_i <1$.
Then
$$K_Y +(1-t)\mu^*L \equiv \mu^*(K_X+L) -E -F$$
and thus
$$H^1(Y, \mu^*(K_X+L) -E+\lceil -F\rceil) =0$$
and we obtain a surjection
$$H^0 (Y, \mu^*(K_X+L) +\lceil -F\rceil) \ra H^0 (E, \mu^*(K_X+L) +\lceil
-F\rceil).$$
The divisor $\lceil - F\rceil$ is effective and any irreducible 
component of $\lceil F\rceil$  is
$\mu$-exceptional therefore $H^0 (Y, \mu^*(K_X+L))=H^0 (Y, 
\mu^*(K_X+L) +\lceil -F\rceil)$ and
we  also have
$$H^0 (Y, \mu^*(K_X+L)) \ra H^0 (E, \mu^*(K_X+L) +\lceil -F\rceil)\to 0.$$

Thus to find a section of $K_X+L$ not vanishing on $W$ it is
sufficient to find a non zero section in $H^0 (E, \mu^*(K_X+L) +\lceil
-F\rceil)$.

\bigskip
The {\bf ideal case} happens when $W = x$ is one point;
in fact in this case $H^0 (E, \mu^*(K_X+L) +\lceil -F\rceil)= \C$
and therefore $K+L$ is base point free at $x$.

\section{Exercises-Examples}

The solution of the next exercise can be found through the book \cite{BS},
even under the milder hypothesis that $L$ is ample and spanned.
We propose it here because we think that the above methods are
convenient to be applied at these problems and because we believe they
should bring to prove the conjecture stated in the item d)
(we do not know how and therefore we adopt the  trick to put it as an 
exercise).

\begin{exercise} Let $L$ be a very ample line bundle on a smooth
projective variety $X$ of dimension $n$. Prove the following:
\item{a)} $K_X +(n+1)L$ is spanned by global sections at each point.
\item{b)} The same is true for $K_X +nL$  unless $X = \Proj ^n$ and $L =
\O(1)$.
\item{c)} If $n \geq 2$ the same is true for $K_X+(n-1)L$
unless $X = \Proj ^n$ and $L = \O(1)$
or $X = \Proj ^2$ and $L = \O(2)$ or $X = \Q^n$ and $L = \O_{\Proj 
^{n+1}}(1)_{\Q^n}$
or $(X,L)$ is a scroll over a curve.
\item{d)} {\rm Conjecture:} If $n \geq 3$ the same is true for 
$K_X+(n-2)L$ as soon as
it is nef and $L^n > 27$.
\end{exercise}

\noindent
{\sl hints} For the question a) let $x\in X$ and take $n$-sections of $L$
meeting transversally
in $x$

For b) use an ''induction procedure''; namely take
a smooth section $D \in |L|$ passing through $x$ (this is Bertini theorem)
and use the exact sequence
$$H^0(X,K_X+nL) \ra H^0(D, K_D+L) \ra 0.$$
One goes down until the dimension of $D$ is $1$, i.e. a curve,
and in this case $K_D+L$ is spanned if and only if $degL \geq 2$.
The only problem is then when $D$ is a line and therefore
  $X = \Proj ^n$ and $L = \O(1)$.

For c), as in the previous step, one can reduce the problem to the 
surface case;
namely $X=S$ is a
smooth surface and one has to prove the spannedness of $K+L$.
In this case there are even stronger theorems (Reider type theorems).

Some comments to the conjecture stated in d): by the inductive procedure
it is enough to prove the proposition for $n =3$. Note that the bound
$L^n > 27$ is necessary since there exists a del Pezzo 3-fold $X$
with $-K_X = 2H$, $H^3 =1$ and $H$ with one base point (take $L = 3H$).

\bigskip
The above exercise is extremely hard when one assumes only ampleness
(and not very ampleness!) of $L$. In fact we have:

\begin{conjecture}
[Fujita conjecture] Let $L$ be an ample line bundle on a smooth projective
variety of
dimension $n$. Then $K_X+mL$ is base point free if $m \geq n+1$ and it is
very ample if
$m\geq n+2$.
\end{conjecture}

\begin{remark} Some very important results toward a proof of the conjecture
were found in recent time. In particular, using an analytic approach,
Demailly, Angern-Siu and Tsuji
proved that if $m \geq {{n+1}\choose{2}}$ then $K_X+mL$ is base
point free
and if $ m \geq {{n+2}\choose{2}}$ then the global sections of
$K_X+mL$ separate points.

The base point free part of the conjecture is true in
the case $n=(1),2,3,4$ by results of Reider,
Ein-Lazarsfeld, Helmke, Kawamata and Fujita
(see \cite{Rei} and \cite{Ka1}).
\end{remark}

\part{Fano-Mori or extremal contractions}
\setcounter{section}{0}

In this Part we first define and give examples of Fano--Mori
spaces. These are exactly
the morphisms constructed in Part \ref{ch:BPF}, and they play a central
role in the Minimal
Model Program.
To study those objects we want to apply an inductive method as in
section \ref{sec:fuj}. A fundamental step is therefore to ensure that we
have base point free
linear systems to slice the fibers. This is the content of Theorem
\ref{aw},
whose proof occupies the last section.

\section[Contractions of rays]{Contractions associated to a ray of the
Kleiman--Mori cone }

A key step in Mori theory, after the description of the structure of
$\overline{NE(X)}$ outlined in a previous section,
is the fact that extremal rays (and in general extremal faces)
give raise to morphisms of the variety.
This is explained in this section.

\begin{Proposition} Let $R$ be an extremal ray of the Kleiman-Mori cone
$\overline{NE(X)}$ such that $R\cdot K_X < 0$. Then there exists
a nef Cartier divisor $H_R$ such that ${H_R}\cdot z= 0$ if and only 
if  $z\in R$.
\label{H}
\end{Proposition}

This proposition is proved for instance in \cite[III.1.4.1]{RaC}.
The proof makes use of the Cone theorem and some easy properties of
closed cones.

Then to a divisor as in the proposition we can associate a morphism
via the following theorem.

\begin{Theorem} [Contraction theorem]
Let $X$ be a variety with log-terminal singularities and let $H$ be a
nef Cartier divisor on $X$.

Assume that $F:= H^{\bot} \cap \overline {NE(X)} \setminus \{0\}$ is contained
in $\{C \in N_1(X): K_X\cdot C <0\}$. Then there exists a projective morphism
$\f :X \ra W$ onto a normal projective variety $W$ which is characterised
by the following properties

\begin{itemize}
\item[i)] For any irreducible curve $C \subset X$, $\f(C)$ is  a point if and
only if $H\cdot C= 0$.
\item[ii)] $\f$ has connected fibers
\item[iii)] $H = \f^*(A)$ for some ample Cartier divisor on $W$.
\end{itemize}
\label{F-Mcontr}
\end{Theorem}

\begin{proof} The proof follows immediately from the theorem \ref{TeoA}
and Zariski's Main theorem once we note that by our assumption and
Kleiman's criterion for ampleness there exists a natural number $a$ such that
$aH-K_X$ is ample.

\end{proof}

\begin{Definition}
A {\sl contraction} is a surjective morphism $f:Y\ra T$, with connected fibers,
between normal varieties.

For a contraction $f:Y\ra T$ the set
$$ E = \{y \in Y; f \hbox{\ is\ not\ an\ isomorphism\ at\  }y\}$$
is the exceptional locus of $f$. Let $\delta = dim E$ where $dim$
denotes as usual the maximum of the dimension of the irreducible components.

\noindent $f$ is called of {\sl fiber type} if $\delta = dim Y$ or 
{\sl birational}
otherwise.

If $f$ is birational and $\delta = dim Y -1$
then it is also called {\sl a divisorial contraction};
if it is birational and $\delta \leq dimY -2$
then it is called a {\sl small contraction}.

For a contraction $f:Y\ra T$ a Cartier divisor $H$ such that
$H = \f^*(A)$ for some ample Cartier divisor on $T$ is called
a {\sl supporting divisor} for the contraction (if $H = H_R$ as in the above
proposition then it is also called a supporting divisor for the
ray $R$).
\end{Definition}

\begin{Definition}
A contraction $f:X\ra W$ as in the above Theorem \ref{F-Mcontr}
is called {\sl Fano--Mori} (F-M) or {\sl extremal}.

A birational contraction $f:X\ra W$ is called {\sl crepant} if
$K_X = f^*K_W$
\label{F-M}
\end{Definition}

\begin{remark} Putting together Theorem \ref{F-Mcontr} and Proposition \ref{H}
we obtain the following. Given an extremal ray of the Kleiman-Mori cone
$R \subset \overline{NE(X)}$ such that $R\cdot K_X < 0$, there exists
a projective morphism with connected fibers
$cont_R :X \ra W$ onto a normal projective variety $W$, which
contracts all (and only) the curves in the ray.

Such a map is also called {\sl the contraction of the extremal ray $R$},
or an {\sl elementary Fano-Mori contraction}.

We stress that Theorem \ref{F-Mcontr} is proved
only in characteristic zero. The existence of this map
is an open problem in positive characteristic.
\label{contR}
\end{remark}

\begin{remark}
It is straightforward to prove that conversely Contraction theorem
implies the Theorem \ref{TeoA}.

Note also that any supporting divisor $H$ for a F-M contraction $\f$
is  of the type $K_X + r L$ with $r$ a rational number
and $L$ an ample Cartier divisor.
In fact let $H$ be a Cartier divisor which is the pull
back of
a sufficiently ample line bundle on $W$. Then $mH-K_X :=L$ is an ample Cartier
divisor for
some rational number $m$ and thus $H= K_X +1/m L$.
\label{supporting}
\end{remark}

\begin{remark}
To construct a divisor as in \ref{F-Mcontr}, and therefore
an associated morphism, one can also use the Rationality theorem
\ref{Ka-rat} as follows. Let $X$ be a variety with at most log-terminal
singularities and let $L$ be a Cartier divisor with nef value $r$.
Then, by the rationality theorem, if $H' := K+rL$ there exists an integer $m$
such that $H:=mH'$ is a Cartier divisor. By definition $H$  satisfies the
assumption in \ref{F-Mcontr}.
\label{renef}
\end{remark}

The following is an important technical result whose proof
may be considered an interesting exercise.

\begin{exercise}(\cite[Proposition 5.1.6]{KMM}).
Let $f: X\ra W$ be a divisorial elementary Fano Mori
contraction
with $X$ smooth or with at most terminal $\Q$-factorial 
singularities. Prove that
the exceptional locus
of $f$ is a unique prime divisor and $W$ has at most terminal $\Q$-factorial
singularities.

\hint Assume by contradiction that there are at least two components. 
Show that a generic
curve in one component cannot be numerically equivalent to a generic 
curve in the other.
\label{divisorial}
\end{exercise}

\subsection {Local Contraction}

In studying F-M contractions it makes sense to fix a
fiber and understand the contraction locally, i.e.
restricting to an affine neighborhood of the fixed fiber.
More general complete F-M contractions can then be obtained by gluing
different local descriptions.

For this we use the local set-up developed by Andreatta--Wi\'sniewski, see
\cite{AW},
which depends on some definitions.

\begin{Definition} Let $f:Y\ra T$ a contraction supported by $K_Y+rL$,
with $r$ rational and $L$ ample and Cartier (i.e. a F-M contraction).
Fix a fiber $F$ of $f$ and take an open affine
$S\subset T$ such that $f(F)\in S$ and $dim f^{-1}(s)\leq dim F$, for $s\in S$.
Let $X=f^{-1}S$
then $f:X\ra S$ will be called a {\sf local contraction around
$F$}. If there is no need to specify fixed fibers then we will simply
say that $f:X\ra S$ is a local contraction.
In particular $S=Spec(H^0(X,\Ol_X))$.
\label{defloc}
\end{Definition}

\begin{Definition} Let $f:X\ra S$ a local F-M contraction
around $F$. Let $r=\mbox{\rm inf}\{t\in \Q: K_Y+tH\nel{f} 0$ for some
ample Cartier divisor 
$H\in Pic(X)\}$. Assume that $K_X+rL\nel{f}\O_X$, that is 
$f$ is supported by $m(K_X+rL)$ for some $m\geq 1$. 
The  Cartier  divisor $L$
will be called {\sf fundamental divisor} of $f$.
Let $G$
a generic non trivial fiber of $f$.

\noindent
	The {\sf dual-index} of $f$ is
$$d(f):=dimG-r,$$
the {\sf character} of $f$ is
$$\gamma(f):=
\left\{\begin{array}{ll}
1 & \mbox{\rm if $dim X>dim S$} \\
0 & \mbox{\rm if $dim X=dim S$}
\end{array}
\right.$$
  and the
{\sf difficulty} of $f$ is
$$\Phi(f)=dim F-r.$$
We will say that
$(d(f),\gamma(f),\Phi(f))$ is the {\sf type} of $f$.
\label{tipo}
\end{Definition}

\section{Examples.}

A large class of examples of F-M contractions is worked out in section 3 of
the paper
\cite{AW3}; we report some of them here, referring the reader for more details
to that paper.

We focus on the case $X$ is smooth,
with the purpose of showing later some classifications
of F-M contractions on a smooth variety.

\begin{example}
Fano varieties (with the constant map $X \ra \{pt\}$,
Scrolls (i.e. $X = \Proj(E) \ra Y$ where $E$ is a vector bundle on
a smooth manifold $Y$),
conic bundles are F-M contractions of fiber type.
\end{example}

\begin{example}
Any blow up of a smooth smooth variety $Y$ along a smooth
subvariety $Z$,  $X:= Bl_{S}Z \ra Z$, is a birational F-M contraction.
\end{example}

\begin{example}

Blow-up a smooth surface in a 4-fold with
an ordinary double point; i.e.
$$S := \{x=z=w=0\} \subset Z := \{xy-zt+w^2\}$$
$$\f : X:= Bl_{S}Z \ra Z.$$
A direct computation shows that $X$ is smooth and that
$\f^{-1} (0) = \Proj^2$.

\medskip
Let $L_1,L_2,L_3$ be three general planes in $\Proj^3$
and let $\Proj^2$ be the base of the net ${\mathcal L} = \Sigma t_i L_i$.
Consider the incidence variety
$$X := \{(p,L); p\in L\}\subset \Proj^3\times\Proj^2.$$

Then the projection $\f:X\ra \Proj^3$ is
a F-M contraction which is a $\Proj^1$-bundle generically and has a fiber
$=\Proj^2$ over the point intersection of the $L_i$.

If we blow-up a smooth surface $S$ in $X$ meeting the general fiber in one
point
we obtain a smooth conic bundle $Y \ra \Proj^3$ with a two dimensional
reducible fiber
and with discriminant locus $\Delta = \f(S)$.

{\sl In coordinates:}
assume $\Proj^3 = [z_0,z_1,z_2,z_3]$, $\Proj^2 = [t_1,t_2,t_3]$ , 
$L_i = z_i$, $i=1,2,3$.
Then $X  = \{ t_1z_1+t_2z_2 +t_3z_3 =0\} \subset \Proj^2 \times \Proj^3$ and
let, for instance, $S= \{t_1 = z_1 = 0\}$. The special two 
dimensional fiber on $Y$ will be
$\F_1 \cup \Proj^2$.

On $Y$ there are two F-M contractions, both of birational type; 
besides the blow up of $X$ along
$S$ we can contract a divisor on $Y$ consisting of the $\Proj^2$ 
component of the two dimensional fiber
and of all the components of the reducible conics not contracted to 
$X$. This is a contraction
as the one described in the first part of the example
(if this is not immediate now, it will be later
when we will give a classification of F-M contractions on smooth $4$-folds).

\end{example}

\begin{example}
\label{examplev}
We now introduce a large class of examples via a standard construction;
for more details see section 3 of \cite{AW3}.

Let ${\mathcal E}$ be a vector bundle over a smooth variety $F$ and let
${\bf V}({\mathcal E}): =  Spec (S({\mathcal E}))$ be the {\sl total space
of the dual}
${\mathcal E}^*$.

If $S^k({\mathcal E})$ is {\sl generated by global sections} for some $k
>0$ let
$$\f :{\bf V}({\mathcal E}) \ra Z = Spec (\bigoplus _{k\geq 0}
H^0(F,S^k(\mathcal E)),$$
be the map associated to the evaluation of $S^k(\mathcal E )$.
$\f$  is a contraction. The collapsing of the zero section , $F_0:= F$,
of the total space
${\bf V}(\mathcal E)$ to the vertex $z$ of the cone $Z$.

It is straightforward to check the following properties:

i) The normal bundle of $F_0$ into ${\bf V}(\mathcal E)$ is ${\mathcal E^*}$.

ii) If $-K_Y-{det \mathcal E}$ is ample then $\f$ is a Fano-Mori
contraction. The map
is birational if the top Segre class of
${\mathcal E}$ is positive (if $rank {\mathcal E} = 2$ then $c_1^2 - c_2 > 0$).

iii) $\Proj(\O \oplus {\mathcal E}) := Proj (S ({\mathcal E} \oplus \O_Y))$
is the projective closure of ${\bf V}({\mathcal E})$.  The map $\f$ is the
restriction of
the map given by the tautological bundle
$\xi$ on $Proj (S ({\mathcal E} \oplus \O_Y))$; $\f$ is birational
if $\xi$ is big.

iv) (Grauert criterion). ${\mathcal E}$ is ample if and only if
the map $\f$ is an isomorphism outside $F_0$.

v) The fiber $F_0$ of the map $\f$ has the fiber structure
(i.e. $\I_{F_0} = \f^{(-1)}m_z \O_X$)
if and only if ${\mathcal E}$ is spanned by global sections.

Let us work out in details the example with $F = \Proj^2$; it
is possible to do the same for a two dimensional quadric,
see \cite{AW3}, or for smooth
del Pezzo surfaces.
Let $\E$ be a rank-2 vector bundle over $\Proj^2$ such that
$\E$ is spanned by global sections and $0\leq c_1(\E)\leq 2$.
These bundles were completely classified in \cite{SW}, and
they are isomorphic to one of the bundles in the following table.
Performing the above construction with them we obtain eight
Fano-Mori contractions with fiber $\Proj^2$;in
the second column we describe these contractions.
In \cite{AW3} it is explained how to obtain these
descriptions; one has to use the results \ref{graded}, \ref{castelnuovo},
\ref{formaln} in the next section.
$$\begin{array}{ll}
\textrm{description of bundle } \E& \textrm{description of } \f
\textrm{ and } Sing(Z)\\

\hline

\E=\O\oplus\O&\textrm{a scroll, $Z$ is smooth}\\

\E=\O\oplus\O(1)&\textrm{a smooth blow-up of a smooth}\\
&\textrm{curve,$Z$ is smooth}\\

\E=T\Proj^2(-1)&\textrm{a generalised scroll,i.e. a fiber type map}
\\&\textrm{general
fiber isomorphic to a line and a }\\
&\textrm{a two dimensional fiber, $Z$ is smooth}\\
\E=\O\oplus\O(2)&\textrm{the blow-up of a smooth curve $C$}\\
&\textrm{$Z$ is singular along $C$ }\\
\E=\O(1)\oplus\O(1)&\textrm{a small contraction, }\\
& \textrm{Z is singular and the flip exists}\\
0\ra \O\ra T\Proj^2(-1)\oplus\O(1)\ra\E\ra 0 &\textrm{the
blow-up of a smooth surface}\\
&\textrm{(passing from) a quadric}\\
&\textrm{singularity of } Z\\

0\ra\O(-1)^{\oplus 2}\ra\O^{\oplus 4}\ra \E\ra 0 &
\textrm{the blow-up of cone over a twisted}
\\&\textrm{cubic in a smooth } Z\\
0\ra\O(-2)\ra\O^{\oplus 3}\ra \E\ra 0 &\textrm{a conic bundle with a two}\\
&\textrm{dimensional fiber, Z is smooth}\\
\end{array}$$
\end{example}

\begin{example}
The existence of F-M birational contractions with
exceptional set of codimension greater than 1 (small contraction)
was proved by P. Francia with a famous example: it is a F-M contraction on a
$3$-fold with terminal singularities and with exceptional
locus $E \cong {\mathbb P} ^1$. The example is worked out
in many books, for instance in (\cite{CKM}, p.33-34).
This is a main difficulty in the MMP over-passed by S. Mori,
with a tremendous work, in dimension $3$ (see \cite{3flip}).

\end{example}

\section[Relative BPF on F--M]{Relative Base point freeness on Fano-Mori
Contractions}
\label{relgood}

A F-M contraction has a supporting divisor
of the type $K_X+ rL$ with $L$ an ample Cartier divisor,as noticed in
\ref{supporting}.

This feature, for which we can also call this {\sl adjoint 
contraction morphism},
allows us to apply an inductive method
which is typical of this theory.
It is a sort of relative {\sl ``Apollonius method''}, see \ref{sec:fuj}, and
in \cite{AW} it is called {\sl horizontal
slicing} argument (sometimes it is called simply {\sl slicing} but we
will need to distinguish it from {\sl vertical slicing}).
It can be briefly summarised as follows.
\par

Consider a general divisor $X'$ from the linear system $\vert L\vert$
(a hyperplane section of $X$ if $L$ is very ample) and assume that it
is a ``good'' variety, i.e. has the same singularities as $X$,
of dimension $n-1$.
By adjunction, $K_{X'}= (K_X+L)_{\vert X'}$ and, by the Vanishing
Theorem \ref{van}, if $r>1$, the linear system $\vert m(K_{X'}+(r-1)L) \vert$
is just the restriction of $\vert m(K_X+rL)\vert$, so that
the {\sl adjoint contraction morphism} of $X'$ can be related
to the one of $X$.
Moreover, fibers of the {\sl adjoint morphism} of $X'$ will be usually
of smaller dimension  and an inductive argument can be applied.
The method will be further outlined in the section \ref{chap:ladder}.
\par
The {\sl horizontal slicing} argument requires therefore the existence
of a ``good'' divisor $X'$
in the linear system $\vert L\vert$ (a {\sl rung} in the language of
\cite{Fu}, see section \ref{sec:fuj}).
The system, however, for an ample (but not very ample) $L$
may a priori be even empty.
To overcome this difficulty we use the local set-up,
described in the previous section, in which the
base of the {\sl contraction morphism} will be affine. We benefit
from this situation also because we may choose effective divisors
which are rationally trivial.

Then the next point
is to ensure that the divisor $X'$
does not contain the whole fiber in question
and has good singularities.
This is the case for instance, via Bertini theorem, if we ensure that base
  locus of
$\vert L\vert$ ($L$ may be changed by adding a divisor
trivial on fibers of $\f$) is empty.
This is what may be called a ``relative good divisor''.
(Now we can explain why we use the word ``horizontal'': we are used to think
about the map $\f: X\ra Z$ as going vertically,
then every divisor from an ample linear system cuts every ``vertical''
fiber of $\f$ of dimension $\geq 1$, so it lies ``horizontally''.)

\smallskip
The above point of view was first exploited in \cite{AW}
where the first part of the following theorem was proved.
The proof used the Base Point Free Theorem method (BPF-method) of
Y. Kawamata, actually a slightly improved version of it by J. Koll\'ar,
see \cite{Ko0}, introduced in the previous section.
The further refinement of the method by Y. Kawamata
in \cite{Ka1} and \cite{Ka2} allowed an improvement of the
theorem in \cite{AW}, this is the second part of the following
theorem and it was proved in \cite{dpf}.

\begin{Theorem} Let $f:X\ra S$ be a local F-M space around $F$ supported by
$K_X+rL$ and let $(d(f),\gamma(f),\Phi(f))$ be the type of $f$
(see \ref{defloc} and \ref{tipo}). Let also $\epsilon$ be
a sufficiently small positive rational number.

Assume one of the two following conditions is satisfied
\begin{itemize}
\item  $dimF < r+1$ or, if $f$ is birational, $dimF \leq r+1$;
equivalently the type of $f$ is
$(*,*,\Phi(f))$, with $\Phi(f)\leq 1-\epsilon\gamma(f)$ (see \cite{AW}),
\item the type of $f$ is $(d,1,1)$, with $d\leq 0$ or with $d=1$
and $F$ is reducible (see \cite{dpf}).
\end{itemize}
Then $L$, the fundamental divisor of the contraction, is relatively spanned,
i.e $Bsl|L|:=Supp(
Coker(f^*f_* L\ra L))$ does not meet $F$.
\label{aw}
\end{Theorem}

In the rest of the section we are going
to prove this theorem.
Let us first roughly summarise the general principles of the proof.
The idea is to proceed by contradiction, we assume therefore that there is
a non empty base locus $V$. Then we produce a log variety non KLT on $V$
(with respect to a divisor in $\delta L$).
Finally we use the method developed in Part
\ref{ch:BPF} to
produce sections of an adjoint line bundle non vanishing along the non KLT part
of the log variety.

To apply this strategy we have a priori a main problem: namely BPF produces
sections
of $K_X+mL$, for $m\gg 0$, while we need sections of $L$ itself.
But in the category of local F--M contractions we have that
$L \nel{f} K_X + (r+1) L$ .
An immediate consequence is the following.

\begin{CO}
In our set up of F-M contraction, we will work with log pairs $(X,D)$ such that
$D\nel{f} \delta L$ and for which $K_X\nel{f} -rL$.
In particular, by subadjunction
formula, see Theorem \ref{clc} part (iii), we have
$$K_W+D_W\nel{f} (\delta-r+\epsilon)L.$$
\end{CO}
So if $W$ is contained in a fiber and
$\delta<r$ then $K_W+D_W$ is antiample.

\begin{Definition}  A {\sf log-Fano variety} is a KLT pair $(X,\Delta)$
  such that for some positive
integer $m$,
$-m(K_X+\Delta)$ is an ample Cartier divisor. The index of a log-Fano
variety $i(X,\Delta):=sup \{t\in \Q: -(K_X+\Delta)\equiv tH$ for some
ample Cartier divisor $H \}$ and the $H$ satisfying
$-(K_X+\Delta)\equiv i(X,\Delta)H$ is called fundamental divisor.
\label{fanodef}
\end{Definition}

 From our point of view these varieties are extremely important because
we have a simple effective non vanishing, directly coming from Hilbert
polynomial.

\begin{Proposition}[\cite{Al},\cite{Am}] Let $(X,\Delta)$ be a
log-Fano n-fold of index $i(X)$, $H$ the fundamental divisor in $X$.
If $i(X)> n-3$ then $h^0(X,H)> 0$, moreover if
$i(X)\geq n-2$ then $h^0(X,H)>1$.
\label{al}
\end{Proposition}
\begin{proof} For simplicity assume that $\Delta=0$
and $i(X)\geq n-2$, the other cases are treated similarly putting some more
effort.
Let $p(t):=\chi(X,tH)=\sum h_j t^j$ the Hilbert polynomial
of $H$ and $d=H^n$ (see section \ref{sec:fuj}).
In particular
$$h_n=d/n!$$ and
$$h_{n-1}=\frac{-K_X\cdot H^{n-1}}{2(n-1)!}=\frac{i(X)d}{2(n-1)!}.$$
By the vanishing theorem \ref{van}
$$H^i(X,tH)=H^i(X,K_X+(tH-K_X))=H^i(X,K_X+(i(X)+t)H)=0$$
for $i>0$ and
$t>-i(X)$.
On the other hand, $H$ is an ample
divisor therefore
$$H^0(X,tH)=0 \text{\ \ for any $t<0$.}$$
Combining the two we obtain that
$$p(t)=0 \text{\ \ for $-i(X)<t<0$},$$
and $p(1)=1$.
Plug this informations into $p(t)$ to get
\begin{eqnarray*}
p(t)&=&\frac{d}{n!}(t+1)(t+2)\ldots(t+n-2)(t^2+at+\frac{n(n-1)}{d})\\
     &=&\frac{d}{n!}t^n+\frac{d}{n!}(a+\frac{(n-2)(n-1)}{2})t^{n-1}+\ldots
\end{eqnarray*}
To determine $a$ use
$$h_{n-1}=\frac{i(X)d}{2(n-1)!}.$$
So that
$$a= \frac{ni(X)-(n-2)(n-1)}{2}.$$
Which yields $h^0(X,H)=p(1)>d/n+(n-1)>1$.

\end{proof}

The next Lemma  translates Proposition \ref{al} in the non vanishing 
theorem we need.

\begin{Lemma} Let $f:X\ra S$ be a local contraction
supported by $K_X+rL$ around $F$. Fix
a subvariety $Z\subset F$, and  a
$\Q$-divisor $D$, with $D\nel{f}
\gamma L$. Assume that $X$ is LT,
$(X,D)$ is LC along $Z$, and
$W\in CLC(X,D)$ is a minimal center contained in $Z$.
Assume that one of the following conditions is satisfied:
\begin{itemize}
\item[i)] $r-\gamma> max\{0,
dim W-3\}$,
\item[ii)] $dim W\leq 1$ and $r-\gamma> -1$.
\end{itemize}
Then there exists a section of $|L|$ not vanishing identically on $W$.
\label{sectame}
\end{Lemma}
\begin{proof}
Since $D$ is LC along $W$ we can assume, up to a perturbation, that there
exists a log resolution $\mu:Y\ra X$  of $(X,D)$  with
$$K_Y-A+E+\Delta+B=\mu^*(K_X+D)-P,$$
where:
\begin{itemize}
\item[-] $E$ is an irreducible integral divisor,
\item[-] $A$ and $B$ are integral divisors,
\item[-] $\Delta$ and $P$ are
$\Q$-divisors.
\end{itemize}
Furthermore these divisors satisfy the following properties:
\begin{itemize}
\item[-] $\mu(E)=W$,
\item[-] $A$ is \hbox{$\mu$-exceptional},
\item[-] $\lfloor \Delta\rfloor=0$,
\item[-] $Z\cap\mu(B)=\emptyset$
\item[-] $P$ is $(f\circ\mu)$-ample.
\end{itemize}

Let
\begin{equation}
\label{bien}
N(t):=\mu^*tL+A-\Delta-E-B-K_Y\nel{f\circ\mu}
\mu^*(t+r-\gamma)L+P,
\end{equation}
  then $N(t)$ is $(f\circ
\mu)$-ample whenever $t+r-\gamma\geq 0$.
In particular if conditions i) or ii) of the Lemma are satisfied, by
vanishing theorem \ref{van}, we have the following surjection

$$
  H^0(Y,\mu^*L+A-B)\ra H^0(E,(\mu^*L+A)_{|E}).
$$

Since $A$ does not contain $E$ and is effective then

$$H^0(W,L_{|W})\hookrightarrow H^0(E,(\mu^*L+A)_{|E}).$$
  In particular
any section of $H^0(W,L_{|W})$ gives rise to a section
in $H^0(X,L)$
not vanishing identically on $W$.
Therefore to conclude the proof it is enough to produce a section in
$H^0(W,L_{|W})$.
By subadjunction formula of
Theorem \ref{clc} there exists a $\Q$-divisor $D_W$ such that
\begin{equation}
K_W+D_W\equiv (K_X+D+\epsilon L)_{|W}\equiv -(r-\gamma-\epsilon)L_{|W},
\label{subadjeq}
\end{equation}
for any $0<\epsilon\ll1$.

In case $(i)$ since $r-\gamma>0$ then by equation (\ref{subadjeq}), for
sufficiently small $\epsilon$,
$(W,D_W)$ is a log Fano variety of index
$i(W,D_W)=r-\gamma-\delta>dim W-3$.
Therefore we can apply Proposition \ref{al}.

If $dim W=1$ then $W$ is smooth. Since $r-\gamma-\epsilon> -1$
  by relation
(\ref{subadjeq})
$$0<L\cdot W\geq 2g(W)-2$$
  thus $h^0(W,L_{|W})>0$ by R--R
formula.
\end{proof}

We have to make the last preliminary to the proof of Theorem \ref{aw}. Till
now we always
worked with LC pairs. Along the proof we use pairs $(X,D)$ which
are not LC. To be able to treat this situation let us introduce the
following definition and make some useful remarks.

\begin{Definition}
The {\sf log canonical threshold related to a
scheme $V\subset X$} of a pair $(X,D)$ is $lct(X,V,D):=\mbox{\rm inf}\{t\in \Q:
  V\cap LLC(X,tD)\neq \emptyset\}$. We will say that $(X,D)$ is LC along
a scheme $V$ if $lct(X,V,D)\geq 1$.
\label{lcsub}
\end{Definition}

\begin{remark}Let $Z\in
CLC(X,lct(X,V,D)D)$ be a center and assume that $Z$ intersects $V$, then
$(X,lct(X,V,D)D)$ is LC on the generic point of $Z$.

\noindent If $(X,D)$ is not LC then
Theorem \ref{clc} is in general false. On the other hand
the first assertion  stays true,
also under the weaker hypothesis that $(X,D)$ is LC on the generic point
of $W_1\cap W_2$. In fact the discrepancy is a concept
  related to a
valuation $\nu$, therefore we can always
substitute the variety $X$ by an affine neighborhood of the generic
point of the center of $\nu$.
\label{rem:clc}
\end{remark}

\begin{proof}[Proof of theorem \ref{aw}]
Let $V=Bsl|L|\cap F$,
  remember that we are in a relative situation, therefore we
need always to consider objects contained in a fixed fiber to fully enjoy
the geometrical consequences of the ample anticanonical class.

Our aim is to derive a contradiction producing a section of $L$
which is not identically vanishing along $V$.
Consider the set $\D=\{ D\}$ of $\Q$-divisors $D$ such that:
\begin{itemize}
\item[-] $D\nel{f} \delta L$,
\item[-] there exists a minimal center $W_D\in CLC(X,D)$ such that 
$W_D\subset F$
and $W_D\cap V\neq\emptyset$,
\item[-] $dim W_D\leq r+1-\delta$,
\item[-] $lct(X,W_D,D)=1$.
\end{itemize}

First observe that $\D$ is non empty. Consider
$D_0=f^*\sum_I l_i(g_i)$, for $g_i$ generic functions on $S$
vanishing at $f(F)$. Then $D_0\nel{f} 0$ and one can choose
$0<l_i\ll1$ such that $lct(X,V,D_0)=1$.

\begin{claim} There exists a $D\in \D$ such that
$W_{D}\subset V$. Furthermore if $D\nel{f} (r+1)L$ one can choose $D$ so that
  $D=D_0+\sum_1^{r+1} H_i,$ with $H_i\in|L|$ generic.
\label{cl:clc}
\end{claim}
\begin{proof}[Proof of the claim] Consider the above $D_0$ and let
$H\in |L|$ be a generic
section. Let
$$c=inf\{t\in \Q^{\geq 0}: LLC({D_0}+tH)\cap V\cap W_{D_0}\neq\emptyset\}.$$

Since $H$ is a Cartier divisor vanishing on $V$,
then $c\leq 1$. Let $D_c=D_0+cH$.

If $c<1$ we assert that there exists a minimal center $W_{D_c}\in 
CLC(X,D_c)$ with $W'\subset V$.
Let us spend a few words on this.
Fix a resolution  $g:Y\to X$ of the singularities of $X$. Let 
$g^*H=H_Y+G$, then by Bertini
Theorem $H_Y$ is smooth outside $Bsl |H_Y|$. Furthermore for any 
$g$-exceptional divisor $A$ such that
$g(A)\not\subset Bsl|L|$ we can choose an $H\in |L|$ such that 
$Supp(H)\not\supset g(A)$.
There are finitely many  $g$-exceptional divisors in $Y$, therefore 
$g(G)\subset Bsl|L|$.
Let now $h:Z\to Y$ be a log resolution of $(Y,H_Y)$, so that 
$f:=g\circ h$ is a log resolution of
$(X,H)$. Let $f^*H=H_Z+\Delta$, then $h(\Delta)\subset Bsl|H_Y|\cup G$. Hence
$f(\Delta)=g(h(\Delta))\subset Bsl|L|$. As a consequence
$LLC(X,D_c)\subset Bsl|L|\cup LLC(X,D_0)$. Furthermore for any 
$0<\epsilon$, $(X,D_c+\epsilon H)$ is
  not LC along $V\cap W_{D_0}$, therefore
there  exists a center $W'\in CLC(X,D_c)$
with $W'\cap (V\cap W_{D_0})\neq\emptyset$ and $W'\cap F\subset V$.

To conclude consider a minimal center $W_{D_c}$ contained in
$W'\cap W_{D_0}\subset V$,
keep in mind Remark \ref{rem:clc}.

If $c=1$ then both $W_{D_0}$ and $H$ are in $CLC(X,D_1)$, and their
intersection is not empty because $W_{D_0}\cap V\neq\emptyset$. Therefore
by Remark \ref{rem:clc} any irreducible component $Z\subset W_{D_0}\cap H$ is
a center. Furthermore  $dim Z=dim W_{D_0}-1$.
This means that $D_1\in \D$  and $dim W_{D_1}<dim W_{D_0}$. Iterating this
procedure we eventually produce $D_{r+1}$ with $W_{D_{r+1}}$ a point
in $V$. Observe that in this case the divisor $D_{r+1}=D_0+\sum_1^{r+1} H_i$.
\end{proof}

Let $D$ be as in the claim, thus
$(X,D)$ is LC along $W_D$. If $r-\delta>-1$ then  we can apply
Lemma \ref{sectame}
  to produce a section of $L$ not vanishing along
$W_D$ and obtain a contradiction.

If $r-\delta=-1$ then $W_D$ is a point in $V$. Moreover, according to 
Claim  \ref{cl:clc},
in this case the divisor
$D$ is of the following type
$$D=D_0+\sum_1^{r+1} H_i,$$
with $H_i\in|L|$ generic.

Let $X_j=X\cap(\cap_1^j H_i)$, then $X_j$ is LT in a neighborhood of
$F_j:=F\cap X_j$ for any $j\leq r+1$.
This assertion is left to be proved to the
reader
as an exercise (\hint the main point to check is
normality. To do it one has to use the fact that terminal singularities are
smooth in codimension 2).

By vanishing theorem \ref{van} we have the following surjection
$$H^0(X,L)\ra H^0(X_j,L_{|X_j}),$$
for any $j\leq r$.

If the type of $f$ is not $(1,1,1)$ then $f_r:X_r\ra S$ is
birational. In particular, by standard vanishing, $Z_r\iso\Proj^1$.
So that $L_{|Z_r}$ is spanned. The idea is to extend a section of
$L_{|Z_r}$ not vanishing on $Z_{r+1}$ to a section of $L_{|X_r}$.
For details on this extension and  about the case of type $(1,1,1)$ we refer
to \cite{dpf}.
\end{proof}

\medskip
We conclude this section with an exercise which follows
easily from the main Theorem \ref {aw} and the method used in the proof
of \ref{al} (a proof can be found for instance \cite[pg 245]{RaC}).

\begin{exercise} Let $X$ be a Fano manifold of index $i(X)$;
then $i(X)\leq dim X+1$; moreover $i(X)=dimX+1$ if and
only if
$X\iso\Proj^{dimX}$ while $i(X)\geq dimX$ if and only if either
$X\iso \Proj^{dimX}$ or $X\iso\Q^{dimX}$.
\label{highindex}
\end{exercise}

\part{Biregular geometry}
\setcounter{section}{0}

Fano-Mori contractions are fundamental tools of the Minimal
Model Program; more generally they are important
in problems of classification of projective varieties.

This part is devoted to the problem of describing F-M contraction.
Unless for few results
at the very beginning we will restrict ourself to the smooth case,
that is we consider F-M contractions of smooth manifolds.

The singular case is very difficult and at the moment very few is known
only in
dimension $3$ (essentially the complete classification of small extremal
contractions on threefolds with at most terminal singularities in
the fundamental papers of Mori \cite{3flip} and of Koll\'ar-Mori
\cite{KoMoflip}).

We will give a complete classification of F-M contractions
of smooth manifolds of dimension $\leq 4$; we collect this classification
in a sequel of theorems in the first section.

We are interested in a local description
of the contraction, in a neighborhood of a given fiber;
in particular we consider a {\sl local contraction around
$F$}, $\f:X\ra Z$, as defined in \ref{defloc}.

We present many steps of the proof
of the classification; each step is important by itself and together they
represent a sort of program for classifying the F-M contractions.
In short they are the following:

1) classify all possible fibers of the F-M contractions;
we will succeed if their dimension is less or equal then two.

2) when the fiber has good singularities (locally complete intersections)
classify the possible normal bundle of these
fibers

3) describe a formal neighborhood of the possible fibers in $X$,
i.e. the local contraction around $F$.

4) find a commutative diagram of morphisms, preferably
blow-up and blow-down, which includes
$\f$ and which can help in understanding $\f$ (
a sort of factorization of $\f$, for example the flip
in the small contraction case).

\smallskip
The results contained in this part are classical for the case
$n = dimX = 2$, and they are due to the Italian school of geometry of
the beginning of the century.

In the case $n= 3$ they were proved by S. Mori in a famous paper,
see \cite{Mo}, who gave rise to the so called Mori theory.

The case $n=4$ was later considered by M. Andreatta and J.A. Wi\'sniewski,
see \cite{AW3}.
\cite{AW2} is a survey of these results on which this part is strongly
based.

\smallskip
In section \ref{sec:An} we present two theorems which characterise
some F-M contractions of a smooth projective variety in higher dimension.

\smallskip
In the last section we outline the biregular classification
of Fano manifolds of high index. These are the building block of F--M
contractions and
their knowledge is the starting point of any further investigation. Also in
this
case we will provide the known general techniques to approach the problem
via adjunction
methods, without any attempt
to be exhaustive in the classification. In particular we will not present
neither Fano--Iskovskikh
approach based on double projections, \cite{Is}, nor Mukai vector bundle
technique, \cite{Mu},
nor Ciliberto--Lopez--Miranda
deformation ideas, \cite{CLM}.

\section[F--M in dimension$\leq 4$]{Fano-Mori contractions on a smooth
$n$-fold
with $n \leq 4$}

Here we describe
all F-M contractions on smooth $n$-folds
with $n \leq 4$.
The case of dimension $4$ is the more
elaborate. Proofs are given in the next sections.

\medskip
\begin{Theorem}
\label{2-folds}
Let $X$ be a smooth projective surface
and $R \subset \overline{NE(X)}$ an extremal ray that is  $R^.K_X < 0$ and $R$ is an edge of the cone.
Then the associated contraction morphism $cont_{R} : X \raa Z$
is one of the following:
\begin{itemize}
\item{(1)} $Z$ is a smooth surface and $X$ is obtained from $Z$
by blowing-up a point; $\rho (Z) = \rho(X) - 1$.
\item{(2)} $Z$ is a smooth curve and $X$ is a minimal ruled surface over $Z$;
$\rho(X) = 2$.
\item{(3)} $Z$ is a point, $\rho(X) = 1$ and $-K_X$ is ample;
in fact $X \cong \Proj ^2$.
\end{itemize}
\label{MMPsurf}
\end{Theorem}

\medskip
\begin{Theorem}
\label{3-folds}
Let $X$ be a smooth projective $3$-fold
and $R \subset \overline{NE(X)}$ an extremal ray.
Then the associated contraction morphism $cont_{R} : X \raa Z$
is one of the following:
\begin{itemize}
\item{(B)} (Birational contractions) $dim Z = 3$, $cont_R$ is a divisorial
contraction
and there are five types of local behavior near the exceptional
divisor $E$:
\begin{itemize}
\item{B1:} $cont_R$ is the (inverse of the) blow-up of a smooth curve
in the smooth threefold $Z$.
\item{B2:} $cont_R$ contracts a smooth $\Proj ^2$ with normal
bundle $\O(-1)$; $cont_R$ is the (inverse of the) blow-up of a smooth point
in the smooth threefold $Z$.
\item{B3:} $cont_R$ contracts a smooth two dimensional quadric, $\F_0$,
with normal bundle $\O(-1)$; $cont_R$ is the (inverse of the) blow-up of
an ordinary double point in
$Z$ (locally analytically, an ordinary double point is given by the equation
$x^2 + y^2 + z^2 +w^2 = 0$).
\item{B4:} $cont_R$ contracts an irreducible singular two dimensional
quadric, $\Sc_2$,
with normal bundle $\O(-1)$; $cont_R$ is the (inverse of the) blow-up of a
point in $Z$
which is locally analytically given by the equation $x^2 + y^2 + z^2 +w^3 = 0$.
\item{B5:} $cont_R$ contracts a smooth $\Proj ^2$
with normal bundle $\O(-2)$;  $cont_R$ is the (inverse of the) blow-up of a
point in $Z$
which is locally analytically given as the quotient of $\C^3$ by the involution
$(x,y,z) \raa (-x,-y,-z)$.
\end{itemize}
\item{(C)} (Conic Bundle) $dimZ = 2$ and $cont_R$ is a fibration whose fibers
are plane conics
(general fibers are of course smooth).
\item{(D)} (del Pezzo fibration) $dimZ = 1$ and $cont_R$ is a fibration whose
general fiber is a del Pezzo surface.
\item{(F)} (Fano threefolds) $dim Z = 0$, $-K_X$ is ample, thus $X$ is a Fano
threefold, and $\rho(X) = 1$.
\end{itemize}
\end{Theorem}

As said in the introduction of this chapter the first Theorem is by
G. Castelnuovo F. Enriques and the second is by S. Mori.
Let us note that they are true actually
over any algebraically closed fields;
the surface case follows from the fact that
the Castelnuovo contraction theorem is true in any characteristic
and the threefold one was proved by Koll\'ar in \cite{KoPoc},
extending Mori's ideas.

\bigskip
The next theorem aims to give the same result
for the case $n=4$; here the situation is much more intricate
and it will take some space to be described.

The result comes from many contributions, the main ones are from
Y. Kawamata, see \cite{flip}, and from M. Andreatta and J.A. Wi\'sniewski,
see \cite{AW3} and \cite{AW4}; in the fiber case, Y. Kachi obtained 
independently
of {\cite{AW3} a similar classification of special two dimensional
fiber of a conic fibration, while in the case of birational contractions
contracting a divisor to a curve (part 3) Takagi obtained the same results
as in section 4 of \cite{AW4}.

\bigskip
\begin{Theorem}
\label{4-folds}
Let $X$ be a smooth projective $4$-fold
and $R \subset \overline{NE(X)}$ an extremal ray.
Let $\f :=  cont_{R} : X \raa Z$ be the associated contraction
morphism.
Let $F=\f^{-1}(z)$ be a (geometric) fiber of $\f$;
we will eventually shrink the morphism $\f$ around $F$, see \ref{defloc}.
Let $E$ be the exceptional locus,
where in the case $\f$ is of fiber type we mean $E = X$.

We divide the classification of these contractions depending on the
couple of numbers $(dimE, dim \f(E))$ which we will call the signature of
the contraction;
note that the pair $(4,b)$ will be given to a fiber type contraction
with $dim Z = b$ and the pairs $(a, b)$ with $b\geq a$ cannot happen.

Note also that if $dim E =3$ then $E$ is irreducible
(see \ref{divisorial}) and so is $\f(E)$,
therefore they are both of pure dimension.

For the notation adopted to describe some special two dimensional fiber
see \ref{surface}.

\smallskip\noindent
{\bf Part 0.}
The are no F-M contraction of a fourfold of signature
$(a,b)$ with $a \leq 1$ and with $a=2$ and $b=1$.

\smallskip\noindent
{\bf Part 1: Small contractions, see \cite{flip}.}
Let $\f$ be a F-M contraction of a fourfold of signature $(2,0)$.
Then $E=F\iso\Proj^{2}$ and its normal bundle is $N_{F/X} = \O(-1) \oplus
\O(-1)$.
The contraction is completely determined in an analytic
neighborhood by this data, see \ref{formaln}
and also \ref{graded}, and locally it is analytically isomorphic
to the contraction given by
(see \ref{examplev})
$$\f :{\bf V}({\mathcal E}) \ra Z = Spec (\bigoplus _{k\geq 0}
H^0(F,S^k(\mathcal E))$$
where ${\mathcal E} = \O(1)\oplus\O(1)$
and the map is associated to the evaluation of $S^k(\mathcal E )$.

In this situation the {\it flip} of $\f$ exists
and it is obtained by blowing up $E$ and then
contracting the exceptional divisor in the other direction.

\smallskip\noindent
{\bf Part 2: Birational; divisor to point.}
Let $\f$ be a F-M contraction of a fourfold of signature $(3,0)$.
Then either $E$ is $\Proj^3$, with normal bundle $\O(-a)$
and $1 \leq a \leq 3$, or a (possibly singular) three dimensional
quadric, with normal bundle $\O(-a)$
and $1 \leq a \leq 2$, or otherwise
$(E; -E_{|E})$ is a del Pezzo threefold,
that is $E$ has Gorenstein singularities,
$-E_{|E}$ is ample and $K_E = 2 E_E$
(these varieties have been classified by T. Fujita,
see \cite{Fu} and \cite{Fu3}).

\smallskip\noindent
{\bf Part 3: Birational; divisor to curve.}
Let $\f$ be a F-M contraction of a fourfold of signature $(3,1)$.
Then
\begin{itemize}

\item{(a)} $C:= \f (E)$ is a smooth curve and $\f: X \ra Z$ is the
blow-up of $Z$
along $C$.

\item{(b)} $g:= \f_{|E} : E \ra C$ is either a $\Proj^2$-bundle or a
quadric bundle.

\item{(c1)} If $E$ is a $\Proj^2$-bundle then the normal bundle of each
fiber in $X$
is either $\O(-1) \oplus \O$ or $\O(-2)\oplus \O$; in particular all fibers
of $\f$ are
reduced and with no embedded components. In the first case $Z$ is smooth
and $\f$ is
the smooth blow-up; in the second $C= Sing Z$ and $Z$ is locally isomorphic to
$S_2 \times \C$ where $S_2$ is the germ of singularity obtained by contracting
the zero section in the total space of the bundle $\O(2)$ over $\Proj^2$.

\item{(c2)} If $E$ is a quadric bundle then the general fiber is
irreducible and
isomorphic to a two dimensional, possibly singular, quadric. Isolated
special fibers can occur
and they are isomorphic either to a singular quadric or to a reduced but
reducible quadric
(i.e. union of two $\Proj^2$ intersecting along a line); in particular
there are no special fibers which are isomorphic to a double plane.
The normal bundle of each fiber is $\O(-1) \oplus \O$.
Locally $Z$ can be described as a hypersurface of $\C^5$; in the following
table
we give a list of possibilities for $Z =V(g) \subset \C^5$ according to the
described
combinations of general and special fibers.
We choose coordinates $(z_1, z_2, z_3, z_4, z_5)$ such that
$C = \{z_1= z_2= z_3= z_4= 0\} \subset \C^5$.

$$
\begin{array}{llll}
\hline
N^0\ & special\  fib. & gen.\  fib. & g = analytic\ equation\ of\ Z
\\
\hline
\\
(1)&\F_0 & \F_0 & z_1^2 + z_2^2 + z_3^2 + z_4^2
\\
(2)&\Sc_2 & \F_0 &z_1^2 + z_2^2 + z_3^2 +z_5^m z_4^2, m\geq 1
\\
(3)&\Sc_2 &\Sc_2 & z_1^2 + z_2^2 + z_3^2 +z_4^3
\\
(4)&\Proj^2 \cup\Proj^2 & \Sc_2 & z_1^2 + z_2^2 +z_3^3 + z_4^3 + z_3^2
z_5^m , m\geq1
\\
(5)&\Proj^2 \cup\Proj^2 & \F_0 &z_1^2 + z_2^2 + z_3^3 +z_4^3 +z_3^2 z_5^m +
z_3 z_4 f(z_5)
\\
  & & &+z_4^2g(z_5) \hbox{\ with\ }  z^mg(z) \not= \frac{f(z)^2}{4}
\\
\hline
\end{array}
$$
\end{itemize}

\smallskip\noindent
{\bf Part 4: Birational; divisor to surface.}
Let $\f$ be a F-M contraction of a fourfold of signature $(3,2)$.
Generically the map is described
by the part 1 of Theorem \ref{Ando}; in particular
$Z$ as well as $S:=\f(E)$ are in general smooth
and $\f$ is a
simple blow-down of the divisor $E$ to the surface $S\subset Z$.

However there can be some special two dimensional fibers $F$.
If this is the case  then the
scheme theoretic fiber structure over $F$ is trivial, that is the
ideal $\I_F$ of $F$ is equal to the inverse image of the maximal
ideal of $z$, that is $\I_F=\f^{-1}(m_z)\cdot\O_X$.

Moreover the fiber $F$ and its conormal bundle $\I_F/\I_F^2$ as well as
the singularity of $Z$ and $S$ at $z$ can be described as follows

$$
\begin{array}{llll}

F  &N^*_{F/X} &Sing Z &Sing S
\\
\hline

\Proj^2 &T(-1)\oplus\O(1)/\O
&cone\  over\  \Q^3&smooth
\\
\Proj^2 &\O^{\oplus 4}/\O(-1)^{\oplus 2}
&smooth &cone\  over\ a \\
&&&twisted \  cubic\\
Quadric &spinor\ bundle\ from\ \Q^4 &smooth&non-normal
\end{array}
$$

The quadric fiber can be singular, even reducible, and in the subsequent
table we present a refined description of its conormal bundle.
The last entry in the table provides information
about the ideal of a suitable surface $S$;
a complete description of these ideals can be find in \cite{AW4}.

$$
\begin{array}{lll}

quadric & conormal\ bundle & \I(S)\ in\ \C[[x,y,z,t]]
\\
\hline

\Proj^1 \times \Proj^1 &\O(1,0)\oplus\O(0,1)
&(xz,xt,yz,yt)
\\
quadric\ cone & 0\ra\O\ra N^*\ra \J_{line}\ra 0
  & generated\ by\ 5\ cubics
\\
\Proj^2 \cup \Proj^2 &T_{\Proj^2} (-1) \cup (\O \oplus \O(1))
& generated\ by\ 6\ quartics
\end{array}
$$

\smallskip\noindent
{\bf Part 5: Conic bundle fibration with possibly special two dimensional
fiber.}
Let $\f$ be a F-M contraction of a fourfold of signature $(4,3)$.
Then $\f$ is a fibration whose general
fibers are plane conics; generically the map is described
by the part 2 of Theorem \ref{Ando}. In particular
$Z$ is in general smooth.

However there can be some special isolated two dimensional
fibers $F$; the possibilities for $F$ are the following:
\begin{itemize}
\item $F\iso \Proj^2$ and
$N^*_{F/X}\iso \O^3/\O(-2)$ or
$T\Proj^2(-1)$. The scheme fiber structure $\tilde F$
is reduced and $Z$ is smooth at $z = \f(F)$.
\item $F$ is an irreducible quadric
and $N^*_{F/X}$ is the pullback of $T\Proj^2(-1)$ via some double covering
of $\Proj^2$. The scheme fiber structure $\tilde F$
is reduced and $Z$ is smooth at $z= \f(F)$.
\item The following other possibilities for $F$ can occur:
$$\Sc_3, \F_1,
\Proj^2\cup\Proj^2, \Proj^2\cup\F_0, \Proj^2\cup_{C_0}\F_1,$$
$$
\Proj^2\cup\Sc_2,\Proj^2\cup\Proj^2\cup\Proj^2,\Proj^2\cup_f{(\F_0)}\cup_{C_0}
\Proj^2,
$$
where any two components intersect along a line (
explicitly indicated by a
subscript, when needed),

\noindent
and the exceptional case of
$\Proj^2\bullet\Proj^2$ when the two components intersect at an isolated point.
\end{itemize}

\smallskip\noindent
{\bf Part 6: del Pezzo and Mukai fibration and Fano fourfolds.}
Let $\f$ be a F-M contraction of a fourfold of signature $(4,d)$,
with $d\leq 2$.
Then  $\f$ is an equidimensional fibration over $Z$. If $d= 2$ the
general fiber is a del Pezzo surface, if $d = 1$ then the
general fiber is a Mukai variety, while if $d= 0$, $-K_X$ is ample, thus
$X$ is a Fano
fourfold, and $\rho(X) = 1$.
\end{Theorem}

\bigskip
Let us add some remarks at this long theorem.

\begin{remark} The case $(3,0)$ is not complete,
in fact it contains many non existing cases.
More precisely let $E$ be a del Pezzo threefold,
i.e. $-K_E = 2 {\mathcal L}$ with ${\mathcal L}$ ample.

If $E$ is smooth then
one can easily construct
a F-M contraction of a smooth fourfold of signature $(3,0)$
and exceptional divisor $E$
by taking (see \ref{examplev})
$$\f :{\bf V}({\mathcal L}) \ra Z = Spec (\bigoplus _{k\geq 0}
H^0(F,({\mathcal L^k}))$$
where $\E = - E_{|E}$
and the map is associated to the evaluation of ${\mathcal L^k}$.

However it is conjecture that there are no F-M contraction of a smooth fourfold
with a non normal exceptional divisor $E$; in section 3 of \cite{Fu3} 
this case is
discussed deeply and a lot of limitation on $E$ are given
(see \ref{conjfu} and the following discussion).

If $E$ is singular but with normal singularities then a list of 
possible $E$ was given in \cite{Be}
but this list contains many redundant case.

\smallskip
The case $(4,3)$ is also not complete. In particular
we have examples of appropriate 2 dimensional fibers
except for the cases $\Proj^2\cup\Sc_2$, $\Proj^2\cup\Proj^2\cup\Proj^2$ and
$\Proj^2\cup_f{(\F_0)}\cup_{C_0}\Proj^2$; we believe these cases
cannot occur.
\end{remark}

\section[F--M with small fibers]{Fano-Mori contractions on a smooth $n$-fold
with fibers of small dimension}
\label{sec:An}

In this section we present two theorems which characterise
some F-M contractions of a smooth projective variety in higher dimension.

The first is due to T. Ando and it deals with F-M contractions
with one dimensional fibers.

\begin{Theorem} (\cite{An}) Let $\f :X \ra Z$ be a (local) Fano-Mori 
contraction of
a smooth variety $X$ of dimension $n$ around a fixed fiber $F=\f^{-1}(z)$
such that $dim F = 1$.
\item{(1)} If $\f$ is birational then $F$ is irreducible, $F \iso \Proj^1$,
$-K_X\cdot F=1$ and its normal bundle is $N_{F/X}=\O(-1) \oplus\O^{(n-2)}$.
The target $Z$ is smooth and $\f$ is a blow-up of a smooth codimension 2
subvariety of $Z$.
\item{(2)} If $\f$ is of fiber type then
$Z$ is smooth and $\f$ is a flat conic bundle. In particular one of the
following is true:
\item{(i)} $F$ is a smooth $\Proj^1$ and $-K_X\cdot F=2$,
$N_{F/X}\iso\O^{(n-1)}$;
\item{(ii)} $F=C_1\cup C_2$ is a union of two smooth rational curves meeting
at one point and  $-K_X\cdot C_i=1$, $(N_{F/X})_{|C_i}\iso \O^{(n-1)}$,
$N_{C_i/X}\iso \O^{(n-2)}\oplus\O(-1)$ for $i=1,\ 2$;
\item{(iii)} $F$ is a smooth $\Proj^1$, $-K_X\cdot F=1$ and the fiber
structure $\tilde F$ on $F$ is of multiplicity 2 (a non reduced conic);
the normal bundle of $\tilde F$ is trivial while $N_{F/X}$ is either
$\O(1)\oplus\O(-1)^{(2)}\oplus\O^{(n-4)}$ or
$\O(1)\oplus\O(-2)\oplus\O^{(n-3)}$
depending on whether the discriminant locus of the conic bundle is smooth
at $z$ or not.
\label{Ando}
\end{Theorem}

The above theorem was
generalised to the case of a variety $X$ with
terminal Gorenstein singularities
by Mori and Koll\'ar (see \cite[4.9 and 4.10.1]{KoMoflip}) for $n\geq 3$.

The case of an extremal contraction of a 3-fold $X$ with terminal non
Gorenstein
singularities is much more difficult; this was discussed in the celebrated
paper
of Mori \cite{3flip} and in \cite{KoMoflip}.

\smallskip
The next theorem is a generalisation of the above theorem of Ando
in the frame of {\sl adjunction theory of projective varieties}
a very classical theory (see \cite{CE}), which was
revitalized and improved in modern time by A.~J.~Sommese and his
school (see \cite{BS}).

One of the goals of this theory is to describe varieties $X$
polarized by an ample line bundle $L$ by means of the Fano-Mori
contraction supported by $K_X + rL$ where $r$ is the nef value of the
pair $(X,L)$.  If $X$ is smooth and $r \geq (n-2)$ then this goal is
achieved and we refer the reader to the book \cite{BS} for an
overview of the theory, see \cite{A1}, \cite{A2} and \cite{Meadj} for the
singular case.

The next theorem, proved in \cite{AW},  shows that it is also 
achieved when the nef value is large
with respect
to the dimension of fibers of $\f$.

\begin{Theorem} [\cite{AW}] Let $\f: X \ra Z$ be a (local) Fano-Mori
contraction
of a smooth variety $X$ and let $F=\f^{-1}(z)$ be a fiber.
Assume that $\f$ is supported by $K_X+rL$, with $L$ a $\f$-ample line
bundle on $X$.
\begin{itemize}
\item{(1)} If $dim F \leq (r-1)$ then  $Z$ is smooth at $z$ and $\f$
is a projective bundle in a neighborhood of $F$.
\item{(2)} If $dimF=r$ then, after possible shrinking of $Z$ and restricting
$\f$ to a neighborhood of $F$, $Z$ is smooth and
\begin{itemize}
\item{(i)} if $\f$ is birational then $\f$ blows a
smooth divisor $E\supset X$
to a smooth codimension $r-1$ subvariety $S\supset Z$,
\item{(ii)} if $\f$ is of fiber type and $dimZ=dimX-r$ then $\f$ is a
quadric bundle,
\item{(iii)} if it is of fiber type and $dimZ = dimX-r+1$ then $r\leq
dimX/2$, $F =\Proj^r$ and the general fiber is $\Proj^{(r-1)}$.
\end{itemize}
\end{itemize}
\end{Theorem}

\medskip
The basic steps of the proofs of these theorems are worked out in the next
sections, together with the proofs of the results in
the previous section.

\section{The fibers of a Fano-Mori contraction}

In this section we will try to give more informations on
the possible fibers of the F-M contractions. In particular
we will classify all possible fibers of dimension less than or equal to two

\subsection{Using the vanishing theorem}

We want to show how the vanishing theorem implies vanishing results on the
fiber. Subsequently we show how these results, via the computation
of the Hilbert polynomial of the (normalization) of the fiber, imply
a bound on the dimension of the fiber.

\smallskip
The proof of the following
proposition can be found in \cite{Mo}, 3.20, 3.25.1,
\cite{Fu}, 11.3, \cite{An} and \cite{AW3}, 1.2.1.

\begin{Proposition} [Vanishing of the highest cohomology]

Let $\f: X \ra Z$ be a local F-M contraction
around $F$ supported by $K_X+rL$ (see \ref{defloc}).
Let $F'$ be a subscheme of $X$
whose support is contained in the fiber
$F$ of $\f$, so that $\f(F') = z$.
If either $t>-r$ or $t =-r$ and $dimF>dimX-dimZ$ then
$$H^{dimF}(F',tL_{|F'})=0.$$
\label{van-fib1}
\end{Proposition}

\begin{Proposition} In the assumption of the above proposition
let also $X'\in |L|$ be the zero locus of
a non-trivial section of $L$. Then we have
$$H^{dim (F \cap X')}(F' \cap X',tL_{|F'\cap X'})=0$$
if either $t > -r +1$
or $t=- r +1$ and $dim (F \cap X') \geq dimX-dimZ$.
\label{van-fib2}
\end{Proposition}

\begin{proof} We like to give here a proof of the second proposition, the
proof of the first is similar.
Note that $H^i(X,tL) = 0 $ for $i > 0$ and $t > -r$
by the theorem (\ref{van}); moreover
we also have $H^i(X,tL) = 0 $ for
$t= -r$ and $i > dimX -dim Z$, for the
so called Grauert-Riemenschneider-Koll\'ar vanishing theorem
(see \cite{KMM}, theorems 1.2.4 and 1.2.7).
Thus from the exact sequence
$$0 \ra -L \ra \O_{X} \ra \O_{X'} \ra 0$$
tensorised by $tL$ we also have
$H^i(X', tL_{X'})=0$ for $i > 0$ and $t > -r+1$ or
$t= -r+1$ and $i > dimX -dim Z$.

Now let $\I_{F'\cap X'}$ be the ideal of $F'\cap X'$ in $X'$ and consider
the sequence
$$0\ra \I_{F'\cap X'}\otimes tL \ra \O_{X'}\otimes tL \ra \O_{F'\cap
X'}\otimes tL \ra 0.$$
Take the associated long exact sequence.
Since
$H^i(X',\I_{F'\cap X'}\otimes tL)_z=0$ for $i> q:= dim F \cap X'$,
the map $H^q(tL_{X'})\ra H^q(F'\cap X',tL_{F'})$
is surjective and the proposition follows from
what we have observed at the beginning.
\end{proof}

The following result is a direct consequence of the above Proposition;
it was proved by T. Fujita, see \cite{Fu0}, following arguments of S. Mori
and T. Ando.

\begin{Theorem} Let $\f: X \ra Z$ be a local F-M contraction
around $F$ supported by $K_X+rL$.
Then
$$dim F \geq (r-1)$$
  and if
$dimF > dimX-dim Z$ then
$$dimF \geq \lfloor r\rfloor$$

\label{boundfiber}
\end{Theorem}

\begin{proof} 
Let $S$ be a component of a fiber $F$ of
dimension $s$ and let
$g: W \ra S$ its desingularization.
By the above proposition \ref{van-fib1} and the Leray spectral sequence
for $g$, exactly as in Lemma 2.4 of \cite{Fu0}, we get
$$H^s(W,g^*(tL))= 0$$
if $ t  > -r$ or $ t = -r$ and  $dimF > dimX-dim W$.
\par
On the other hand, since $g^*(L)$ is nef and big on $W$, by the
Kawamata-Viehweg
vanishing theorem we have that $H^i(W,g^*(tL)) = 0$ for $t \geq -r$ and
$0< i < s$. Moreover, since $L$ is ample,
we have also that $H^0 (W,g^*(tL)) = 0$ for $t<0$.
\par
Consider now the Hilbert polynomial $\chi(t) := \chi(W,g^*(tL))$;
it is a polynomial in $t$ of degree equal to $dim W = dimS$.
By what proved above $\chi$ is zero
for all integers $t$ such that $0 > t  > -r$; if $dimF > dimX-dimZ$
and $r$ is an integer then $\chi$ is zero also for $t = - r$.
The inequalities follows then immediately since
$deg\chi \geq $ number of its zeros.

\end{proof}

\subsubsection{Exercises-Examples}

\begin{exercise} (see \cite{Fu0}) Let $(X,L)$ be a polarized variety; 
$K_X+nL$ is nef
except when $(X,L) = (\Proj^n, \O(1))$. Also, if $n \geq 3$,  $K_X+(n-1)L$
is nef unless
$(X,L)$ is one of the following: $ (\Proj^n, \O(1))$,  $X$ is the quadric in
$\Proj^n$ and $L$ is a hyperplane section, $X$ is the projectivization of a
rank $n$ vector bundle over a smooth curve $A$ and $L= \O(1)$.
\end{exercise}

\noindent
This result was proved in the paper \cite{Fu0} with the use of
the above
Theorems; a completely different proof, which make use of the deformation
of rational curves
and which works in all characteristic, has been recently given in \cite{KK}.

\subsection{Existence of a ladder for a fiber of a F-M contraction;
horizontal slicing}
\label{chap:ladder}

Here we will develop in more details the method outlined
in the introduction of section \ref{relgood}. With the use of the theorem
\ref{aw}
we will study the fibers of a
F-M contraction inductively.
More precisely if $F$ is a fiber of a Fano-Mori contraction of sufficiently
high dimension (i.e. with
''small'' difficulty) then we can construct a ladder for the pair $(F,L_F)$
and prove that $\Delta(F,L) = 0$.

In order to do this we first start with a
Bertini type theorem.

\begin{PD}[\cite{AW}, Lemma 2.6 , \cite{dpf}, Lemma 1.3](Horizontal slicing)
Let $\f:X\ra S$ be
a local contraction around $\{F\}$
supported by $K_X+rL$. Let $H_i\in |L|$
generic divisors and $X_k=\cap_1^k H_i$, a scheme theoretic intersection;
assume that
$dim X_k = n-k\ (> 0)$ and that $(r-k) \geq 0$;
note that since $X_k$ is a complete intersection
it is $\Q$-Gorenstein, i.e. $K_{X_k}$ is $\Q$-Cartier.
\begin{itemize}
\item[i)] Let
$\f_{|X_k}=g\circ \f_k$ the Stein factorisation of
$\f_{|X_k}:X_k\ra S$; then $\f_k:X_k\ra S_k$ is a morphism
with connected fiber, around
$\{F\cap(\cap^k_1 H_i)\}$, supported by
$K_{X_k}+(r-k)L_{|X_k}$
and $S_k$ is affine. In particular if $X_k$ is normal then $\f_k$
is a local contraction.
\end{itemize}
Assume that $X$ has LT singularities and, if
$\epsilon$ is a sufficiently small positive rational number,
that $r\geq \epsilon\gamma(\f)$ and $k\leq r+1-\epsilon\gamma(\f)$.
\begin{itemize}
\item[ii)] Outside of $Bsl|L|$ $X_k$ has singularities
which are of the same type of the ones of $X$ and
any section of $L$ on $X_k$ extends
to a section of $L$ on $X$.
\label{horizontal}
\end{itemize}
\end{PD}

\begin{proof}
See \cite{dpf}.

$(i)$ is just Stein factorisation  (see \cite[III.11.5]{Ha})
and adjunction formula,
once noticed that $f_{|X_k}(X_k)=Spec(H^0(H,\O_{X_k}))$ and
that there is
a morphism $S_k\ra S$ induced by the ring morphism
$H^0(X,\O_X) \ra H^0(X_k,\O_{X_k})$.

For $ii)$ the first statement is just Bertini Theorem, while
for the latter consider the exact sequences
\begin{eqnarray*}
  0\ra \O_{X_i}(-L)\ra \O_{X_i}\ra \O_{X_{i+1}}\ra 0\\
  0\ra \O_{X_i}\ra \O_{X_i}(L)\ra \O_{X_{i+1}}(L)\ra 0.
\end{eqnarray*}
Thus to prove the assert it is enough to prove that $H^1(X_i,\O_{X_i})=0$,
for $i\leq r-\epsilon(dimX-dimS)$. But this is equivalent, using inductively
the first sequence tensored,
to $H^1(X,-iL)=0$, for $i\leq r-\epsilon(dimX-dimS)$,
which follows from the vanishing theorem \ref{van}.
\end{proof}

\begin{Theorem}  Let $\f: X \ra Z$ be a local F-M contraction
around $F$ supported by $K_X+rL$ of type
$(d(\f),\gamma(\f),\Phi(\f))$.
Let $S$ be any component of $F_{red}$ and $\epsilon$
a sufficiently small positive rational number.
If $\Phi(\f) \leq 1-\epsilon
\gamma(\f)$
or the type of $\f$ is $(d,1,1)$, with $d\leq 0$,
then $\Delta(S,L_S) =0$.

If $\Phi(\f) \leq -\epsilon \gamma(\f)$ or if the type of $\f$ is $(-1,1,0)$,
then $F_{red}$ is irreducible and isomorphic to
$\Proj^{dimF}$.
\label{delta}
\end{Theorem}

\begin{proof}
We present the proof which is contained in
the paper (\cite{dpf}, Theorem 2.17); this is a slight
generalisation of the proof in (\cite{AW3}, Proposition 4.2.1) and
(\cite{AW2}, Theorem 1.10)

To prove the first part of the Theorem
let $S$ be an irreducible component of $F_{red}$ of dimension $s$
and let $\delta := L^s.S$.
We have to prove that $h^0(S,L_S) \geq \delta +r +1$.
This follows obviously if we will prove that there are at least 
$\delta +r +1$ independent
sections
of $H^0(X,L)$ not vanishing identically on $S$.

By the propositions \ref{horizontal} and \ref{aw} we reduce
to the case of a contraction $\f: X \ra T$  with one
dimensional fiber $F$.
Then, by assumption, we can use again
propositions \ref{horizontal} and \ref{aw} and
go one step further with a section $H \in |L|$;
$\f_{|H }:H \ra T$ is finite and by \ref{horizontal} all
section of $|L_{|H}|$ extend to sections of $|L|$ proving the thesis.

Finally assume that $\Phi(\f) \leq -\epsilon \gamma(\f)$ and
assume, by contradiction, that the fiber has (at least) two
irreducible components intersecting in a subvariety of dimension $t
\leq (r-1)$. By the base point freeness of $L$, we can choose $t + 1$
sections of $L$ intersecting transversally in a variety with log
terminal singularities and meeting the two irreducible components not
in their intersection.  By construction the map $\f$ restricted to this
variety has non connected fibers and this is in contradiction with
\ref{horizontal}.
Similarly one can prove that $L^r\cdot F=1$
(we can slice to points and still have the connectedness, but then we
must have only one point...)
and thus that $F = \Proj^{dimF}$.
\end{proof}

An immediate corollary in the case of two dimensional fiber is the
following.

\begin{Corollary} Let $F$ be a two dimensional fiber of a F-M contraction
$\f : X \ra Z$ of a Gorenstein variety $X$
and let $F'$ be any component.
Assume that $\f$ is birational or that the general non
trivial fiber has dimension $1$. Then $F'$ is normal and
the pair $(F',L_{|F'})$ has sectional (and thus Fujita $\Delta$)-genus $0$
and therefore (see \ref{surface}) it is among the following:
\item{(1)} $(\Proj^2,\O(e))$, with $e = 1,2$,
\item{(2)} $(\F_r, C_0 + kf)$ with $k \geq r+1$, $r\geq 0$,
\item{(3)} $(\Sc_r,\O_{\Sc_r}(1))$ with $r\geq 2$.

Moreover $F$ is Cohen-Macaulay unless the zero locus of
a general section in $|L_F|$ is disconnected.
\label{deltasup}
\end{Corollary}
We will see however in the next subsection that not all the possibilities
can occur
if the domain $X$ is smooth (or has very good singularities).
Another type of
argument is needed to get rid of some cases.

To conclude the section we will mention another Bertini
type theorem which has to do with the sections of the supporting divisor
of the F-M contraction.

\begin{PD}[\cite{AW}](Vertical slicing) Let $\f:X\ra S$ be a
local contraction
supported by $K_X+rL$, with $r\geq-1+\epsilon\gamma(\f)$
and $\epsilon$ a sufficiently small positive rational number.
Assume that $X$ has LT singularities and let
$h$ be a general function on $S$. Let $X_h=\f^*(h)$ then the singularities
of $X_h$ are not worse than these of $X$ and any section of $L$ on
$X_h$ extends to $X$.
\label{vertical}
\end{PD}

\medskip
\subsubsection{Related topics and further results}

\begin{exercise}[Lifting a contraction]
Let $X$ be a smooth complex projective variety
of dimension $n$ and $L$ be an ample line bundle
with a section $D\in |L|$ with good singularities (smooth, KLT).
(More generally let $\E$ be an ample vector bundle of
rank $r$ on $X$ such that there exists
a section $s \in \Gamma(\E)$ whose zero locus, $D= (s=0)$, is a smooth
submanifold  of
the expected dimension $\dim D= \dim X -r = n-r$.)

A classical and natural problem is to ascend the
geometric properties of $D$ to get informations on the geometry of $X$;
a very good account on this problem can be found in \cite[Chapter 5]{BS}.
In \cite{AO1} and \cite{AO3} the problem was considered from the point of view
of Mori theory, posing the following question: assume that $D$ is not
minimal, i.e.
$Z$ has at least one extremal ray in the negative part of the Mori cone;
does this ray (or the associated extremal contraction) determine a ray (or a
contraction) in $X$,
and if so, does this new ray determine the structure of $X$?\par

For instance assume that $D$ is $\Proj^s$ or a scroll.
\end{exercise}

\begin{exercise}[Construct F-M contractions] Find a local Fano-Mori contraction
around $F$ supported by $K_X+rL$ of type
$(d(\f),\gamma(\f),\Phi(\f))$
with $\Phi(\f) > 1-\epsilon \gamma(\f)$
and for which $F$ is not Cohen Macaulay
or in general   non normal.
Can you find such an example with $X$ smooth?
(the examples (1.18) in \cite{AW2} and
(3.6) in \cite{AW3}).
\end{exercise}

However there is the following
\begin{conjecture} Let $\f:X \ra Z$ be a Fano-Mori contraction
of a manifold of dimension $\leq 4$.
Then all the fibers are normal, with
the exception in (3.6) of \cite{AW3}.
\label{conjfu}
\end{conjecture}
Note that the conjecture, after Mori's and  Andreatta-Wi\'sniewski's 
work, is open only
for the case in which $\f:X \ra Z$ is a birational Fano-Mori contraction
of a manifold of dimension $4$ which contracts an irreducible divisor $E$
to a point. Moreover, by the paper \cite{Fu3}, $E$, which is a del Pezzo
threefold,
is (possibly) non normal only if $(-K_E)^3 = 7$, $Sing(E) \cong {\mathbb P}^2$,
the normalization of $E$ is ${\mathbb P}(\O_{{\mathbb P}^1}(1) \oplus
\O_{{\mathbb P}^1}(1) \oplus\O_{{\mathbb P}^1}(5))$ and some other conditions.

\subsection{Rational curves on the fiber of F-M contractions}
\label{sec:ratcur}

In this subsection we use another
fundamental feature of a Fano-Mori contraction
in order to complete the classification of possible two dimensional fibers:
namely the existence of
rational curves in its fibers.

\subsubsection{General Facts.}
We have in fact the
following existence theorem due to Mori \cite{HFconj}, \cite{Mo} in the
smooth case
and extended by Kawamata to the log terminal case in \cite{Ka}.
We recall
that a rational curve is a curve whose normalization is $\Proj^1$.
(Although we work over $\C$, we would like to note that if $X$ is smooth
then the existence theorem is true also in positive characteristic; this
concerns also the subsequent results obtained via deformation methods.)

\begin{Theorem}[(Existence of rational curves)]
Let $\f:X\ra Z$ be a Fano-Mori contraction of a variety with log terminal
singularities.  Then the exceptional locus of $\f$ is covered by
rational curves contracted by $\f$.
\end{Theorem}

In this section we study deformations of rational curves
following ideas started with the paper of Mori \cite{HFconj}.  We
discuss only some of the results,
concentrating on the case of smooth $X$. We refer the reader to the
book of Koll\'ar \cite{RaC} for general results concerning deformation
of curves.  The following result is from \cite[II.1.14]{RaC}.

\begin{Theorem} Let $C$ be a (possibly reducible) connected curve
such that $H^1(C,\O_C)=0$ and assume that $C$ is smoothable (see
\cite[II.1.10]{RaC}, for the definition;
an example of smoothable curve is a tree of smooth rational curves,
i.e. $C=\cup_i R_i$
where: (i) any $R_i$ is a smooth rational curve
(ii) $R_i$ intersects $\sum_{j = 1}^{i-1} R_j$ in
a single point which is an ordinary node of $C$,
see \cite[II.1.12]{RaC}).
Suppose that $f:C\ra X$ is an
immersion of $C$ into a smooth variety $X$. Then any component of the
Hilbert scheme containing $f(C)$ has dimension
at least $-K_X\cdot C + (n-3)$.
\label{hilb}
\end{Theorem}

The above result has several different versions. For example, Mori
\cite{HFconj}
proved a version of it for maps of rational curves with fixed points.
An important part of the Mori's proof of the existence of rational
curves is a technique of deforming rational curves with a fixed
0-dimensional subscheme (to ``bend'' these curves) in order to
produce rational curves of lower degree with respect to a fixed ample
divisor (to ``break'' them). In short: if a rational curve can be
deformed inside $X$ with two points fixed then it has to break.

Mori's bend-and-break technique was used by Ionescu and
Wi\'sniewski (see \cite{Io}, 0.4, and \cite{Wi}, 1.1)
to prove a bound on the dimension of the fiber. The reader can
compare this bound with the one obtained in Theorem \ref{boundfiber}.

\begin{Theorem} Let $\f :X \ra Z$ be a Fano-Mori contraction
of an extremal ray $R$
of a smooth variety $X$. Let $E$ be the exceptional locus of $\f$
(if $\f$ is of fiber type then $E:=X$) and
let $S$ be an irreducible component of a (non trivial) fiber $F$.
Let $l$ = min $\{ -K_X\cdot C$: $C$ is a rational curve in $S\}$.

Then $dim S + dim E \geq dimX + l -1$.
\label{IoWi}
\end{Theorem}

\begin{Corollary} Let $\f :X \ra Z$ be a Fano-Mori contraction
of a smooth variety $X$ supported by $K_X+rL$. Let $E$ be the exceptional
locus of $\f$
and let $S$ be an irreducible component of a (non trivial) fiber $F$.

Then $dim S + dim E \geq dimX + r -1$.
\label{ineq}
\end{Corollary}

\begin{Proposition} \cite[Lemma (1.1)]{ABW}. In the hypothesis
of the above corollary, if the
equality holds for an irreducible component
then the normalization of $S$ is $\Proj^s$.
\end{Proposition}

\smallskip
Thus one can propose the following conjecture
(it was actually posed in \cite{AW2})

\begin{conjecture}
In the hypothesis of Theorem \ref{IoWi}, if equality holds
for an irreducible component $S$ then its normalization
is isomorphic to $\Proj^s$.
\end{conjecture}

A step toward the conjecture was given by the following
theorem proved in \cite{AO2}.

\begin{Theorem} If the contraction $\f:X \ra Z$ is divisorial then the
conjecture
is true and $\f$ is actually a smooth blow-up (i.e.of a smooth submanifold
of $Z$
which is also smooth).
\end{Theorem}

The following result, which was a long lasting
conjecture, has been recently proved; it is not difficult to show that
this proves a part of the above conjecture.

\begin{Theorem} \cite{CMS} \cite{Keb}
If $X$ is a smooth projective variety of dimension $n$ such that
$K_X\cdot C\leq -n-1$ for any complete curve $C\subset X$
then $X\iso\Proj^n$.
\end{Theorem}

In the Cho-Miyaoka-Shepherd Barron's preprint a more general version
is stated; in particular the variety $X$ can have normal singularities.
This version should imply the  above conjecture.

Let us notice that the last theorem is a very nice generalization of
the famous theorem of S. Mori, i.e. the proof of the Hartshorne 
Frenkel conjecture.

\begin{Theorem} \cite{HFconj}
If $X$ is a smooth variety with ample tangent bundle
then $X\iso\Proj^n$.
\end{Theorem}

In a slight different direction it has recently been proved also the following
generalization of Mori's theorem.

\begin{Theorem} \cite{AW5}
If $X$ is a smooth variety which has an ample locally free
subsheaf of the tangent bundle then $X\iso\Proj^n$.
\end{Theorem}

\medskip
\subsubsection{Rational curves on fibers of a F-M contraction of dimension
$\leq 2$}
Now we work out a complete classification
of fiber $F$ of dimension $\leq 2$ of a F-M contraction of a smooth variety
$X$.

\begin{Lemma} If a fiber $F$ of a Fano-Mori contraction of a smooth
$n$-fold $X$ contains a component of dimension 1
then $F$ is of pure dimension
$1$ and $-K_X\cdot F\leq 2$. In particular
$F$ is a line or a conic (with respect to the relative very ample
line bundle $ -K_X$), the last possibly reducible or non reduced.
If the contraction is birational then $F$ is a line.
\label{onefib}
\end{Lemma}

\begin{proof} Let $F'$ be a 1-dimensional component of $F$
($F'$ is a rational curve because of \ref{delta}). Then, by \ref{hilb}
$dim_{[F']}HilbX\geq -K_X\cdot F'+(n-3)$ and therefore small
deformation of $F'$ sweep out at least a divisor. More precisely:
taking a small analytic neighborhood of $[F']$ in $Hilb$ and the
incidence variety of curves we can produce an analytic subvariety
$E\subset X$ which is proper over $Z$ such that $F\cap E=F'$ and
$dimE\geq n-1$. This implies that all components of $F$ meeting $F'$
are of dimension 1 and by connectedness of $F$ we see that $F$ is of
pure dimension 1. The bound on the degree can be obtained similarly,
(note that because of the base point free theorem $-K_X$ is $\f$-very
ample so
that one can apply \ref{hilb} to a curve consisting of two components).
\end{proof}

\bigskip
Suppose now
that $\f:X\ra Z$ is a local Fano-Mori contraction of a smooth variety and
$F$ is
an {\sl isolated} fiber of $\f$ of dimension $\geq 2$;
isolated means that all the neighboring fibers are of dimension $\leq 1$.
That is, because $Z$
is affine, we can assume that all the fibers of $\f$ except $F$ are of
dimension $\leq 1$.

Note that by the base point freeness
theorem $L:=-K_X$ is $\f$-very ample (see \ref{aw}; the theorem states
only the relative base point freeness of $L$, but as noticed in \cite{AW3},
Proposition 1.3.4,
after possibly shrinking the affine variety $Z$, the same proof yields the
relative
very ampleness of $L$).

By Lemma \ref{onefib} all
1-dimensional fibers of $\f$ are of degree 1 (lines), or $\leq 2$ 
(conics), with
respect to $-K_X$, if $\f $ is birational or of fiber type,
respectively.

Let now $C\subset F$ be a rational curve or an immersed
image of a smoothable curve of genus 0. If the degree of $C$ with
respect to $-K_X$ is bigger than that of 1-dimensional fibers of
$\f$, then deformations of $C$ in $X$ {\bf must remain inside $F$} 
which, in view of
\ref{hilb}, provides us with the following useful observation:

\begin{Lemma} In the above situation
$$dim_{[C]} Hilb(F)\geq -K_X\cdot C+(n-3).$$
\label{hilbfib}
\end{Lemma}

Let us explain why this simple observation is useful for understanding
the structure of the fiber $F$.
Lemma \ref{hilbfib} can in fact be used to rule out
many redundant case in the list \ref{deltasup} of possible components of $F$.
We note that the very ampleness of $L = -K_X$ as well as
the precise description of the components of the fiber in \ref{deltasup}
allow us
to choose properly the curve which satisfies
the assumptions in \ref{hilbfib}.

We will give just an example
which explains our argument.
All possible cases are discussed in details in
\cite{AW3}, section 4.

Suppose that $S$ is a component of $F$ and $S\iso\Sc_r$
where $r\geq 3$. Then as the curve $C$ we take the union of general
$r+1$ lines passing through the vertex of $\Sc_r$ (the lines are general so
that none of them is contained in any other component of $F$).
Then Lemma \ref{hilbfib} implies that $\Sc_r$ can not be a component of $F$,
for $r\geq 3$ if $\f$ is birational and for $r\geq 4$ if $\f$ is of
fiber type.

\smallskip
Also this way, using \ref{hilbfib}, for a reducible fiber $F$ we can
limit the possible combination of irreducible components of $F$.
To show how let us consider the following
situation.

\begin{Lemma}
Let $F$ be an isolated fiber of dimension $\geq 2$ of a
birational contraction $\f: X\ra Z$ of a smooth $n$-fold $X$.
Suppose that the exceptional locus of
$\f$ is covered by rational curves which are lines with respect to $-K_X$.
If there exists a nontrivial decomposition $F=F_1\cup F_2$
then $F_1\cap F_2$ does not contain 0-dimensional components.
\label{interfib}
\end{Lemma}

\begin{proof}
Let $x\in F_1\cap F_2$ be an isolated point
of the intersection. Since $X$ is smooth
$dim_xF_1+dim_xF_2\leq n$.  For $i=1,\ 2$ let $C_i\subset F_i$ be
a line containing $x$.
The variety parametrising deformations
of $C_i$ inside of $F_i$ with $x$ fixed is of dimension $\leq
dim_xF_i-1$.
Indeed, take a point $y\in F_i$, then by the Bend and Break
argument of Mori, see the section II.5 of \cite{RaC},
there is only one curve of the family passing through both $x$ and $y$
(i.e through two distinct points passes only one line,
with respect to any ample line bundle).

Let us take $C=C_1\cup C_2$.
Then
$$dim_{[C]} Hilb(F)\leq dim_xF_1+dim_xF_2-2\leq n-2$$ and because
$-K_X\cdot C = 2$ we arrive to the contradiction with \ref{hilbfib}.
\end{proof}

\begin{remark}  Let us note that the above conclusion of
\ref{interfib} is no longer true if we do not
assume that $\f$ is birational, see \cite[Example (2.11.2)]{AW2}.
\end{remark}

\smallskip
With a combination of the above arguments, all based on the Lemma
\ref{hilbfib}, and through a long list of cases, in the section 4 of 
\cite{AW3},
the following has been proved.

\begin{Proposition}[ \cite{AW3}, Sect.~4)]
Let $\f: X\ra Z$ be a Fano-Mori contraction
of a smooth $n$-fold $X$ with an isolated 2-dimensional fiber $F$; let
$L=-K_X$.

If $\f$ is birational we have the following possibilities
for the pair $(F, L_F)$

$$
\begin{array}{lll}
n\geq 5&  n = 4 \ \ & n=3\\

\hline

\\
(\Proj^2,\O(1))&(\Proj^2,\O(1)) & (\Proj^2, \O(1))
\\
(n\leq 6) &(\F_0,C_0+f) &(\Proj^2,\O(2))
\\
&(\Sc_2,\O_{\Sc_2}(1)) & (\Sc_2,\O_{\Sc_2}(1))
\\
& (\Proj^2\cup \Proj^2, \O(1)) & (\F_0 , C_0 + f)
\\
& & (\F_1, C_0 + 2f)
\\
& & \Proj^2 \cup_{C_0} \F_2,
\\
& & \hbox {with\ } L_{|\Proj^2} = \O(1),
L_{|\F_2} = C_0 + 3f

\end{array}
$$

If $\f$ is of fiber type and $L$ is $\f$-spanned then
we have the following possibilities for the pair $(F, L_F)$

$$
\begin{array}{llll}
n\geq 5 &n = 4, irreducible\ \ \ \   &n=4, reducible &\ \  n=3
\\

\hline

\\
(\Proj^2,\O(1)) & (\Proj^2,\O(1)) & \Proj^2\cup\Proj^2 & (\F_0,
C_0 + 2f)
\\
(n\leq 7) & (\Proj^2,\O(2)) & \Proj^2\bullet\Proj^2 &
\F_0 \cup \F_1,
\\
(\F_0,C_0+f)& (\Sc_2,\O(1)) &\Proj^2\cup\F_0
& L_{\F_0}=C_0+f
\\
(n=5)& (\Sc_3,\O(1)) &\Proj^2\cup_{C_0}\F_1 &L_{\F_1}=C_0+ 2f
\\
(\Proj^2\cup \Proj^2, \O(1))&(\F_1,C_0+2f) & \Proj^2\cup\Sc_2
\\
(n=5)&(\F_0,C_0+f) & \Proj^2\cup\Proj^2\cup\Proj^2
\\
&  & \Proj^2\cup_{C_0}\F_0\cup_{f}\Proj^2

\end{array}$$

\medskip\noindent
In the above list the components of reducible fibers have a common line
(in some cases we point out which line is it) with the unique exception
of two $\Proj^2$-s which meet at a point --- we denote this union by $\bullet$.
(We suppress the description of $L$ whenever it is clear.)
\label{twofiber}
\end{Proposition}

\begin{remark}
Let us say again that for almost all the above possibilities we can
construct examples with appropriate isolated 2-dimensional fiber, see
Section 3 of \cite{AW3}.  However, there are some
exceptions
for which we were unable to construct examples and we do not expect
that all of them exist.  This concerns only fiber type contractions
and reducible
fibers: $\Proj^2\cup\Proj^2$ for $n =5$ and  $\Proj^2\cup\Sc_2$,
$\Proj^2\cup\Proj^2\cup\Proj^2$ and $\Proj^2\cup\F_0\cup\Proj^2$ for $n=4$.
\end{remark}

\section[Normal bundle of a fiber]{The description of the normal bundle of
a fiber of a F-M contraction}

In order to describe a contraction locally, after having determined the fiber,
one has to find out what are the possible normal bundles of these fibers;
of course when this is possible, that is when the fiber is a local
complete intersection in $X$.

This can be considered as a second order
type problem and it is very hard compared to
the determination of the fiber. If the fiber is one dimensional
it was considered by S. Mori in the case of $3$-folds
and by T. Ando in general. The case of two dimensional fiber is one
of the main achievement of the paper \cite{AW3}.

\smallskip
If the fiber is a divisor in $X$
its normal bundle is already given by adjunction formula
(since we know ${K_X}_{|F}$); in general
this gives only the first Chern class of the normal bundle.

\smallskip
Let us start with an easy lemma, which however gives a broader picture
of what we are actually going to prove, namely the base point freeness of
the normal
bundle.

\begin{Lemma} \cite[Proposition (3.5)]{AW3}
Let $\f: X\ra Z$ be a Fano-Mori or crepant contraction of a smooth
variety with a fiber $F=\f^{-1}(z)$. Assume that $F$ is locally
complete intersection and that the blow up $\beta:\hat X\ra X$ of $X$
along $F$ has log terminal singularities. By $\hat F$ we denote the
exceptional divisor of the blow-up.
Then the following properties are equivalent:
\par
\item{(a)} the bundle $N^*_{F/X}$ is generated by global sections on $F$,
\item{(b)} the invertible sheaf
$\O_{\hat X}(-\hat F)$ is generated by global sections
at any point of $\hat F$.
\item{(c)} $\f^{-1}m_z\cdot\O_X=\I_F$ or, equivalently,
the scheme-theoretic fiber structure of $F$ is reduced, i.e.~$\tilde F=F$.
\item{(d)} there exists a Fano-Mori contraction
$\hat\f : \hat X \ra \hat Z=Proj_Z(\bigoplus_k m_z^k)$
onto a blow-up of $Z$ at the maximal ideal of $z$,
and $\f^*(\O_{\hat Z}(1)) =
\O_{\hat X}(1)$.
\label{normal}
\end{Lemma}

\medskip
\subsubsection{The normal bundle of a 1-dimensional fiber.}

The case in which $F$ is a fiber of dimension 1 was mainly studied,
after S. Mori, by T. Ando, \cite{An}; we will report
some of his results and we will introduce an
alternative proof, as done in \cite{AW3}.

\smallskip
Let $C$ be an irreducible component of $F$.
As we saw in Lemma \ref{onefib},
$C$ is a rational curve and it can be
either a line or a conic with the respect to
$-K_X$. In the latter $\f$ is of fiber type.
\par
Let $\I$ be the ideal
of $C\subset X$ (with the reduced structure) and consider the exact sequence
$$0 \ra \I/\I^2 \ra \O_X /\I^2
  \ra \O_X /\I \ra 0.$$
In the long cohomology sequence associated
the map of global sections $H^0(\O_X /\I^2)\ra H^0({\O_X /\I})$
is surjective. Moreover,
by \ref{van-fib1}, we have the vanishing $H^1({\O_X /\I^2})=0$.
Therefore $H^1({\I/\I^2}) = 0$ which gives a bound
on $N^*_{C/X}=\I/\I^2$. Namely if $N_{C/X} = \oplus \O(a_i)$ then
$a_i < 2$.  On the other hand, by adjunction,
$det(N_{C/X})=  \Sigma a_i =\O(-2-K_X\cdot C)$ and thus the list of possible values
of $(a_1,...,a_{n-1})$ is finite.
\par
If $\f$ is a good birational
contraction then we have even a better bound because,
similarly as above and using \ref{van-fib1}, we actually
get $H^1(N^*_{C/X}\otimes\O(K_X\cdot C))=0.$
Therefore, since $K_X\cdot C =1$,
there is only one possibility, namely $N_{C/X}=\O(-1) \oplus\O^{(n-2)} $.
\par
If $\f$ is of fibre type then the estimate coming
from this technique is not sufficient and one has to use other
arguments. More precisely, one has to deal with a scheme associated
to a double structure on $C$ --- see \cite{An}.

It is also convenient to
use arguments coming from the deformation theory.
Namely, the possibilities which can occur from the above vanishing,
if $n=3$, are the following:
$$\begin{array}{llll}
\O\oplus\O, & \O\oplus\O(-1), & \O(1)\oplus\O(-2), &
\O(1)\oplus\O(-1).\\
\end{array}$$
We will show that the last possibility does not occur using an argument
related to the deformation technique. It can
be used to deal with 2 dimensional fibers too.

\begin{Lemma}
The normal bundle $N_{C/X}$ cannot be  $\O(1)\oplus\O(-1)$.
\label{onenorm}
\end{Lemma}

\begin{proof} Assume the contrary and let $\psi : \hat X \ra X$ be the
blow-up of
$X$ along $C$.
Let $E:=\Proj(\O(1)\oplus\O(-1))$ be the exceptional divisor.
Let $C_0$ be the curve contained in $E$ which is the section of the ruled
surface
$E \ra C$ corresponding to the surjective map
$(\O(1)\oplus\O(-1)) \ra\O(-1)$. We have immediately
that
$E\cdot C_0 = 1$ and that $\psi_{C_0}$ is a $1-1$ map from
$C_0$ to $C$. Therefore $K_{\hat X}\cdot C_0 = K_X\cdot C +E\cdot C_0 = -1$. In
particular this
implies that $C_0$ moves at least in a $1$-dimensional family on $\hat X$
(see \ref{hilb}). Since it
does not move on $E$ it means that it goes out of $E$. Since $C_0$ is
contracted by $\f \circ \psi$ it implies that $E\cdot C_0 = 0$,
but this is a contradiction since $E\cdot C_0 = 1$.
\end{proof}

\subsubsection{The normal bundle of a two-dimensional fiber.}

The case in which $F$ is a fiber of dimension 2 was studied,
by M. Andreatta and J.A. Wi\'sniewski, \cite{AW3}.
We will report here the main results contained in section 5.7 of \cite{AW3},
referring to it for proofs and more details.

\smallskip
Let us present in general the point of view of \cite{AW3}.
To understand higher dimensional fibers of Fano-Mori contractions we will
slice them down. Thus we will need some kind of "ascending property".

Suppose that $\f: X\ra Z $ is a Fano-Mori contraction of a smooth variety,
${\mathcal L}$ is an ample line bundle on $X$ such that 
$-K_X-{\mathcal L}$ is $\f$-(nef\& big);
for instance if $\f$ is birational one can take ${\mathcal L} = L := -K_X$.
Let $F=\f^{-1}(z)$ be a (geometric) fiber of $\f$. Suppose that
$F$ is locally complete intersection. Let $X'\in |{\mathcal L}|$ be a
normal divisor which does not contain any component of $F$.
Then the restriction of $\f$ to $X'$, call it $\f'$,
is a contraction, either Fano-Mori or crepant (see \ref{horizontal}).
The intersection $F'=X'\cap F$ is then a fiber of $\f'$.
The regular sequence of local generators $(g_1,\dots,g_r)$
of the ideal of the fiber $F$ in $X$ descends
to a regular sequence in the local ring of $X'$ which defines
a  subscheme $F\cdot X'$ supported on $F'=F\cap X'$,
call it $\bar F'$. Let us note that
if the divisor $X'$ has multiplicity 1 along each of the components of
$F$ then, since a locally complete intersection has no
embedded components, we get $\bar F'=F'$.

\begin{Lemma}  The scheme $\bar F'$
is locally complete intersection in $X'$ and
$$N^*_{\bar F'/X'}\otimes_{\O_{\bar F'}}\O_{F'}\iso (N^*_{F/X})_{|F'}.$$
If moreover $X'$ is smooth, ${\mathcal L}$ is spanned
and $dimF'=1$, then
$$H^1(F',(N^*_{F/X})_{|F'})=0$$.
\label{slicing}
\end{Lemma}

\begin{proof} The first part of the lemma follows from the preceding discussion
so it is enough to prove the vanishing.
Let $\J$ be the ideal of $\bar F'$ in $X'$.
 From \ref{van-fib1} we know that $H^1(\bar F',\O_{X'}/\J^2)=0$ and since
we have an exact sequence
$$0\raa \J/\J^2=N^*_{\bar F'/X'}\raa
\O_{X'}/J^2\raa \O_{X'}/\J=\O_{\bar F'}\raa 0$$
then we will be done if we show
$H^0(\bar F',\O_{\bar F'})=\C$.
Since $H^1(\bar F',\O_{\bar F'})=0$ then this is equivalent to
$\chi(\O_{\bar F'})=1$.
The equality $H^0(\bar F',\O_{\bar F'})=\C$ is clear
if $\bar F'$ is reduced. But since ${\mathcal L}$ is spanned and $F$ is locally
complete
intersection
then there exists a flat deformation of $\bar F'$ to another
intersection $F\cdot X''$ which is reduced. This is what we need,
because flat deformation
preserves Euler characteristic.
\end{proof}

\medskip
Now let us consider the following ascending property.
Take a point $x\in F'$. Suppose that the ideal of $F'$,
or equivalently $N^*_{F'/X'}$, is generated by global functions from
$X'$. That is, there exist global functions $g'_1, \dots g'_r\in
\Gamma(X',\O_{X'})$
which define $F'$ at $x$. Then, since $H^1(X,-{\mathcal L})=0$ these 
functions extend
to $g_1, \dots g_r\in \Gamma(X,\O_{X})$ which define $F$. Thus passing from
the ideal $\I$ to its quotient $\I/\I^2$ we get the first part of

\begin{Lemma}
If $N^*_{F'/X'}$ is spanned at a point $x\in F'$
by global functions from $\Gamma(X',\O_{X'})$
then $N^*_{F/X}$ is spanned at $x$ by functions from
$\Gamma(X,\O_{X})$.
If $N^*_{F'/X'}$ is spanned by global functions from $\Gamma(X',\O_{X'})$
everywhere on $F'$ then $N^*_{F/X}$ is nef.
\label{ascending}
\end{Lemma}

\begin{proof} We are only to proof the second claim of the Lemma. Since
$F'\subset F$ is an ample section then the set where $N^*_{F/X}$ is not
generate by global sections is finite in $F$. Therefore the restriction
$(N^*_{F/X})_{|C}$ is spanned generically for any curve $C\subset F$ and
consequently it is nef.
\end{proof}

\smallskip
If the fiber is of dimension 2 then we have a better extension property.

\begin{Lemma} Let $\f:X\ra Z$ be a Fano-Mori
birational contraction of a smooth variety
with a 2-dimensional
fiber $F$ which is a locally complete intersection.
Let $L=-K_X$; it is a $\f$-ample line bundle which can be
assumed $\f$-very ample (see \ref{aw} and \cite[Proposition (1.3.4)]{AW3}).
Then the following conditions are equivalent:
\item{(a)}$N^*_{F/X}$ is generated by global sections at any point of $F$
\item{(b)} for a generic (smooth) divisor $X'\in |L|$ the bundle
$N^*_{F'/X'}$ is generated by global sections at a generic point of
any component of $F'$.
\end{Lemma}

\begin{proof} The implication (a)$\Rightarrow$(b) is clear. To prove the
converse
we assume the contrary. Let $S$ denote the set of points on $F$ where
$N^*_{F/X}$ is not spanned. Because of the extension property \ref{ascending}
and the fact that for a one dimensional fiber of a F-M contraction
the spannedness of the normal bundle
is equivalent to spannedness at a generic point (see \cite[Corollary
5.6.2]{AW3},),
the set does not contain $F'$ and thus it is finite. Now we choose another
smooth section
$X'_1\in |L|$ which meets $F$ along a (reduced)
curve $F'_1$ containing a point of $S$.
(We can do it because $L$ is $\f$-very ample.)
The bundle $N^*_{F'_1/X'_1}$
is generated on a generic point of $F'_1$ so it is generated everywhere but
this, because
of the extension property, implies that $N^*_{F/X}$ is generated at some
point of $S$,
a contradiction.
\end{proof}

\begin{Lemma}
Let $\f:X\ra Z$ be a Fano-Mori birational contraction of a smooth
4-fold with a 2-dimensional fiber $F=\f^{-1}(z)$. As usually $L=-K_X$
is a $\f$-ample line bundle which may be assumed to be $\f$-very
ample. Suppose moreover that either $F$ is irreducible or $L^2\cdot F\leq 2$
(which is the case when $F$ is an isolated 2 dimensional fiber).
Then the scheme fiber structure $\tilde F$ is reduced
unless one of the following occurs:
\item{(a)} the fiber $F$ is irreducible and the restriction of
$N_F$ to any smooth curve $C\in |L_{|F}|$ is isomorphic to
$\O(-3)\oplus\O(1)$,
\item{(b)} $F=\Proj^2\cup\Proj^2$ and the restriction of
$N_F$ to any line
in one of the components is isomorphic to
$\O(-2)\oplus\O(1)$.
\end{Lemma}

\begin{proof}  Let us consider a curve $C\in |L_{|F}|$.
Since $L$ is $\f$-very ample we can take a smooth $X'\in |L|$ such that
$\f'=\f_{|X}$ is a crepant contraction and $C=F\cap X'$
(see \ref{horizontal}). Then, considering the
embeddings $C=F\cap X'\subset F\subset X$ and $C=F\cap X'\subset X'\subset
X$, we get
$$N_{C/X}=N_{C/X'}\oplus L_C=(N_{F/X})_{|C}\oplus L_C$$
and therefore $N_{C/X'}=(N_{F/X})_{|C}$.

Now the normal bundles $N_{C/X'}$ of the crepant contraction $\f'$
can be easily described, in a way similar to the previous section
(see (\cite{AW3}, 5.6.1)).
In particular it follows that if neither (a) nor (b)
occurs then the fiber structure of the contraction $\f'$ is reduced.
Thus, using the previous lemma and
the equivalence in \ref{normal}, we conclude that $\tilde F=F$.
\end{proof}

\medskip
Now, one has to discuss the possible exceptions described in the above lemma.
This is done extensively in \cite{AW3} and the following was proved:

\begin{Theorem} (\cite{AW3}, theorems 5.7.5 and 5.7.6)
Let $\f:X\ra Z$ be a Fano-Mori birational contraction of a smooth 4-fold
with an
isolated 2-dimensional fiber $F=\f^{-1}(z)$.
Then the fiber structure $\tilde F$ coincides with the geometric structure $F$
and the conormal bundle $N^*_{F/X}$ is spanned by global sections.

Moreover if $F=\Proj^2$ then $N^*_{F/X}$ is either
$\O(1)\oplus\O(1)$ or $T(-1)\oplus\O(1)/\O$, or
$\O^{\oplus 4}/\O(-1)^{\oplus2}$.
If $F$ is a quadric (possibly singular or even reducible) then $N^*_{F/X}$ is
the spinor
bundle ${\mathcal S}(1)$.
\label{twonormal}
\end{Theorem}

\medskip
In some respects the above results about the fiber structure of a
2-dimensional fiber are nicer than one may expect. Namely, there is
no multiple fiber structure, the conormal bundle is nef and the normal
bundle of the
geometric isolated fiber has no section. Thus the situation is better than
for 1-dimensional isolated fibers in dimension 2 and 3: the
fundamental cycle of a Du Val $ADE$ surface singularity is
non-reduced and in dimension 3 one may contract an isolated $\Proj^1$
with the normal bundle $\O(1)\oplus\O(-3)$.  On the other hand, using a
double covering construction (see \cite[Examples(3.5)]{AW2})
in dimension 5 one may contract a
quadric fibration over a smooth 3-dimensional base with an isolated
fiber equal to $\Proj^2$, scheme theoretically the fiber is a double
$\Proj^2$. Using the sequence of normal bundles and the deformation
of lines argument, one may verify that in this case $N_F\iso \O(1)\oplus\E^*$
where $\E$ is a rank 2 spanned vector bundle with $c_1=2$, $c_2=4$,
so that $dimH^1(\E^*)=-\chi(\E^*)=3$.

\smallskip
Let us also note that for divisorial fibers we have the following:

If $F=\bigcup F_i$ is a divisorial fiber of a surjective map
$X\ra Y$, where $X$ is smooth and $dimY\geq 2$ then for some $k>0$
the line bundle $\O_{F_i}(-kF_i)$ has non-trivial section and thus
no multiple of $\O_{F_i}(F_i)$ has a section.
In particular, if $rank(Pic(X/Y))=1$ then $\O_F(-F)$ is ample.

\smallskip
One can then try to conjecture that if $F$ is an isolated fiber
of a (Fano-Mori) contraction
which is locally complete intersection and with ``small'' codimension
then $H^0(F,N_F)=0$.

\bigskip
The above result on contractions of 4-folds can be generalised for
the adjunction mappings of an $n$-fold. Namely, suppose that
$\f: X\ra Y$ is a Fano-Mori contraction of a smooth $n$-fold $X$ supported
by a divisor $K_X+(n-3)H$, where $H$ is a $\f$-ample divisor on $X$.
Since we are interested in the local description of $\f$ around a non-trivial
fiber $F=\f^{-1}(z)$, we may assume that the variety $Z$ is affine.

\begin{Corollary}
Let us assume that $\f$ is birational and that $F$
is an isolated fiber of dimension $n-2$. If $n\geq 5$ then
the contraction $\f$ is small and
$F$ is an isolated non-trivial fiber. More precisely
$F\iso\Proj^{n-2}$ and $N_{F/X}\iso\O(-1)\oplus\O(-1)$, and there exists
a flip of $\f$ (see [A-B-W]).
\end{Corollary}

\bigskip
Note that the preceding arguments, which led to the classification
of the birational 4-dimensional case, depend on the isomorphism
$\f'_*\O_{X'}\iso\O_Z\iso\f_*\O_X$. This fails to be true if
$\f$ is of fiber type. Namely, let $\f:X\ra Z$ be a conic fibration,
i.e.~a Fano-Mori contraction such that $dimZ=dimX-1$. As usually,
we will assume that $F$ is an isolated
  2 dimensional fiber of $\f$ and $L=-K_X$ is
$\f$-spanned. Then the restriction of $\f$ to
a general section $X'\in |L|$ is generically $2:1$
covering of $Z$.
A different argument is developed for this case in \cite{AW2}, where the
following theorem was proved.

\begin{Theorem} [\cite{AW3}, proposition 5.9.5 and theorem 5.9.6]\

Let $\f:X\ra Z$ be a conic fibration of a smooth 4-fold.
Suppose that $F$ is an irreducible isolated
2-dimensional fiber of $\f$ which is either a projective
plane or a quadric. Then the conormal bundle $N^*_{F/X}$
is nef.

Moreover if $F\iso \Proj^2$ then
$N^*_{F/X}\iso \O^3/\O(-2)$ or
$T\Proj^2(-1)$. If $F$ is an irreducible quadric
then $N^*_{F/X}$ is the pullback of $T\Proj^2(-1)$ via some double covering
of $\Proj^2$. In both cases the scheme fiber structure $\tilde F$
is reduced and $Z$ is smooth at $z$.
\label{twonormalfiber}
\end{Theorem}

\section[Normal bundle and formal completion]{When the normal of a 
fiber determines
locally the contraction?}

We illustrate some results which show how the normal
bundle can give all the informations we want on the contraction.
That is when the second order approximation actually determines completely
the formal neighborhood. This is
of course not always the case.

In some situations the knowledge of the normal bundle $N_{F/X}$ allows
to determine the singularity of $Z$ at $z = \f(F)$.
Let us recall that for a local ring $\O_{Z,z}$ with the maximal ideal
$m_z$ one defines the graded $\C$-algebra $gr(\O_{Z,z}):=\bigoplus_k
m_z^k/m_z^{k+1}$.  The knowledge of the ring $gr(\O_{Z,z})$ allows
sometimes to describe the completion ring $\hat\O_{Z,z}$.  Also, we
will say that a spanned vector bundle $\E$ on a projective variety
$Y$ is p.n.-spanned (p.n. stands for projectively normal) if for any
$k> 0$ the natural morphism $S^kH^0(Y,\E)\ra H^0(Y,S^k\E)$ is
surjective. As we noted while discussing the contraction to the
vertex, projective normality allows us to compare gradings of rings
``upstairs'' and ``downstairs''.

\begin{Proposition}[\cite{Mo}]
Let $\f:X\ra Z$ be a contraction as above.  Suppose moreover that
$N^*_{F/X}$ is p.n.-spanned. Then $\f_*(\I^k_F)=m^k_z$,\ \
$\f^{-1}(m^k_z)\cdot\O_X=\I^k_F$ and there is a natural isomorphism
of graded $\C$-algebras:
$$gr(\O_{Z,z})\iso\bigoplus_kH^0(F,S^k(N^*_{F/X})).$$
\label{graded}
\end{Proposition}

We omit the proof of the above result referring to Mori, \cite{Mo},
p.164, who proved it in case when $F$ is a divisor, the
generalisation is straightforward.

The next is a version of a
theorem of Mori \cite{Mo}, 3.33, which is a generalisation of
a Grauert-Hironaka-Rossi result:

\begin{Proposition}
Suppose that $F$ is a smooth fiber of a Fano-Mori or crepant
contraction $\f:X\ra Z$ and assume that its conormal bundle
$N^*=N^*_{F/X}$ is nef. If $H^1(F,T_F\otimes S^i(N^*))=H^1(F,N\otimes
S^{i}(N^*))=0$ for $i\geq 1$ then the formal neighborhood
of $F$ in $X$ is determined uniquely and it the same
as the formal neighborhood of the zero section in
the total space of the bundle $N$.
\label{formaln}
\end{Proposition}

Also the following assertion is a straightforward generalisation of the
celebrated Castelnuovo contraction criterion for surfaces; its proof
is similar to the one of \cite{Ha}, V.5.7, see also \cite{AW2}.

\begin{Proposition} [Castelnuovo criterion]
Let $\f :X \ra Z$ be a projective morphism from a smooth variety $X$
onto a normal variety $Z$ with connected fibers (a contraction). Suppose
that $z \in Z$ is a point of $Z$ and $F = \f^{-1} (z)$ is the
geometric fiber over $z$ which is locally complete intersection in
$X$ with the conormal bundle $N^*_{F/X}$.  Assume that for any
positive integer $k$ we have $H^1( F, S^k(N^*_{F/X})) = 0$ (note that
this assumption is fulfilled if $\f$ is Fano-Mori contraction,
$N^*_{F/X}$ is nef
and the blow-up of $X$ at $F$ has log terminal singularities). If for
any $k\geq 1$ it is
$S^k H^0(F,N^*_{F/X})\iso H^0(F ,S^k(N^*_{F/X}))$
then $z$ is a smooth point of $Z$ and $dimZ=dimH^0(F,N^*_{F/X})$.
\label{castelnuovo}
\end{Proposition}

\section[Remarks on smooth contractions]{Concluding remarks on the
classification of Fano-Mori contractions on a smooth $n$-fold with $n\leq 4$}

The proofs of the theorems announced in the first two sections of this part
can be given applying the numerous results we have given up to now
(very often not in a unique way).
This may be not so trivial, so in this section we
will give some possible schemes of proof.

\smallskip
A good starting point is to use the Ionescu-Wi\'sniewski inequality in
\ref{IoWi}, this gives the possibilities for the dimension
of the fibers and of the exceptional locus. In particular
it says that there are no small contractions (i.e. contractions
whose exceptional locus has codimension $\geq 2$)
on a smooth threefolds
and also it proves the part 0 of the Theorem \ref{4-folds}.

\medskip \noindent
{\bf  Description of the F-M contractions around a fiber $F$ of dimension $1$
(general case)}.

This is given in \ref{Ando}; this covers almost all Theorem
\ref{2-folds} and part of Theorems \ref{3-folds} and \ref{4-folds}.
The proof of \ref{Ando} follows from \ref{onefib} (which describes
the possibilities for $F$),  \ref{onenorm} and the discussion before it
(which describes the possibilities for the normal bundles; in \ref{onenorm}
only the case $n=3$ is discussed in details)
and \ref{castelnuovo}.

Note that in the case of surface if $\f = cont_R$ is a conic bundle
then actually $\f$ gives the structure of a minimal ruled surface.
In order to prove this we have to show that there are no
reducible or non reduced fiber of $\f$
In fact if, by contradiction,  $F$ is such a fiber then
$F = \sum a_i C_i = [C]$ with $[C]  \in R$.
But since $R$ is extremal this implies that $C_i \in R$ for every $i$.
Thus $C_i^2 = 0$, since a general fiber of $\f$ is a smooth irreducible and
reduced
curve in the ray, and ${C_i}\ ^. K_X < 0$. By the adjunction formula this
implies that
$C_i \iso \Proj ^1$ and ${C_i}\ ^. K_X = -2$. Thus
$$-2 = (C\ ^. K_X) = \sum a_i({C_i}\ ^. K_X) = -2 \sum a_i,$$
which gives a contradiction.
Furthermore using Tsen's Theorem, one can prove that $X$ 
 is the projectivization of a rank two vector bundle on $\Proj ^1$, see \cite[C.4.2]{Reid}.

\medskip \noindent
{\bf  Description of the F-M contractions around a fiber $F$ of dimension $2$
(3-folds and 4-folds)}.

In the case of surface we have that the contraction
$\f$ contracts $X$ to a point; that is that $X$ has
$Pic = \Z$ and $-K_X $ is a ample, i.e. $X$ is a Fano surface of
Picard number $1$. Then one can prove
that $X = \Proj ^2$; see for instance \cite[ p. 21]{CKM}.

In the threefolds case we have that either
$\f$ is a contraction of fiber type contracting $X$ to a curve, with all
fibers of dimension two,
or $\f$ is a birational divisorial contraction which contracts a unique
prime divisor equal to $F$
to a point: in fact, by Exercise \ref{divisorial}, $\f$ cannot be a 
contraction to a
surface with
some isolated two dimensional jumping fibers and if $\f$ is birational then
the
exceptional divisor is prime.

In the first case, by adjunction, the general fiber is a smooth del Pezzo
surface.

In the birational case we can apply \ref{twofiber} which gives all the
possibilities
(namely $3$) for the exceptional divisor $F$; by adjunction we easily
compute the normal
bundle and we prove the uniqueness of the analytic neighborhood of the
contraction around $F$
by using \ref{castelnuovo}, \ref{graded} and \ref{formaln}.
Actually for the uniqueness in the case of a fiber isomorphic to the singular
quadric
an extra argument is needed (see for this \cite[p. 165]{Mo}).
Note also that all these cases exist and can be constructed via the basic
example \ref{examplev}
except the case with a fiber isomorphic to the singular quadric;
in this case if we work as in \ref{examplev} we construct a singular
threefold $X$. Moreover the case of the smooth quadric can be constructed
via \ref{examplev} but not as an elementary contraction, i.e. as a contraction
of a single ray.
Two good examples were given in \cite[3.44.2 and 3.44.3]{Mo}.

\smallskip
We consider then Fano Mori elementary contraction from a fourfold,
in a neighborhood of a two dimensional fiber $F$.
We start with the birational case:
the fiber can be an isolated two dimensional fiber or can stay in a
one dimensional family of two dimensional fibers
(the other possibilities are ruled out with Exercise \ref{divisorial}).

The first case is described in Part 4 of the Theorem \ref{4-folds};
to prove it we can first apply the Theorem
\ref{twofiber}, which gives 4 possibilities for the fiber $F$, namely $F$
can be
$\Proj^2$ or a reduced quadric (smooth, singular or reducible).

Then we apply the Theorem \ref{twonormal} which describes the possible
normal bundles.
If $F$ is smooth then we can prove the uniqueness of a formal neighborhood
of the fiber using \ref{formaln} and construct an example using \ref{examplev}.

If $F$ is a singular quadric or a reducible one (union of two $\Proj ^2$
meeting along a line)
then the situation is more complicated. In \cite{AW4} one can find
good examples for these situations; moreover in the case $F$ is reducible
there are at least two possible analytically non equivalent formal 
neighborhood
of $F$. That is, as one expects, the fiber and its normal does not always
determine the analytic neighborhood. The (open) question is whether in this
case these two
neighborhoods are the only possible ones (up to analytic equivalence).

Part 3 of the Theorem \ref{4-folds} describes the case in which $F$
stays in a one dimensional family of two dimensional fibers.
The proof of it is essentially different from the previous one
and it was given in \cite{AW4}.
Using a vertical slicing, see the proposition \ref{vertical}, one can prove
that the general two dimensional
fiber of this family is either $\Proj^2$ or a irreducible and reduced quadric.
In fact, in the notation of \ref{vertical} (vertical slicing),
$f_{X_h} :X_h \ra f(X_h)$ is a
Fano Mori contraction
from a smooth threefold which contracts a general fiber to a point; thus we
can apply the
result on threefolds , i.e. \ref{3-folds}
(this is actually a bit quick: in fact
$f_{X_h}$ can be non elementary, i.e. the contraction of a face, not
of a ray. In \cite{AW4} this is in fact ruled out).
It is easy to prove that if the general fiber is $\Proj^2$ then the same is
for the special fiber.
If the general fiber is a quadric then the special one is also a quadric,
but very likely
more singular. It turns out that there are no non reduced quadric, that is
double $\Proj ^2$,
as special
fiber; the other possibilities all occur. We refer to \cite{AW3} for
further details and examples.

\smallskip
We finally pass to the case of fiber type Fano-Mori contractions with two
dimensional fibers.
If the two dimensional fiber is not isolated then the contraction is to a
surface,
by Exercise \ref{divisorial} it is equidimensional, i.e. all fibers 
are of dimension
two,
and its general fiber is a del Pezzo surface (by adjunction formula).

If the two dimensional fiber $F$
is isolated then it is one of those described in the Theorem \ref{twonormal};
not all the given possibilities are locally complete intersections.
If $F$ is a smooth $\Proj^2$ or a smooth quadric then the normal bundle is
computed
in the Theorem \ref{twonormalfiber}.
In these two cases example can be constructed using \ref{examplev};
may be one can also prove the uniqueness of a formal neighborhood
of the fiber using \ref{formaln} (this can be a hard computation!).
Some examples were constructed for other possible fibers
but for some of them we cannot even construct an example (see
the remark after the Theorem \ref{twofiber}).

\medskip \noindent
{\bf  Description of the F-M contractions around a fiber $F$ of dimension $3$
(4-folds)}.

If the contraction is birational by Exercise \ref{divisorial} $F$ is the unique
prime divisor
equal to the exceptional locus. It is immediate to see,
using adjunction formula, that $F$ is a del Pezzo threefolds. The
problem here is to prove that $F$ is normal and which normal non smooth
del Pezzo threefolds can actually occur (not all of them!)
(see Part 2 of the Theorem \ref{4-folds} and the following remark).

If the contraction is of fiber type then, again by \ref{divisorial},
all fibers are three dimensional and the generic one is a Mukai manifold.

\section[Fano Manifolds]{Classification of Fano manifolds of high index}

This section is devoted to the study of F-M contractions of a smooth
manifolds with target a point, i.e. to Fano manifolds.

We already noticed that Fano varieties
with high index, with respect to the dimension,
are easier to be understood.
Namely for $i(X)\geq \dim X$
the informations given by the Hilbert polynomial are already sufficient to
give a complete description of all possible cases, see Exercise
\ref{highindex}.

It is not surprising that for
lower indexes
the world is wilder. The best known way to go further is, again, the following
adjunction procedure.

Let $X$ be a Fano manifold of index $i(X)=r$ and fundamental divisor $L$
(see the Definitions \ref{fanodef} and also \ref{tipo}).
Assume that
$|L|$ is not empty. Let $H\in|L|$ a generic member. Then by adjunction formula
$$-K_H=-(K+L)_{|H}\sim (r-1)L_{|H}.$$
In other words whenever $r>1$ the section $H$ is a Fano variety of the same
dual index.
So that, if one is able to control the  singularities of $H$, then it is
possible to
study $X$ through $H$.
More generally we can pose the following.
\begin{Definition} Let $f:X\ra S$ be a local contraction of type
$(d,\gamma,\Phi)$, supported by
$K_X+rL$. Then we will say that $f$ has {\sf good divisors} if,
after maybe shrinking S,
the generic
element $H\in|L|$
has at worst the same singularities as $X$ and
$f_{|H}:H\ra S_H$ is of type $(*,*,\Phi)$.
\end{Definition}

Assume that the good divisor problem has a positive answer for a fixed
index. Then
the above observation allows to classify all Fano manifolds of fixed dual
index in an
inductive way, starting from the lower dimensional ones.

This is what Fujita did, see \cite{Fu}, for del Pezzo manifold, i.e.
Fano manifold with $i(X)=\dim X-1$.
His results can be summarised in the following way.
\begin{Theorem}[\cite{Fu}] Let $X$ be a del Pezzo manifold of dimension $n$.
If $n=2$ then $X$ is  either $\Q^2$ or $\Proj^2$ blown up in $r\leq 8$
general points.
Assume $n\geq 3$. Let $d=H^n$ the degree of $X$. Then we have the following
cases:
\begin{itemize}
\item[$d=1$] $X_6\subset\Proj(1^{n-1},2,3)$, i.e. a hypersurface of degree 6
in the weighted projective space with weights $(1,\ldots,1,2,3)$;
\item[$d=2$] $X_4\subset\Proj(1^{n},2)$;
\item[$d=3$] $X_3\subset\Proj^{n+1}$;
\item[$d=4$] $X_{2,2}\subset\Proj^{n+2}$, i.e. a complete intersection of
two Quadrics
\item[$d=5$] a linear section of $\G(1,4)\subset\Proj^9$;
\item[$d=6$] $X$ is either $\Proj^1\times\Proj^1\times\Proj^1$ or
$\Proj^2\times\Proj^2$ or
$\Proj_{\Proj^2}(T\Proj^2)$;
\item[$d=7$] the blow up of $\Proj^3$ in one point
\end{itemize}
\end{Theorem}

\medskip
A tutorial case is the surface one, which is given through the 
following exercise
entirely based on Minimal Model techniques.

\begin{exercise} Let $S$ be a del Pezzo surface of degree $d$;
then the Kleiman-Mori cone is spanned by extremal rays, see Theorem 
\ref{cone}.
We also have a precise description of the contraction associated to each
extremal ray, see \ref{MMPsurf}.
\begin{itemize}
\item Assume that $S$ has an extremal ray $R$ whose associated
contraction is birational; that is $\f:S\ra S'$ is the contraction
of a $(-1)$ curve C. Prove that $S'$ is a del Pezzo surface of degree $d+1$.
\end{itemize}

Since going from $S$ to $S'$ we decrease the Picard number of one,
after finitely many contractions we
have a del Pezzo surface
$S_k$ with only fiber type extremal
rays.

\begin{itemize}
\item Show that $S_k$ is either $\Q^2$ or $\Proj^2$
(use \ref{MMPsurf}).
\end{itemize}

We have thus established that any del Pezzo surface $S$ is either $\Q^2$ or the
blow up of $\Proj^2$
in a finite number of points (a blow up of $\Q^2$ is in fact a
blow up of $\Proj^2$).

\begin{itemize}
\item Show that you cannot blow up more than 8 points with the following
restrictions: no 3 are on a line and no 6 on a conic.

\hint evaluate the self intersection of $K_S$
and the intersection of
$K_S$ with the strict transform of line or conic.
\end{itemize}

We have now to conclude.
\begin{itemize}
\item A surface obtained by blowing up $r$ points of $\Proj^2$, with the above
restrictions,
is a del Pezzo surface.

\hint Study either the combinatorics of the cone of effective curves,
or the linear systems of cubics with imposed conditions.
\end{itemize}

One can improve the knowledge of these surfaces observing that whenever
$d\geq 3$ then
$|-K_S|$ is very ample and embeds $S\subset\Proj^d$. For for $d=2$ 
the complete linear system
$|-K_S|$ is spanned and gives
a double cover of $\Proj^2$ ramified along a quartic. While for $d=1$ 
the system
$|-2K_S|$ is spanned and gives
a double cover of a singular quadric ramified along a sextic and the vertex.
\end{exercise}

\medskip
In higher dimension the idea of the proof is the following.
First show that the good divisor problem has a
positive answer, i.e. prove the following exercise.

\begin{exercise} Let $X$ be a del Pezzo manifold. Prove that $X$ has good
divisors.
\label{delpezzoex}
\end{exercise}

hint: Read first the proof of the following Theorem \ref{good_mukai}.

\smallskip
Then all the information on del Pezzo
surface can
be extended.

The first two cases can be obtained with the machinery of
graded rings,
\cite[\S 3]{MoCI}.
We would like to stress here the following property.
If a variety $X$, of dimension at least 4, contains an hyperplane
section which is a weighted complete intersection then the variety 
itself is a weighted complete
intersection, \cite[Corollary 3.8]{MoCI}.

If degree $d\geq 3$ then $|-K_X|$
is very ample. So that $d=3,4$ are immediate while for $d\geq 5$ the study
is more subtle and we leave
it for the interested reader, see \cite{Fu}.

\bigskip
The next case is the one of Mukai manifold, i.e. manifold with $i(X)=dim X-2$.
These varieties are named after S. Mukai who first
announced their classification, \cite{Mu}, assuming the existence of good
divisors.

This assumption is proved in \cite{Me}, where the Base Point 
Free technique is applied to answer the good divisor problem for Mukai varieties.

The idea is simple. Let $X$ be a Fano manifold and $|L|$ the 
fundamental divisor of $X$.
Let $D\equiv \delta L$ a $\Q$-divisor with $\delta<1$.  By BPF
technique  there is a section
of $|L|$ non vanishing identically on
$LLC(X,D)$. If we sum up with  Bertini Theorem we immediately get that the
generic section
of $|L|$ cannot have worse than LC singularities. We actually prove the
following.

\begin{Theorem}[\cite{Me}]Let $X$ be a
Mukai
variety with at worst log terminal singularities.
Then $X$ has good divisors except in the following cases:
\begin{itemize}
\item[i)]  $X$ is a singular terminal Gorenstein 3-fold which is
  a ``special'' (see \cite{Me}) complete intersection of a quadric
and a sextic in $\Proj(1,1,1,1,2,3)$
\item[ii)] let $Y\subset \Proj(1,1,1,1,1,2)$ be a ``special'' complete
intersection of a quadric cone and a quartic; let $\sigma$ be the 
involution on $\Proj(1,1,1,1,1,2)$
given by $(x_0:x_1:x_2:x_3:x_4:x_5)\mapsto(x_0:x_1:x_2:-x_3:-x_4:-x_5)$
and let $\pi$ be the map to the quotient space.
  Then $X=\pi(Y)$ is a terminal not Gorenstein 3-fold.
\end{itemize}
In both exceptional cases the generic element of the
fundamental divisor has canonical singularities and $Bsl|-K_X|$ is a 
singular point.
It has to be stressed that the generic 3-fold in i) and ii) has good 
divisors, but there are
``special'' complete intersections whose quotient has a singular 
point in the base locus
of the anticanonical class, see \cite[Examples 2.7, 2.8]{Me} for details.
\label{good_mukai}
\end{Theorem}
\begin{proof} We prove the theorem in four steps, from the
more singular to the smooth case. As usual we argue by contradiction.

\begin{claim} Assume that $X$ has log terminal singularities. Then 
$X$ has good divisors.
\label{cl:LT}
\end{claim}

\begin{proof}[Proof of the Claim] Proposition \ref{al} ensures that 
$dim |L|\geq 1$.
Let $H\in|L|$ a generic section and assume that $H$ has worse
singularities than LT. Let
$\gamma=lct(X,H)$ and $Z\in CLC(X,\gamma H)$ a minimal center.
\begin{exercise}$\gamma\leq 1$ and $cod Z\geq 2$.

\hint if $cod Z=1$ then $Z$ is a fixed component
of $Bsl|L|$, but $dim |L|>0$ and $H$ is ample therefore connected. In 
particular if
$H$ is reducible then there exists a codimension two in $Sing(H)$.
\end{exercise}

This is enough to derive a contradiction.
By Bertini Theorem $Z\subset Bsl|L|$, while by Lemma \ref{sectame}
there exists a section of $L$ non vanishing on $Z$.
\end{proof}

\begin{claim}$X$ has canonical singularities. Then $X$ has good divisors.
\label{cl:ca}
\end{claim}

\begin{proof}[Proof of the Claim] By Claim \ref{cl:LT} $H$ has LT 
singularities.
Let $\mu:Y\ra X$ a log resolution of $(X,H)$, with $\mu^*H=\overline{H}+
\sum r_i E_i$, where $|\overline{H}|$ is base point free,
  and $K_Y=\mu^*K_X+\sum a_i E_i$. Let us assume that $H$ has
not canonical singularities, then, maybe after reordering the indexes,
we have $a_0<r_0$. Since $H$ is generic then
$\mu(E_i)\subset Bsl|L|$, for all $i$ with $r_i>0$, see the proof of 
Claim at page \pageref{cl:clc}.
  Let $D=H+H_1$, with $H_1\in |L|$ a generic section.
First observe that $\mu$ is a log resolution of $(X,D)$ as well.
Let  $\mu^*H_1=\overline{H_1}+\sum r_i E_i$, since $H_1$ is a 
Cartier divisor then the
$r_i$ are positive integers. In particular  $a_0+1<r_0+r_0$, hence
$(X,D)$ is not LC.
Let $\gamma=lct(X,D)<1$ and $W$ a minimal center of
$CLC(X,\gamma D)$.

\begin{exercise} Prove that $cod W\geq 3$.

\hint $H$ is LT and canonical singularities are
Gorenstein in codimension 2.
\end{exercise}

We are again in the condition to apply Lemma \ref{sectame} to derive 
a contradiction.
\end{proof}

\begin{claim} If $X$ has terminal singularities then $X$ has good 
divisors unless $X$ is as in
either i) or ii).
\label{cl:ter}
\end{claim}

\begin{proof}[Proof of the Claim] If $X$ is a terminal Mukai variety 
of dimension$\geq 4$
then by the above Claim \ref{cl:ca} $H$ has canonical singularities. 
Furthermore
we can apply Claim \ref{cl:ca} to $H_{|H}$, to deduce that even 
$H_{|H}$ has canonical singularities.
Let $f:Y\to X$ a log resolution for $(X,H)$. Assume that 
$K_Y=f^*K_X+\sum a_iE_i$ and
$f^*H=H_Y+\sum r_iE_i$. Then $K_{H_Y}=f^*K_H+\sum(a_i-r_i)E_i$ and, 
with obvious notations,
$K_{H_{|HY}}=f^*K_{H_H}+\sum(a_i-2r_i)$. We just observed that 
$a_i-2r_i\geq 0$, therefore
$a_i-r_i>0$ whenever $r_i>0$. This proves that $H$ is terminal on the 
base locus of $|L|$ and we conclude by
Bertini Theorem that $H$ is terminal.

The case of terminal 3-folds with $-K_X\equiv L$ is left; this goes a 
bit beyond the
techniques we developed and so here we only state the result
(remember that terminal surface singularities are smooth points):

\begin{Theorem}[\cite{Me}] Let $X$ be a terminal Mukai 3-fold, assume that
all the divisor in the linear system $|L|$ are singular , then $X$ is one
of the two exceptions in the Theorem \ref{good_mukai}.
They actually exist.
\label{th:ter}
\end{Theorem}
\end{proof}

\begin{claim}If $X$ is smooth then $X$ has good divisors.
\end{claim}

\begin{proof}[Proof of the Claim] Assume that the generic element in
$|L|$ is not smooth. Then a 3-fold section $T\subset X$ is one of the
two exceptions to Claim \ref{cl:ter}.
Then by the usual vanishing theorem
$$H^0(X,L)\ra H^0(T,L_{|T})\ra 0,$$
and by Theorem \ref{th:ter} $Bsl|L|=Bsl|L_{|T}|$ is just one point, say $x$.
Let $H_i\in |L|$, for $i=1,\ldots, n-1$ generic elements and
$D=H_1+\cdots+H_{(n-1)}$, then the minimal center of $CLC(X,D)$ is $x$ and
$(X,D)$ is not LC at $x$, since $2(n-1)>n$.
We, therefore,
derive a contradiction by Lemma \ref{sectame}.
\end{proof}
\end{proof}

With similar arguments one can prove a good divisor problem for other
F-M contractions, for details and related results see \cite{dpf};
for instance we have the following.

\begin{Theorem}[\cite{dpf}] Let $f:X\ra S$ a local
contraction of type $(1,1,1)$.
Assume that $X$ is smooth.
Then
  $f$ has good divisors.
\label{secsm}
\end{Theorem}

\begin{remark} Could this be the starting point of a relative
analogue of Fujita classification?
The above theorem reduces this study to that of fibrations of
surfaces.
The main problem to solve is
a base point free statement in a neighborhood of an irreducible
non reduced fiber.  With this one could provide
a structure theorem as in the absolute case, embedding this spaces
in some relative (weighted) projective space. Then
one should try to extend Andreatta-Wi\'sniewski theory one step further.
Unfortunately as far as we can say this is quite hard and requires
lot of unknown results on vector bundles on del Pezzo surfaces.
\end{remark}

Let $X$ be a smooth Fano n-fold of index $r=n-2$, i.e. a Mukai manifold.
Let $|H|$ the linear
system of fundamental divisors.
The integer
$$ g=\frac12 H^n+1,$$
is called the genus of $X$, the reason will be clear after Proposition\ref{isk}; by Riemann--Roch
$$\dim H^0(X,H)=n+g-1.$$
By Theorem \ref{good_mukai} the generic element $S\in |H|$ is smooth. As
observed above this
allows an inductive argument toward three dimensional Fano's.

\noindent
{\bf Part I: }If $rk Pic(X)=1$ we use Iskovkikh's results on Fano 
3-folds, \cite{Is},
and we have
\begin{Proposition} Let $X$ be a smooth Mukai n-fold with $rk 
Pic(X)=1$. Then $|H|$ is base
point free and one
of the following is true:
\begin{itemize}
\item[i)] $|H|$ is very ample and embeds $X$ in $\Proj^{g+n-2}$. In
particular $X\subset \Proj^{g+n-2}$
  has a smooth curve section canonically embedded.
\item[ii)] the morphism associated to sections of $|H|$ is a finite
morphism of degree 2 either onto
$\Proj^n$ (in case $g=2$) or onto $\Q^n\subset \Proj^{n+1}$ (in case $g=3$).
\end{itemize}
\label{isk}
\end{Proposition}

We are therefore restricted to study projective varieties
with smooth canonical
curve section; the actual point of the classification is to 
understand all of them.
(For a somewhat backward approach, see also \cite{CLM}.)

\begin{Theorem}[\cite{Mu}] Let $X_{2g-2}\subset \Proj^{n+g-2}$ be as in
point $i)$ of Proposition \ref{isk}.
If $g\leq 5$ then $X$ is a complete intersection.
Assume that $10\geq g\geq 6$, then we have the following picture
$$
\begin{array}{ccc}
g&n(g)&X^{n(g)}_{2g-2}\subset \Proj^{g+n(g)-2}\\
6&6&C(\G(1,4))\cap\Q\subset\Proj^{10}\\
7&10&SO(10,\C)/P\subset\Proj^{15}\\
8&8&\G(1,5)\subset\Proj^{14}\\
9&6&Sp(6,\C)/P\subset\Proj^{13}\\
10&5&G_2/P\subset\Proj^{13}
\end{array}
$$
where $C(X)$ is the cone over $X$ and $n(g)$ is the maximal dimension for a
Mukai variety of that type (as
observed before any
hyperplane section is then a Mukai variety).
If $g>10$ then $g=12$ and $n(g)=3$. Then some special $X_{22}^3$
can be seen as a smooth equivariant compactification
of $SO(3,\C)/(\mbox{\rm icosahedral group})$, \cite{MU}.
In general it is possible to give a description using net of
Quadrics.
\end{Theorem}
We remark that, as before, the case of genus$\leq 5$ is easily 
obtained, while the study of
the remaining
cases is the heart of the proof.

\medskip
\noindent
{\bf Part II:} If $rk Pic(X) >1$ then there are at least two extremal
rays on $X$.
In the following exercise
we collect all the crucial informations to complete the
biregular classification in this case, see \cite{MoMu} \cite{Mu}.

\begin{exercise}
\begin{itemize}
\item[1)] Prove that the only possible F-M birational contractions are
those of a divisor to a point or to a smooth curve.
\item[2)] Which divisors are then possible ?
\item[3)] In case of extremal rays with associated contraction of fiber type then
the base of the contraction is smooth.
\end{itemize}
hint: Use the classification provided by Theorem \ref{3-folds}
\end{exercise}
\part{Birational geometry}
\label{ch:biraz}\setcounter{section}{0}

One of the main goal of Algebraic Geometry is to achieve a
birational classification of projective algebraic varieties.
The Minimal Model Program, or Mori's Program,
is an attempt to get (part of) this classification.
In the first section of this chapter we want to introduce the
general philosophy of the MMP.

The developed techniques however allow to
treat other birational aspects as well.
In this realm we would like to focus on two
different settings:\#-Minimal Model and  Sarkisov program.

The former is a polarized minimal model program
that enables to study special uniruled 3-folds and
it will be described in section
\ref{sharp}.

The latter is used to investigate birational morphisms between
Mori spaces, see
Definition \ref{MoriSpaces}, and its application to the projective plane is
outlined in section \ref{S_2}.

\section{Minimal Model Program philosophy}

Let us present the approach of the MMP toward
the birational classification of Algebraic Varieties.

The case of smooth curves is clearly settled
by the Riemann uniformization theorem.
The case of surfaces is more complicated and
it was developed by Italian algebraic geometers at the beginning
of the twentieth century. This big achievement was used as a model
for the higher dimensional case and at the end it was included
in the more general philosophy of Minimal Model Theory as we will
see in the next.

Consider a smooth projective variety $X$. The aim of Minimal Model
Theory is to distinguish, inside the set of varieties which are
birational to $X$, a special ``minimal''
member $\tilde{X}$ so to reduce the study
of the birational geometry of $X$ to that of $\tilde{X}$.
The first basic fact is therefore to define what it means to be minimal.
This is absolutely a non trivial problem and the following definition
is the result of hard work of persons like Mori and Reid in the
late 70's.

\begin{Definition} A variety $\tilde{X}$ is minimal if
\begin{itemize}
\item[-] $\tilde{X}$ has $\Q$-factorial terminal singularities
\item[-] $K_{\tilde{X}}$ is nef
\end{itemize}
\end{Definition}

Let us make some observations on this definition.
The second condition wants to express the fact that the
minimal variety is {\sl (semi) negatively curved}.
We note in fact that if $det T_X = -K_X$ admits a metric
with semi-negative curvature then $K_X$ is nef.
The converse is actually an open problem (true in the case of surfaces
and in general it may be considered as a conjecture).

The condition on the singularities is the real break-through of the definition.
The point of view should be the following: we are in principle
interested in smooth varieties but we will find out that
there are smooth varieties which does not admit smooth minimal models.
  However we can find such a model
if we admit very mild singularities, the ones stated in the definition.
Note also that terminal singularities are smooth in the surface case.

It happens that in the birational 
class of a given variety there is not a  minimal model,
think for instance of rational varieties.
But the MMP hopes to make a list of them all.

\smallskip
Given the definition of minimal variety we want now
to show how, starting from $X$, one can determine a
corresponding minimal model $\tilde{X}$.
In view of \ref{cone} and \ref{F-Mcontr} (or \ref{contR})
the way to do it is quite natural.
Namely, if $K_X$ is not nef, then by \ref{cone} there exists an extremal ray
(on which $K_X$ is negative) and by \ref{F-Mcontr} (or \ref{contR})
we can construct an elementary (Fano-Mori) contractions $f:X\ra X^{\prime}$
which contracts all curves in this ray
into a normal projective variety $X^{\prime}$.

A naive idea at this point would be the following.

If $f$ is of fiber type, i.e. $dimX' < dim X$, then one hopes to recover
a description of $X$ via $f$. Indeed, by induction on the dimension, one
should know a description of $X'$ and of the fibers of $f$, which 
are, at least generically, Fano varieties.
We will say something more on this case in the last part of the section
(see \ref{MoriSpaces}).

If $f$ is birational then one thinks to substitute
$X$ with $X^{\prime}$ and proceed inductively.

The problem is of course that the Theorem \ref{F-Mcontr}
says very little about the singularities of
$X^{\prime}$ (now it starts to be clear that the choice of the singularities
in the above definition is crucial). It says only that it has
normal singularities.

\smallskip
However in the surface case the situation is optimal,
namely Theorem \ref{2-folds} first of all says
that if $cont_R$ is birational then the image is again
a smooth surface (see \ref{2-folds}.1). Then apply recursively
\ref{2-folds}.1 and obtain that after finitely many blow downs of
$(-1)$-curves one reaches a smooth
surface $S^{\prime}$ with  either $K_{S^{\prime}}$ nef
or with an extremal ray of fiber type.
Note that while performing the MMP we stay in the category of smooth surfaces.
If $cont_R$ is of fiber type then, again by \ref{2-folds},  its 
description is very precise.
We have proved the following.

\begin{Theorem} [Minimal Model Program for surfaces] Let $S$ be a smooth
surface.
After finitely many blow downs of $(-1)$-curves one reaches a smooth
surface $S^{\prime}$ satisfying one of the following:
\begin{itemize}
\item[1)] $K_{S^{\prime}}$ is nef i.e. $S^{\prime}$ is a minimal model
\item[2)] $S^{\prime}$ is a ruled surface
\item[3)] $S^{\prime}\iso\Proj^2$.
\end{itemize}
\end{Theorem}

\smallskip
In higher dimensions the requirement on the singularity starts to play.
In particular we note that cases B3), B4) and B5) in the Theorem \ref{3-folds}
lead to singular 3-folds; the case B5) leads to a
2-Gorenstein singularity.
However all this singularities are terminal and $\Q$-factorial.
The fact that the Cone Theorem \ref{cone} and the Contraction Theorem
\ref{F-Mcontr}
hold in the more general case of variety with terminal singularities seems
to give some hopes.

Moreover the good property of birational contractions in the surface case
ascends in higher dimension to the fact that if an elementary F-M
contraction of
a smooth (or terminal $\Q$-factorial) variety is
divisorial then the target has at worst terminal $\Q$-factorial singularities,
see Exercise \ref{divisorial}.

But a very serious problem is now coming up. Namely if we consider
varieties with terminal singularities then they can have birational
F-M contractions which are not divisorial! This was first noticed by
P. Francia with a famous example, see for instance \cite[p.33-34]{CKM}.

Let us see why an elementary contraction
$f:X\ra Y$ of a variety $X$ with $\Q$-Gorenstein singularities and
with exceptional locus $E$ such that $cod E \geq 2$ gives problems.
Let $U=X\setminus E$ then $f_U:U\to f(U)=V$ is an isomorphism. In particular it
is clear that
$K_{X|U}\iso f_U^*(K_V)$. Let $M$ be
the extension of $f_U^*(K_V)$ to $X$. Then the codimension assumption
yields $M\iso K_X$. On the other hand
$-K_X$ is $f$-ample therefore  $M$ cannot be the pullback of any 
$\Q$-Cartier divisor on $Y$.
In other words $K_Y$ is not $\Q$-Cartier!

In particular on such a $Y$ even the definition of minimal model
does not make sense. Our naive solution came abruptly to a stop and new
solutions are needed.
The principal ideas are summarised in the following Flip conjectures, \cite{KMM},
which very roughly says that
instead of contracting the exceptional locus of this ``small'' rays
we have to make a codimension 2 surgery, called flip,
that replace the curve with another one which has a different normal sheaf.

This is the precise statement.

\begin{conjecture}[Flip Conjecture]\ 

Let $X$ be a terminal $\Q$-factorial
variety and assume that
there exists an extremal ray $\R^+[C]\subset \NE(X)$ with associated
elementary contraction $f:X\ra W$; assume also that $codim (Exc(f)) \geq 2$.
Then there exists a terminal $\Q$-factorial variety $X^{\scriptstyle +}$
and a map $f^+:X^+ \ra W$ such that
\begin{itemize}
\item[1)]  $K_{X^+}$ is $f^+$ ample,
\item[2)]  $Exc(f^*)$ has codimension at least two in $X^+$,
\item[3)] the following diagram is commutative
$$\diagram
X\drto_f\rrdashed^{\Phi}|>\tip&&X^+\dlto^{f^+}\\
&W&
\enddiagram
$$
That is $X$ is isomorphic to $X^+$ in codimension 1.
\end{itemize}
The rational map $\Phi$ 
is called {\sl flip}.
\end{conjecture}

The following Theorem is the breakthrough of Mori Theory for threefolds
which over-passed the flip problem.
The proof is very intricate and it is based on a careful classification
of all possible small contractions occurring on a terminal 3-fold.

\begin{Theorem}[\cite{3flip}] The Flip conjecture holds for threefolds.
\end{Theorem}

\begin{remark} After Mori's proof of the existence of flips, different proofs of
flip, even log-flip, for 3-folds were obtained mainly by Shokurov and Koll\'ar. 
The best account of them is
in \cite{U2}. Very recently still a new approach of Shokurov simplified greatly the 
3-fold proof and is very promising in higher dimensions.
\end{remark} 

Then, assuming the flip conjecture, one asks if this sort
of inductive procedure will come to an end, namely we
need a kind of termination for these birational modification.

If $f$ is a divisorial contraction then the Picard number
drops by one
so there cannot be an infinite number of those.

For flips there is not such a straightforward criterion
and so the following termination conjecture raises.

\begin{conjecture}[Termination Conjecture]\ 

 Let $X$ be a terminal $\Q$-factorial variety which
is not minimal. Then
after finitely many flips there is an extremal ray whose exceptional locus
is of codimension$\leq 1$.
\label{termination}
\end{conjecture}

\begin{Proposition} [\cite{KMM} ,Theorem 5.1.15]
The Termination conjecture is true for $n$-folds with $n \leq 4$.
\end{Proposition}

\bigskip
So that, assuming both flips and termination conjectures, after finitely many birational modifications
we reach either a Minimal Model or we encounter an elementary extremal contraction of fiber type.
For this we give a definition.

\begin{Definition} A Mori space is a terminal $\Q$-factorial Fano-Mori
contraction
$\pi:X\ra S$ such that $dim S< dim X$ and $rkPic(X/S)=1$.\label{MoriSpaces}
\end{Definition}

The goal in this case is, like in the surfaces,
to get a classification of Mori spaces.

In general the Mori space associated to a variety by the MMP is not
uniquely determined. This problem arise when two extremal rays
have not disjoint exceptional loci.
A very simple example. Let $T=E\times {\bf F}_1$, where
$E$ is a smooth curve of genus $g>0$.
Then there are two extremal rays, one of
divisorial type and the other of fiber type. In this case the
order in which the rays are contracted determines the
F-M space. In one case it is a $\Proj^1$-bundle in the other
a $\Proj^2$-bundle.

Let us make some observation on these spaces.
We note that the generic fiber is a variety with $\Q$-factorial
terminal singularities with anticanonical bundle ample. Thus no multiple
of the canonical bundle of $X$ has a section, that is the Kodaira
dimension of $X$ is $-\infty$. On the other hand we already noticed, see section \ref{sec:ratcur},
 that such an $X$ is covered by rational curves. There is a deep, and still not entirely understood,
 relation between these two facts. 

\begin{Definition} A variety $X$ is {\sf uniruled} if there exists
a generically finite surjective map $Y\times\Proj^1\flip X$.
\label{unir}
\end{Definition}

\begin{Proposition} [\cite{KMM}, Corollary 5.1.4] Let $\pi:X\ra S$ be 
a Mori Space, then
$X$ is a uniruled variety
\end{Proposition}

The proposition is a consequence of
the following Theorem proved by means of
the theory of deformation of rational curves
on smooth varieties.

\begin{Theorem}[\cite{MiMo}] Let $X$ be a projective variety.
  Assume that for
a general $x\in X$ there is a smooth proper
curve $C$ and a morphism $f:C\ra X$ such that
\begin{itemize}
\item[-] $x\in f(C)$
\item[-] $X$ is smooth along $f(C)$
\item[-] $deg_Cf^*K_X<0$
\end{itemize}
Then $X$ is uniruled.
\end{Theorem}

\medskip
To conclude the above discussion let us write a Minimal Model Conjecture.

\begin{conjecture} Let $X$ be a projective variety with at most terminal
$\Q$-factorial singularities.
Then there exists a minimal model $X'$ birational to $X$ if and only if
$X$ is not uniruled.
\end{conjecture}

In dimension 3 the above conjecture is now a Theorem. To prove it
we still need the following observations.

\begin{Proposition}[\cite{RaC}] Let $X$ be a smooth uniruled variety.
Then there is a dense family of rational curves with negative
intersection with the canonical class. In particular $X$ has negative
Kodaira dimension.
\label{uniconj}
\end{Proposition}

\begin{exercise} Let $X$ be a threefold with terminal singularities
covered by curves negative
with respect to the canonical class. Prove that,
after finitely
many birational modifications, instead of reaching
a minimal model one
encounters an extremal ray whose exceptional locus covers the whole
variety. In other words one reaches a Mori space.

This can be proved also for n-folds as soon as one assumes the flip 
conjecture and the
termination conjecture.
\label{mmpmos}
\end{exercise}

Summing things up we get.

\begin{Theorem} Let $X$ be a threefold with terminal singularities, then
the Minimal Model conjecture holds.
\label{th:3mmp}
\end{Theorem}

Let us also mention the following result.

\begin{remark} The converse of Theorem \ref{uniconj} is a deep and
challenging problem. It is conjectured that all smooth varieties with
negative Kodaira dimension are uniruled. But a positive answer is only known up
to dimension 3, as a byproduct of MMP, see \cite{Mi}.
\end{remark}

As a consequence one can formulate the Minimal Model Conjecture in 
this straightened form,
which is also true in dimension $3$.

\begin{conjecture} Let $X$ be a projective variety with at most terminal
singularities.
Then there is a minimal model $X'$ birational to $X$ if and only if $k(X)
\not= -\infty$.
\end{conjecture}


\section{The birational geometry of the plane}
\label{S_2}

Sarkisov program is devoted to study the possible birational,
not biregular, maps between Mori spaces.

We do not want here to outline the complete program and its applications,
for this we refer the reader to
\cite{Co1}, \cite{Co2}.
However we like to give an idea of its techniques and possible
applications in the  simpler set up of surfaces;
for this, using Sarkisov dictionary, we prove
the following beautiful Theorem.

\begin{Theorem}[Noether-Castelnuovo]\
The group of birational transformations of
the projective plane is generated by linear transformations and the 
standard Cremona
transformation, that is
$$(x_0:x_1:x_2)\mapsto(x_1x_2:x_0x_2:x_0x_1),$$
where $(x_0:x_1:x_2)$ are the coordinates of $\Proj^2$.
\label{th:NC}
\end{Theorem}

Let $\chi:\Proj^2\flip\Proj^2$ a birational map which is not an
isomorphism.
To study the map $\chi$ we start factorizing it with simpler birational maps, ``elementary links'',
between Mori Spaces (these maps will be either the blow up of a point
in $\Proj^2$, or an elementary transformation of a rational ruled surface, see diagram (\ref{elt})).
Consider 
 $\H=\chi_*^{-1}\O(1)$, the strict transform of lines in $\Proj^2$; then
$\H$ is without fixed components and  $\H\subset |\O(n)|$ for some $n>1$. Our point of view
is to consider
 the general element $H\in\H$ as a twisted line. 
The factorisation we are aiming  should ``untwist'' $H$  step by step
so to give back the original line hence the starting $\Proj^2$. 
 Observe that
the fact that $\chi$ is not biregular is encoded in the base locus of $\H$, therefore the untwisting
is clearly related to the singularities of the log pair
$(\Proj^2,\H)$, where by the pair
$(\Proj^2,\H)$ we understand the pair $(\Proj^2,H)$ were $H\in\H$ is 
a general element.

\begin{Theorem}Let $\H\subset |\O(n)|$ be as above;
then the pair $(\Proj^2,(3/n)\H)$
has not canonical singularities.

In particular there is a point $x\in \Proj^2$ such that
\begin{equation}
\label{mult_p2}
\mult_x\H>n/3.
\end{equation}
\label{th:NF2}
\end{Theorem}

\begin{proof}[Proof of Theorem \ref{th:NF2}]
Take a resolution of $\chi$
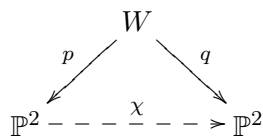
$$
\diagram
&W\dlto_p\drto^q&\\
\Proj^2\rrdashed|>\tip^{\chi}&&\Proj^2
\enddiagram
$$
and pull back the divisor $K_{\Proj^2}+(3/n)\H$ and $K_{\Proj^2}+(3/n)\O(1)$ 
via
$p$ and $q$ respectively.

We have
\begin{eqnarray*}
K_W+(3/n)\H_W &=&\\
p^*\O_{\Proj^2}+\sum_i a_i^{\prime}E_i+\sum_h c_hG_h&=&
q^*\O_{\Proj^2}(3(1/n-1))+\sum_i a_i E_i+\sum_j b_j F_j
\end{eqnarray*}

where $E_i$ are $p$ and $q$ exceptional divisors, while $F_j$ are $q$ 
but not $p$ exceptional divisors
and $G_h$ are $p$ but not $q$ exceptional divisors.
Observe that, since $\O(1)$ is base point free, the $a_i$'s and $b_j$'s are
positive integers.

Let $l\subset\Proj^2$ a general line in the right hand side plane. In 
particular $q$ is an isomorphism
on $l$ and therefore $E_i\cdot q^*l=F_j\cdot q^*l=0$ for all $i$ and $j$.

The crucial point is that on the
right hand side we have some negativity coming from the non effective divisor
$K_{\Proj^2}+(3/n)\O(1)$
that has to be compensated by some non effective exceptional divisor 
on the other side.

More precisely, since $n>1$, we have on one hand that
$$(K_W+(3/n)\H_W)\cdot q^*l=(q^*\O_{\Proj^2}(3(1/n-1))+\sum_i a_i 
E_i+\sum_j b_j F_j)\cdot q^*l<0,$$
and on the other hand that
$$0>(K_W+(3/n)\H_W)\cdot q^*l=(p^*\O_{\Proj^2}+\sum_i 
a_i^{\prime}E_i+\sum_h c_h G_h)\cdot q^*l.$$
So that
  $c_h<0$ for some $h$, that is $(\Proj^2,(3/n)\H)$ is not canonical.

We leave to the reader to justify equation (\ref{mult_p2}); remember that
one can resolve the base locus
of $\H$ blowing up smooth points only.
\end{proof}

The above proof can be generalised to the following set up. Let $\pi:X\ra
S$ and $\f:Y\ra W$
be two Mori spaces of dimension $\leq 3$. Let $\chi:X\flip Y$ a birational
not biregular map. Choose
$\H_Y$ a very ample linear system on $Y$. Let $\H=\chi_*^{-1}\H_Y$ then by
the definition of
Mori space there exists a $\mu\in \Q$ such that $K_X+(1/\mu)\H\nel{\pi}0$.

\begin{Theorem}[Noether--Fano inequalities]\cite{Co1} In the above notation,
in particular with $\chi$ non biregular and $K_X+(1/\mu)\H\nel{\pi}0$,
then either $(X,(1/\mu)\H)$ has not canonical singularities
or $K_X+(1/\mu) \H$ is not nef.
\label{NFineq}
\end{Theorem}

We are now ready to start the factorisation of $\chi$.

For this let $x\in \Proj^2$ be a point such that $(\Proj^2,(3/n)\H)$ is not
canonical at $x$.
Such a point exists by  Theorem \ref{th:NF2}
and let $\nu:\F_1\ra \Proj^2$ be the
blow up of $x$, with exceptional divisor
$C_0$.

In the context of Sarkisov theory it is important to look at this blow up in the following way.
\begin{Definition} A terminal extraction is a birational morphism with
connected fibers
$f:Y\supset E\ra X\ni x$.
Such that:
\begin{itemize}
\item[-] $X$ and $Y$ are terminal varieties, $Y$ is $\Q$-factorial
\item[-] the exceptional locus is an irreducible divisor
$E$, with $f(E)\ni x$
\item[-] $-K_Y$ is $f$-ample.
\end{itemize}
\end{Definition}
\begin{exercise}Prove that the only terminal extraction
from a smooth point of a surface is the blow up of the maximal ideal of a
point.

\hint For a surface terminal is equivalent to smooth. This is a just a restatement of 
Theorem \ref{2-folds}.
\end{exercise}
\begin{remark} 
More generally whenever a log pair $(X,(1/\mu)\H)$ is not canonical then there exists a terminal
extraction, see \cite{Co1}.
\end{remark}

Let us come back to the proof.
Observe that the natural map $\pi_1:\F_1\ra \Proj^1$ is a Mori space
structure.
The map $\nu:\F_1\ra \Proj^2$, the blow up of $\Proj^2$,
is the first elementary link we define. 

Let $\chi^{\prime}=\chi\circ\nu:\F_1\flip \Proj^2$ and
$\H^{\prime}=(\chi^{\prime})_*^{-1}\O(1)$. Let $n'=n-mult_x\H$, then
$$K_{\F_1}+(2/n^{\prime})\H^{\prime}\nel{\pi_1}0.$$

We are in the conditions to apply Theorem \ref{NFineq} to the pair 
$(\F_1,(2/n')\H^{\prime})$.
Let us first notice that
$K_{\F_1}+(2/n^{\prime})\H^{\prime}$ is nef.
In fact let $f\subset\F_1$ a generic fiber of the ruled structure. Then
$$K_{\F_1}+(2/n^{\prime})\H^{\prime}\cdot f=0,$$
by definition. On the other hand
\begin{eqnarray*}
(K_{\F_1}+(2/n^{\prime})\H^{\prime})\cdot C_0& =&-1+(2/n^{\prime})mult_{x}\H\\
&=& -n+3mult_{x}\H>0
\end{eqnarray*}
where the last inequality comes directly from equation (\ref{mult_p2}),
that is the existence of non canonical singularities for $(\Proj^2,(3/n)\H)$.

Thus, by Theorem \ref{NFineq},  $K_{\F_1}+(2/n^{\prime})\H^{\prime}$ is not canonical and
therefore the linear system $\H^{\prime}$ admits a point 
$x^{\prime}\in \F_1$
with  multiplicity greater than $2/n'$.

The next step is a terminal extraction from $x'$. Let 
$$\psi:Z\supset E\ra \F_1\ni x^{\prime},$$
the blow up of $x'$.

This time $Z$ is not a Mori space, but the strict transform of the
fiber of $\F_1$ containing
$x^{\prime}$ is now a
$(-1)$-curve which can then be contracted: $\f:Z\to S$.
\begin{equation}
\label{elt}
{\diagram
&Z\dlto_{\psi}\drto^{\f}&\\
\F_1\rrdashed|>\tip&&S
\enddiagram}
\end{equation}

This modification is known as an {\sl elementary transformation} of 
ruled surfaces.

\begin{exercise} Prove that $S$ is either a quadric, $\F_0$, or $\F_2$.
\end{exercise}
\noindent
\hint It depends on the
position of the point with respect to $C_0$.

\smallskip
Let $x_2\subset S$ be the exceptional locus of $\f^{-1}$ and $\H_2$
be the strict transform of $\H'$.
Observe the following two facts:
\begin{itemize}
\item[i)] $(K_{S}+(2/n^{\prime})\H_2)\cdot f=0$,
where, by abuse of notation, $f$ is the strict transform of
$f\subset\F_1$,
\item[ii)] since $mult_{x^{\prime}}\H'>\frac{\H'\cdot f}2,$
then $(S,(2/n^{\prime})\H_2)$ has terminal singularities at $x_2$.
\end{itemize}

By i)
we can apply Theorem \ref{NFineq} to  the log pair $(S,(2/n^{\prime})\H_2)$. 
Moreover by ii) we did not introduce any new canonical singularities
since the point $x_2$ is a terminal singularity for this pair. 
This is very important because proves that after finitely many elementary
transformations we reach a pair $(\F_k,(2/n^{\prime})\H_r)$ with
canonical singularities
  such that
$$K_{\F_k}+(2/n')\H_r\nel{\pi_k}0.$$
Then, again by Theorem  \ref{NFineq}, the pair $(\F_k,(2/n^{\prime})\H_r)$ cannot be nef.

Observe that $NE(\F_k)$ is a two dimensional cone. In particular it 
has only two
rays. One is spanned by $f$, a fiber of $\pi_k$. Let $Z$ an effective 
irreducible curve in the other ray.
Then
\begin{equation}
\label{eq:neg}
(K_{\F_k}+(2/n')\H_r)\cdot Z<0.
\end{equation}
Since $\H_r$ has not fixed components then $\F_k$ is a del Pezzo surface and
the only possibilities are therefore $k=0,1$.

In case $k=1$ then what is left is to simply blow
down the exceptional curve $\nu:\F_1\to\Proj^2$, and reach $\Proj^2$ 
together with a linear
system $\nu_*\H_2=:\tilde{\H}\subset |\O(j)|$.  Note that in this 
case, by equation (\ref{eq:neg}),
$$K_{\F_1}+(2/n')\H_r=\nu^*(K_{\Proj^2}+(2/n')\tilde{\H})+\delta C_0,$$
for some positive $\delta$. Therefore $K_{\Proj^2}+(2/n')\tilde{\H}$ 
is not nef.
In other terms
$$(2/n')j<3,$$
and
$$j<\frac{3(n-mult_x\H)}2<n.$$

Now we iterate the above argument,
i.e. we restart at the beginning of the proof but with the pair 
$(\Proj^2,(3/j)\tilde{\H})$;
the above strict inequality $j < n$, tells us that
and after finitely many steps we
untwist the map $\chi$, i.e. we reach $\Proj^2$ with a linear system 
$\H= |\O(1)|$.

In case $k=0$ observe that $\F_0\iso \Q^2$ is a Mori space for two 
different fibrations, let
$f_0$ and $f_1$ the general fibers of these two fibrations.
Moreover by equation (\ref{eq:neg})
$$(K_{\F_0}+(2/n')\H_r)\cdot f_1<0.$$
That is there exists an
\begin{equation}
\label{eq:n}n_1<n'
\end{equation}
  such that
$$(K_{\F_0}+(2/n_1)\H_r)\cdot f_0>0,$$
and
$$(K_{\F_0}+(2/n_1)\H_r)\cdot f_1=0.$$
Again by Theorem \ref{NFineq}, this time applied with respect to the fibration
with fiber $f_1$, this implies that
$(\F_0,(2/n_1)\H_r)$ is not canonical and we iterate the procedure.
As in the previous case the strict inequality of equation (\ref{eq:n})
implies a termination after finitely many steps.

Thus we have factorised any birational, not biregular, self-map of 
$\Proj^2$ with
a sequence of ``elementary links'', namely elementary transformations and blow ups of $\P^2$ at a
point.

The next step  is to interpret a standard Cremona transformation in this
new language, i.e. in term of the elementary links we have introduced above.

\begin{exercise} Prove that a standard Cremona transformation is given by the
following links
$$
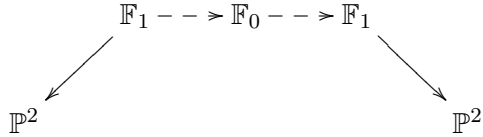
\diagram
&\F_1\dlto\rdashed|>\tip&\F_0\rdashed|>\tip&\F_1\drto\\
\Proj^2&&&&\Proj^2
\enddiagram
$$
Vice-versa any map of type
$$
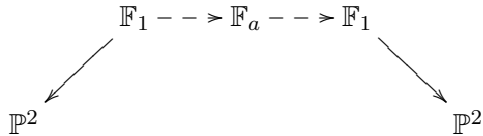
\diagram
&\F_1\dlto\rdashed|>\tip&\F_a\rdashed|>\tip&\F_1\drto\\
\Proj^2&&&&\Proj^2
\enddiagram
$$
can be factorised by Cremona
transformations.
\label{ex:cre}

\noindent
\hint A standard Cremona transformation is given by conics trough 3 
non collinear
points. The link above is possible only for $a=0,2$. They represents 
birational maps
given by conics with either 3 base points or 2 base point plus a 
tangent direction.
Try to factorise
the following map
$$(x_0:x_1:x_2)\to (x_1x_2:x_0x_2:x_1x_2+x_0x_2+x_0^2).$$
with Cremona transformations.
\end{exercise}

\begin{proof}[Proof of Theorem \ref{th:NC}]
Let $\chi:\Proj^2\flip\Proj^2$ a birational map and
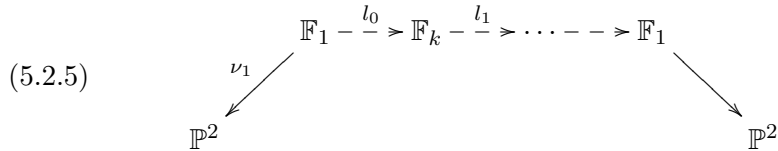
\begin{equation}
\label{eq:diag}
{\diagram
&\F_1\dlto_{\nu_1}\rdashed^{l_0}|>\tip&\F_k\rdashed^{l_1}|>\tip&\ldots\rdashed|>
\tip&\F_1\drto&\\
\Proj^2&&&&&\Proj^2
\enddiagram}
\end{equation}
the factorisation in elementary links obtained above.
Let us first make the following observation.
If there is a link leading to an $\F_1$ then we can break the birational 
map simply blowing
down the $(-1)$-curve. That is substitute $\chi$ with the following two pieces
$$
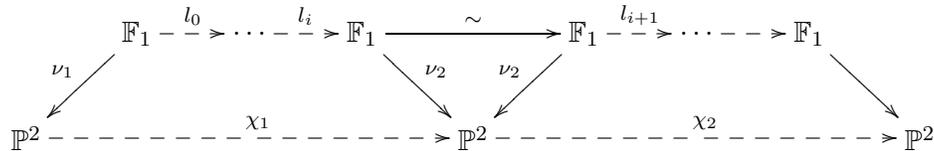
\diagram
&\F_1\dlto_{\nu_1}\rdashed^{l_0}|>\tip&\ldots\rdashed^{l_i}|>\tip&
\F_1\drto^{\nu_2}\rrto^{\sim}&&\F_1\dlto_{\nu_2}\rdashed^{l_{i+1}}|>\tip&\ldots\rdashed|>\tip&\F
_1\drto\\
\Proj^2\rrdashed&&\rrdashed^<{\chi_1}|>\tip&&\Proj^2\rrdashed&&\rrdashed^<{\chi_2}|>\tip&&\Proj^2
\enddiagram
$$
So that we can assume
\begin{equation}\label{eq:spezza}
\mbox{\rm  there are no links leading to $\F_1$ ``inside'' the factorisation.}
\end{equation}

Let
$$d(\chi)=\max\{k:\mbox{\rm there is an $F_k$ in the factorisation}\}.$$
If $d(\chi)\leq 2$ we can factorise it by exercise \ref{ex:cre}.

We now prove the Theorem by induction on $d(\chi)$.
Consider the left part of the factorisation (\ref{eq:diag}). Since 
$d(\chi)>2$, by assumption
(\ref{eq:spezza}), then
$l_0$ is of type $\F_1\flip \F_2$ and $l_1$ is of type $\F_2\flip \F_3$.
Then we force Cremona like diagrams in it,
at the cost of introducing new singularities.
Let
$$
\diagram
&\F_1\dlto_{\nu_1}\rdashed^{\alpha}|>\tip&\F_0\rdashed^{l_0}|>\tip&
\F_1\drto^{\nu_2}&&\F_1\dlto_{\nu_2}\rdashed^{\alpha^{-1}}|>\tip&\F_0\rdashed^{l_1}|>\tip&\ldots\F
_1\drto\\
\Proj^2&&&&\Proj^2\rrdashed&&\rrdashed^<{\chi'}|>\tip&&\Proj^2
\enddiagram
$$
where $\alpha:\F_1\flip\F_0$ is an elementary transformation centered 
at a general point of $\F_1$,
and $exc(\alpha^{-1})=\{y_0\}$.
So that $\alpha_{*}(\H')$ has an ordinary singularity at $y_0$.
Then $l_0$ is exactly the same modification but leads to an $\F_1$ and 
$\nu_2$ is the
blow down of the exceptional curve of this $\F_1$. Observe that 
neither $\alpha_0$ nor $\nu_2$
are links in the Sarkisov category, in general.
  Nonetheless the first part can be factorised by standard Cremona 
transformations. Let $\chi'=\chi_1\circ\ldots\circ\chi_k$ a decomposition
of $\chi$ in pieces satisfying (\ref{eq:spezza}). Then
$d(\chi_i)<d(\chi)$ for all $i=1,\ldots,k$. 
Therefore by induction hypothesis also $\chi'$ 
can be factorised by Cremona
transformations. Hence $\chi$ is factorised by Cremona transformations
\end{proof}

\section{\#-Minimal Model}
\label{sharp}

We already pointed out that Minimal Model Program allows to attach a Mori
space to
a uniruled 3-fold (see \ref{mmpmos}).
How can we use it to study the birational
geometry of $X$ ?

The main difficulty here is that the birational modifications
occurring along the MMP are difficult to follow and usually
it is almost impossible to guess what is the output.
We want to rephrase, after \cite{R2},
the standard minimal model program for
uniruled varieties
using a polarizing divisor; this is called \#-minimal model.
Under strong assumptions on
the variety studied, we are able to govern the program and understand 
its output.

\begin{Definition}[\cite{sdt}]
Let $T$ be a terminal $\Q$-factorial uniruled
3-fold and $\h$ a  movable linear system,
i.e. $\dim|mH|>0$ for $m\gg0$,
with generic element $H\in \h$ on $T$.
Assume that $H$ is nef, then
   $$\rho=\rho_{\h}=\rho(T,\h)=:\mbox{  \rm sup   }\{m\in \Q|H+mK_T
\mbox{
   \rm is an effective $\Q$-divisor  }\},$$
is the threshold of the pair $(T,\h)$.

Since we are assuming that $dim\h\geq 0$
then $\rho\geq 0$.
  A pair $(\T,\hh)$ is called {\sf \#-minimal model} of $(T,\h)$ if:
\begin{itemize}
\item[i)]$\T$ has a Mori fiber space structure $\pi:\T\ra W$ and
$\hh$ is a movable Weil divisor,
\item[ii)] there exists a birational map $\psi:T\dasharrow\T$ such 
that $\hh=\psi_*\h$,
\item[iii)] let $\Hs\in \hh$ a general member, then
$\rho(T,\h)K_{\T}+\Hs\nel{\pi} \Os_{\T}.$
\end{itemize}
\end{Definition}

To find a {\sf \#-minimal model} of a given pair $(T,\h)$
let us proceed in the following way.

Let  $(T_0,\h_0)=(T,\h)$, $H_0$ is nef by hypothesis
and $T_0$ is uniruled;  therefore to
$(T_0,H_0)$ is naturally associated the nef value
$t_0=sup\{ m\in \Q | mK_{T_0}+H_0$ is nef $\}$
and a rational map $\f_0:T_0\flip T_1$, which is either an extremal contraction
or (if the extremal contraction is small) a flip,
of an extremal ray in the face spanned by $t_0K_{T_0}+\h_0$
(see section (3.1) and in particular \ref{renef}).

Consequently on $T_1$ one defines a movable linear system
by $\h_1:=\f_{0*}\h_0$. That is to say
$\f_{0*}H_0\neq 0$.

Note that, by construction, $t_{0}K_{T_{1}}+H_{1}$ is nef,
thus one inductively defines
$\f_i:T_i\flip T_{i+1}$ and $(T_{i+1},\h_{i+1})$ as follows.

Let  $\delta=sup\{ d\in \Q
| dK_{T_i}+(t_{i-1}K_{T_i}+H_i)$ is nef $\}$ and define
$t_{i}:=\delta+t_{i-1}$.

\begin{exercise} Prove that there always exists an extremal ray
$[C_i]\subset \NE(T_i)$ in the face supported by $t_iK_{T_i}+H_i$.
\end{exercise}

Thus let us define
  $\f_i:T_i\flip T_{i+1}$ the birational modification associated to the extremal
ray $[C_i]\subset \NE(T_i)$, and
  $\h_{i+1}:=\f_{i*}\h_i$.

The inductive process is therefore composed by divisorial contractions
and flips. Since $T_0$ is uniruled it has not a minimal model, see 
Theorem \ref{th:3mmp}.
After finitely many of these birational modifications, we get
a Mori fiber space.

\begin{exercise} Prove that the output $(T_k,\h_k)$
is a \#-minimal model, that is
$$\rho(T,\h)K_{T_k}+H_k\nel{\pi} \Os_{T_k}.$$
\end{exercise}

\begin{remark}Note that $\hh$ is relatively nef. Furthermore if
the rational map defined by $|mH|$ is birational
then $\hh$ is relatively ample.
\end{remark}

The presence of a
polarization in the  \#-program  allows to control the steps if we are able
to impose restrictions on the
threshold.
Let $(T,\h)$ as above and assume moreover that $\rho_{\h}<1$ and
that there exists a smooth
surface $S\in \h$. Notice that the latter hypothesis
is not as strong as it seems, see
Definition \ref{def:lmr}. Under this assumption it
is
possible to describe in detail the \#-process in a neighborhood of
the surface $S$, see also \cite{CF}.

\begin{Proposition}[\cite{sdt}]
Let $\f_i:T_i\flip T_{i+1}$ be a birational modification
  in the \#-program
relative to $(T,\h)$ with  $\rho_{\h}<1$.
Assume that $S\in \h_i$ is  a smooth surface.
Then $\f_i(S)=\overline{S}$ is a smooth surface and $\f_{i|S}:S\ra
\overline{S}$ is either an isomorphism or the contraction of a disjoint
union of (-1)-curves.
\label{benino}
\end{Proposition}
\begin{proof}[Sketch of proof]
  Since $S$ is smooth and $T_i$ is terminal $\Q$-factorial
then $S$ is out of $Sing(T_i)$. In
particular $H_i$ is a Cartier divisor. We have the following cases.
\begin{itemize}
\item[-][$\f_i$ contracts a divisor $E$ onto a curve]
  Then $H_i\nel{\f_i}0$ and $S\cap E$ is the disjoint union
of (-1)-curves.
\item[-][$\f_i$ is a flip]
  $S$ is disjoint from the flipping curve
\item[-][$\f_i$ contracts a divisor $E$ to a point]
  $\f_{i|S}$ is birational and is either an isomorphism or the contraction of a
(-1)-curve.
$(E,E_{|E})\iso(\Proj^2,\O(-1))$ and $H_{|E}\sim \O(1)$.
\end{itemize}
\end{proof}

Using the above Proposition we can control the \# minimal model and its output.

\begin{Corollary}  Let $T$ be a terminal $\Q$-factorial
  uniruled 3-fold and $\h$ a movable and nef
linear system and  $(\T,\hh)$ a \#-minimal model of $(T,\h)$.
Assume that
$\rho_{\h}<1$ and $\h$ is base point free then
\hbox{$\Hs\in Pic(\T)$}, $\hh$ has at most base points and
$\Hs$ is smooth.
\label{bene}
\end{Corollary}
\begin{proof} By Bertini Theorem $H$ is smooth therefore we can apply
  Proposition \ref{benino} in an inductive way up to reach a
model $(\T,\hh)$.
\end{proof}

We need a relative version of Corollary \ref{bene},
and for this we first give a definition.

\begin{Definition}[\cite{sdt}] \label{def:lmr}Let $T$ be a 3-fold and $\h$  a 
movable linear system,
with $dim\h\geq 1$.
Assume that $H=M+F$, where ${\mathcal M}$ is a movable
linear system without fixed component and
$F$ is the fixed component.
  A pair $(T_1,\h_1)$ is called a log minimal
resolution of the pair $(T,\h)$ if
  there is a
morphism $\mu:T_1\ra T$, with the following
properties:
  \begin{itemize}
\item[-] $T_1$ is terminal $\Q$-factorial
\item[-] $\mu^{-1}_*M=H_1$, where
$H_1$ is a  Cartier divisor, $dimBsl(\h_1)\leq 0$
\item[-] a general element $H_{1}\in \h_1$ is a minimal
resolution of a
general element $M\in{\mathcal M}$.
\end{itemize}
\end{Definition}

\begin{Corollary} For any pair $(T,\h)$ with $T$ an irreducible
$\Q$-factorial 3-fold
and  $\h$ a movable linear system with $dim\h\geq 1$,
  there exists a log minimal
resolution. \label{lmr}
\end{Corollary}

\begin{remark}Using Corollary \ref{lmr} we can study any irreducible 3-fold  $T$ equipped
with a movable linear system $\h$ with $dim\h\geq 1$.
Indeed we consider a log minimal resolution of $(\T,\h)$ and then a  \#-Minimal Model of it.
Note that this is well defined only up to birational
equivalence.
\end{remark}

\section{Applications of the \#-program}

We want now to apply the \# theory to some concrete situations.
Despite the assumptions in the previous section,
are  quite strong, they are geometric in nature
and therefore of easy interpretation.

\subsection{3-folds with a uniruled movable system}

\begin{Definition} Let $T$ be a terminal $\Q$-factorial 3-fold and 
$\h$ a movable linear system.
We say that
$(T,\h)$ is a pair with a big uniruled system if
$H\in \h$ is nef and big and $H$ is a smooth surface of negative 
Kodaira dimension.
\end{Definition}

\begin{exercise} Let $(T,\h)$ be a pair with a big uniruled system. Then
$T$ is uniruled and $\rho(T,\h)<1$.
\end{exercise}

Using \#-MMP techniques we can
describe in detail
the possibilities that occur under these conditions.

\begin{Theorem}[\cite{sdt}]\label{esse}
  Let $(T,\h)$ be a pair with a big uniruled system. Then $(\T,\hh)$ 
is one of the following:
   \begin{itemize}

   \item[i)] a $\Q$-Fano 3-fold of index $1/\rho>1$,
with $K_{\T}\sim-1/\rho\Hs$ and
$\Phi_{|\Hs|}$ birational,
the complete classification is given in \cite{CF} and \cite{Sa}:
\begin{itemize}
\item[] $(\Proj(1,1,2,3),\Os(6))$
\item[] $(X_6\subset \Proj(1,1,2,3,a), X_6\cap \{x_4=0\})$, with $3\leq a
\leq 5$
\item[] $(X_6\subset\Proj(1,1,2,2,3),X_6\cap\{x_3=0\})$
\item[] $(X_6\subset\Proj(1,1,1,2,3),X_6\cap\{x_0=0\})$
\item[] $(\Proj(1,1,1,2),\Os(4))$
\item[] $(X_4\subset\Proj(1,1,1,1,2),X_4\cap\{x_0=0\})$
\item[] $(X_4\subset \Proj(1,1,1,2,a), X_4\cap \{x_4=0\})$, with $2\leq a
\leq 3$
\item[] $(\Proj^3,\Os(a))$, with $a\leq 3$, $(\Q^3,\Os(b))$, with $b\leq 2$
\item[] $(X_3\subset \Proj(1,1,1,1,2),X_3\cap \{x_4=0\})$,
$(X_3\subset \Proj^4,\Os(1))$
\item[] $(X_{2,2}\subset \Proj^5,\Os(1))$
\item[] a linear section of the Grassmann variety
parametrising
   lines in $\Proj^4$, embedded in $\Proj^9$ by Plucker
coordinates
\item[] $(\Proj(1,1,1,2),\Os(2))$, the cone over the Veronese surface
\end{itemize}
\item[ii)] a bundle over a smooth curve with at most finitely many fibers
$(G,\Hs_{|G})\iso
(\Sn_4,\Os(1))$, and  generic fiber $(F,\Hs_{|F})
\iso(\Proj^2,\Os(2))$. Where $\Sn_4$ is the cone over the normal 
quartic curve and
the vertex sits over an hyper-quotient singularity of type $1/2(1,-1,1)$ with
$f=xy-z^2+t^k$, for $k\geq 1$, \cite{YPG},
\item[iii)]
   a quadric bundle with at most
$cA_1$ singularities of type $f=x^2+y^2+z^2+t^k$, for $k\geq 2$,
and $\Hs_{|F}\sim \Os(1)$,
\item[iv)]$(\Proj(E),\Os(1))$ where $E$ is a rk 3  vector bundle over a
smooth curve,
   \item[v)] $(\Proj(E),\Os(1))$ where $E$ is a rk 2  vector
bundle over a surface of negative Kodaira dimension.
   \end{itemize}
\end{Theorem}

\begin{remark} The above Theorem allows to extend the result of \cite{CF}
  to
  3-folds $T$ which contain a smooth surface
$H$ of negative
Kodaira dimension such that $H$ is nef and big.  Ciro Ciliberto
pointed out to us that the Theorem completes 
a research suggested by
Castelnuovo, \cite[pg 187]{Cas}, to study linear systems of rational surfaces.
\end{remark}

\begin{exercise}Prove the following.
  Let $T$ be a terminal 3-fold and $H\subset T$
a smooth surface
of negative Kodaira dimension. Assume that $H$ is nef and big,
then $T$ is birational to one of the following:
\begin{itemize}
\item[-] $\Proj^3$
\item[-] $H\times \Proj^1$,
\item[-] a terminal sextic in either $\Proj(1,1,1,2,3)$ or $\Proj(1,1,2,2,3)$,
\item[-] a terminal quartic in $\Proj(1,1,1,1,2)$,
\item[-] a terminal cubic in $\Proj^4$.
\end{itemize}
\end{exercise}

There exists a natural geometric interpretation of the conditions
imposed in Theorem \ref{esse}.

\begin{Theorem}[\cite{sdt}]
   \label{main}
   Let $T_d\subset \Proj^n$ be a degree $d$ non degenerate irreducible 
3-fold. Suppose that
   $d<2n-4$, then any \#-Minimal Model $(\T,\hh)$ of $(T_d,\Os(1))$
is in the list of Theorem \ref{esse}.
\end{Theorem}

\begin{proof}
Let $\nu:X\ra T$ a resolution of singularities and $\h=\nu^*\Os(1)$.
First prove that  $(K_X+H)\cdot H^2<0$.
We argue comparing Castelnuovo bound on the genus of $C:=H^2$
  and the genus formula
on the surface $H$. From the latter we obtain that
$g(C)=1+d/2+(K_X+H)\cdot C/2$. For the former let $m=\lfloor
\frac{d-1}{n-3}\rfloor$ then by Castelnuovo inequality,
\cite[pg 527]{GH}, we have that
$$g(C)\leq \frac{m(m-1)}{2}(n-3)+m(d-1-m(n-3)).$$
It is therefore enough to impose that
\begin{equation}
1+d/2>\frac{m(m-1)}{2}(n-3)+m(d-1-m(n-3)),
\label{cab}
\end{equation}
after a small calculation one verifies that this is true
whenever $d<2n-4$.
Then
by adjunction formula,
$H\in\h$ is a smooth surface of negative Kodaira
dimension and Theorem \ref{esse} applies.
\end{proof}

\begin{remark}This Theorem can be interpreted as the
three dimensional counterpart of the classical result that a
non degenerate surface $S\subset \Proj^n$ of degree $d\leq n-1$ is
birational either to a rational scroll or to a projective plane, 
\cite[pg 525]{GH}.
Observe that all the listed 3-folds admit an embedding
satisfying the numerical criterion.
\end{remark}

By means of adjunction theory on terminal varieties, see \cite{Meadj}, one
can prove the following
higher dimensional
analog of Theorem \ref{main}.

\begin{Theorem}[\cite{sdt}] Let $X_d\subset\Proj^n$ a non degenerate
$k$-fold with $k>3$ and only $\Q$-factorial terminal singularities.
Assume that $d<2(n-k)-2$ then a \#-minimal model $(\X,\Hs)$ of
$(X,\Os(1))$, in adjunction theory language $(\X,\Hs)$ is the first reduction,
is one of the following:
\begin{itemize}
   \item[i)] a $\Q$-Fano n-fold of Fano index $1/\rho>k-2$,
with $K_{\T}\sim-1/\rho\Hs$ and
$\Phi_{|\Hs|}$ birational,
the complete classification is given in \cite{Fu} if $\X$ is Gorenstein
  and in \cite{CF} and \cite{Sa} in the non-Gorenstein case.
\item[ii)] a
   projective bundle over a curve with fibers $(F,\Hs_{|F})
\iso(\Proj^{k-1},\Os(1))$,  or a quadric bundle with at most
$cA_1$ singularities, with $\Hs_{|F}\sim \Os(1)$;
   \item[iii)] $(\Proj(E),\Os(1))$ where $E$ is a $rk(k-1)$ ample vector
bundle either
  on $\Proj^2$ or on a ruled surface.
   \end{itemize}
\end{Theorem}

\begin{remark}Note that since $X_d$ has terminal singularities
then, assuming the minimal model conjecture,
the above is the classification
of \#-models of those varieties.
\end{remark}

\subsection{General elephants of $\Q$-Fano threefolds}

Another direction is the study of the birational class of $\Q$-Fano whose
generic section has worse than canonical singularities.

\begin{conjecture}[Reid] Let $X$ be a $\Q$-Fano threefold and $H\in|-K_X|$
a generic section of the anticanonical divisor. Then $H$ has at worst
canonical singularities.
\end{conjecture}

The motivation of this conjecture is that to classify Fano variety,
as we learned, one uses sections of the fundamental divisor. For
non Gorenstein 3-fold with index$<1$,
this invariant is quite meaningless and one tries
to use directly sections of $|-K_X|$. So more than a conjecture it is
an hope that things are not too bad in this corner of the world.
It has to be said that the most recent techniques to study $\Q$-Fano do not
rely completely on this procedure.
The \#-program allows to understand the birational nature of these
strange objects.

\begin{Theorem}[\cite{sdt}] Let $T$  be a $\Q$-Fano 3-fold.

\noindent Assume that
$dim\phi_{|-K_T|}(T)=3$ and the general element in $|-K_T|$ has
worse than Du Val singularities. Then $T$ is
birational to a smooth Fano 3-fold $\T$ of index $\geq 2$.
\label{mainFano}
\end{Theorem}

The rough idea  is to take a log minimal resolution of $(T,|-K_T|)$
and control the output. By the singularity requirement the generic element in
$|-K_T|$ is a uniruled surface therefore we can apply all results
of previous sections. For more details see \cite{sdt}.

\end{document}